\documentclass[journal]{new-aiaa}
\usepackage[utf8]{inputenc}
\usepackage{textcomp}
\usepackage{chngcntr}
\usepackage{graphicx}
\usepackage{subcaption}
\usepackage{graphicx}
\usepackage{floatpag}
\usepackage{float}
\usepackage{afterpage}
\usepackage{amsmath}
\usepackage[version=4]{mhchem}
\usepackage{siunitx}
\usepackage{longtable,tabularx}
\usepackage{subcaption}
\usepackage{nomencl}
\makenomenclature
\usepackage{etoolbox}
\renewcommand\nomgroup[1]{%
  \item[\bfseries
  \ifstrequal{#1}{A}{Sets, Indices, and Graph}{%
  \ifstrequal{#1}{B}{Parameters}{%
  \ifstrequal{#1}{C}{Decision Variables}{}}}%
]}
\usepackage{algorithm}
\usepackage{algpseudocode}
\usepackage{geometry}
\usepackage{atbegshi}
\usepackage{etoolbox}
\usepackage{caption}
\usepackage{soul}
\usepackage{booktabs} 
\usepackage{array} 
\usepackage{caption}
\usepackage{longtable}

\usepackage{ifthen}

\usepackage{nomencl}
\usepackage{xpatch}

\newlength{\nomitemorigsep}
\title{An Integrated Supply Chain Network Design for Advanced Air Mobility Aircraft Manufacturing Using Stochastic Optimization}

\author{Esrat Farhana Dulia \footnote{Graduate Research Assistant, College of Aeronautics and Engineering, edulia@kent.edu.} and Syed A. M. Shihab \footnote{Assistant Professor, College of Aeronautics and Engineering, sshihab@kent.edu.}}
\affil{Kent State University, Kent, OH, USA  44242}

\begin{document}
\newgeometry{left=1.65cm,right=1.65cm,top=1.65cm,bottom=1.65cm}
\maketitle

\begin{abstract}
\begin{spacing}{1}







{
Electric vertical takeoff and landing (eVTOL) aircraft manufacturers await numerous pre-orders for eVTOLs and expect demand for such advanced air mobility (AAM) aircraft to rise dramatically soon. However, eVTOL manufacturers (EMs) cannot commence mass production of commercial eVTOLs due to a lack of supply chain planning for eVTOL manufacturing. The eVTOL supply chain differs from traditional ones due to stringent quality standards and limited suppliers for eVTOL parts, shortages in skilled labor and machinery, and contract renegotiations with major aerospace suppliers. The emerging AAM aircraft market introduces uncertainties in supplier pricing and capacities, eVTOL manufacturing costs, and eVTOL demand, further compounding the supply chain planning challenges for EMs. Despite this critical need, no study has been conducted to develop a comprehensive supply chain planning model for EMs. To address this research gap, we propose a stochastic optimization model for integrated supply chain planning of EMs while maximizing their operating profits under the abovementioned uncertainties. We conduct various numerical cases to analyze the impact of 1) endogenous eVTOL demand influenced by the quality of eVTOLs, 2) supply chain disruptions caused by geopolitical conflicts and resource scarcity, and 3) high-volume eVTOL demand similar to that experienced by automotive manufacturers, on EM supply chain planning. The results indicate that our proposed model is adaptable in all cases and outperforms established benchmark stochastic models. The findings suggest that EMs can commence mass eVTOL production with our model, enabling them to make optimal decisions and profits even under potential disruptions.

}

\end{spacing}
\end{abstract}

\begin{spacing}{1.35}

\section*{Keywords}
Supply chain network design;  stochastic optimization;  advanced air mobility;  electric vertical takeoff and landing aircraft;  rolling horizon approach.


\section{Introduction} \label{intro}


\subsection{Emerging eVTOL Industry}



\begin{figure}[tb]
\centering
\begin{minipage}{0.5\textwidth}
  \centering
  \includegraphics[width=\textwidth]{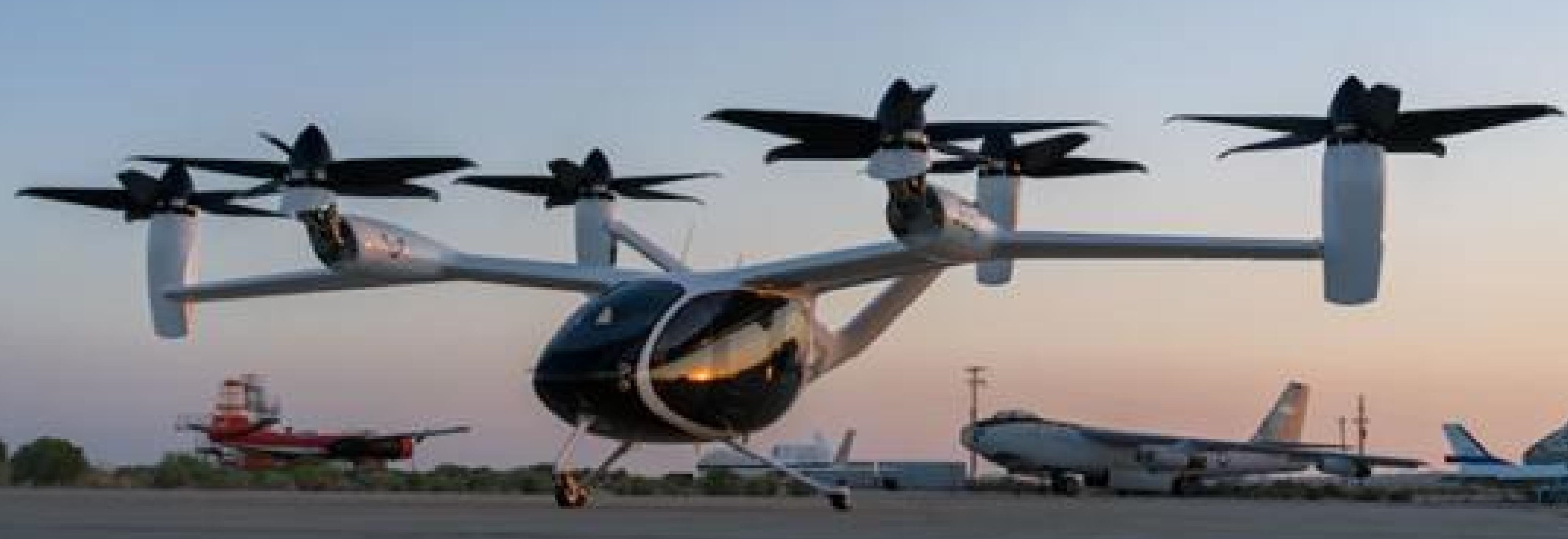}
  \caption{Joby Aviation eVTOL aircraft.}
  \label{fig_evtol1}
\end{minipage}
\hfill
\begin{minipage}{0.45\textwidth}
  \centering
  \includegraphics[width=\textwidth]{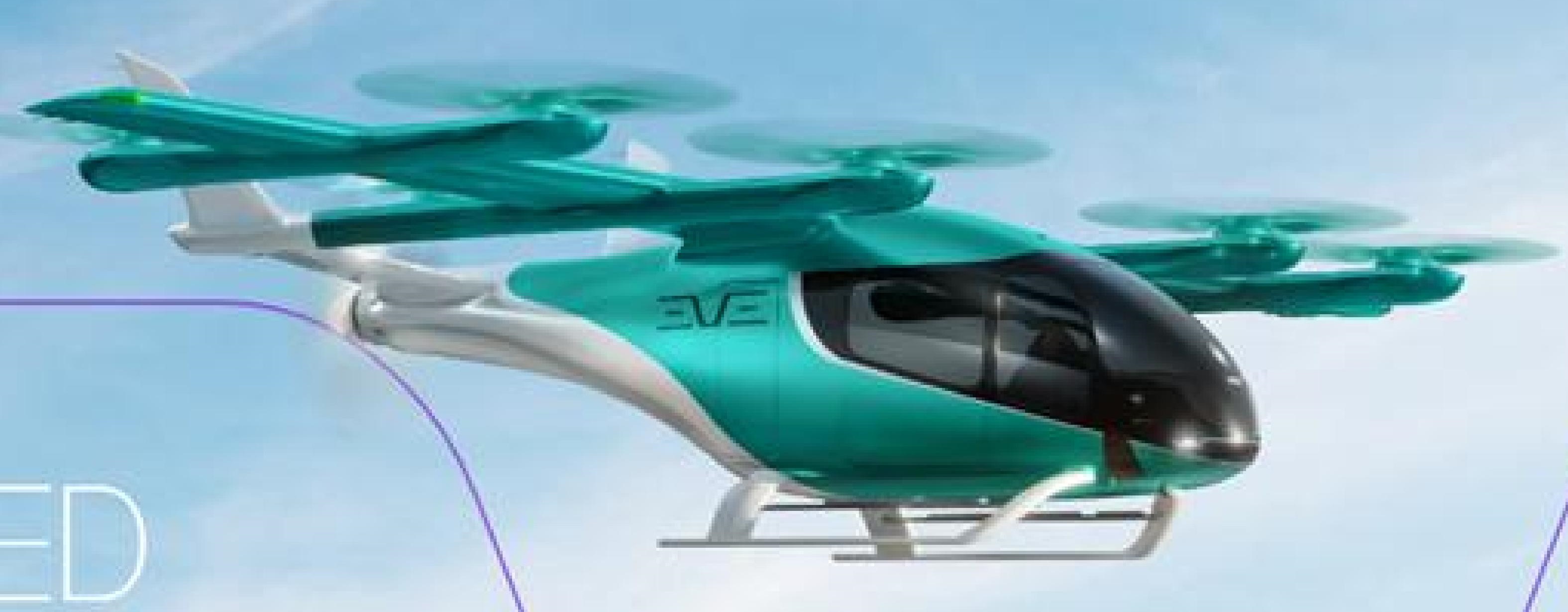}  
  \caption{Eve Air Mobility eVTOL aircraft.}
  \label{fig_evtol2}
\end{minipage}
\caption{Different models of eVTOLs.} 
\label{fig_evtol}
\end{figure}

Advanced Air Mobility (AAM) is a rapidly developing transportation system that aims to enable safe and efficient operations of novel and emerging aircraft, such as electric vertical takeoff and landing aircraft (eVTOL), in low-altitude airspace \cite{FAAConOpsv2}. AAM envisions a future where eVTOLs are utilized for diverse applications, including passenger transportation and cargo delivery, emergency medical services, search and rescue missions, aerial surveys, and infrastructure inspection \cite{dulia2021benefits}. These aircraft promise to be a faster, safer, and more sustainable means of  transportation \cite{garrow2021urban, vieira2019electric}. Various eVTOL models have been developed so far by different manufacturers, which include Joby Aviation, Airbus, Lilium, Volocopter, Beta Technologies, Archer Aviation, Eve Air Mobility, Vertical Aerospace, and many more. As of August 2022, there are 347 eVTOL companies in the market \cite{VFS_eVTOL_companies}. Two examples of such eVTOLs are shown in Figure \ref{fig_evtol}, sourced from \cite{joby_design, eve}, respectively. 

\begin{figure}[htb] 
    \centering
    \includegraphics[width=0.45\linewidth]{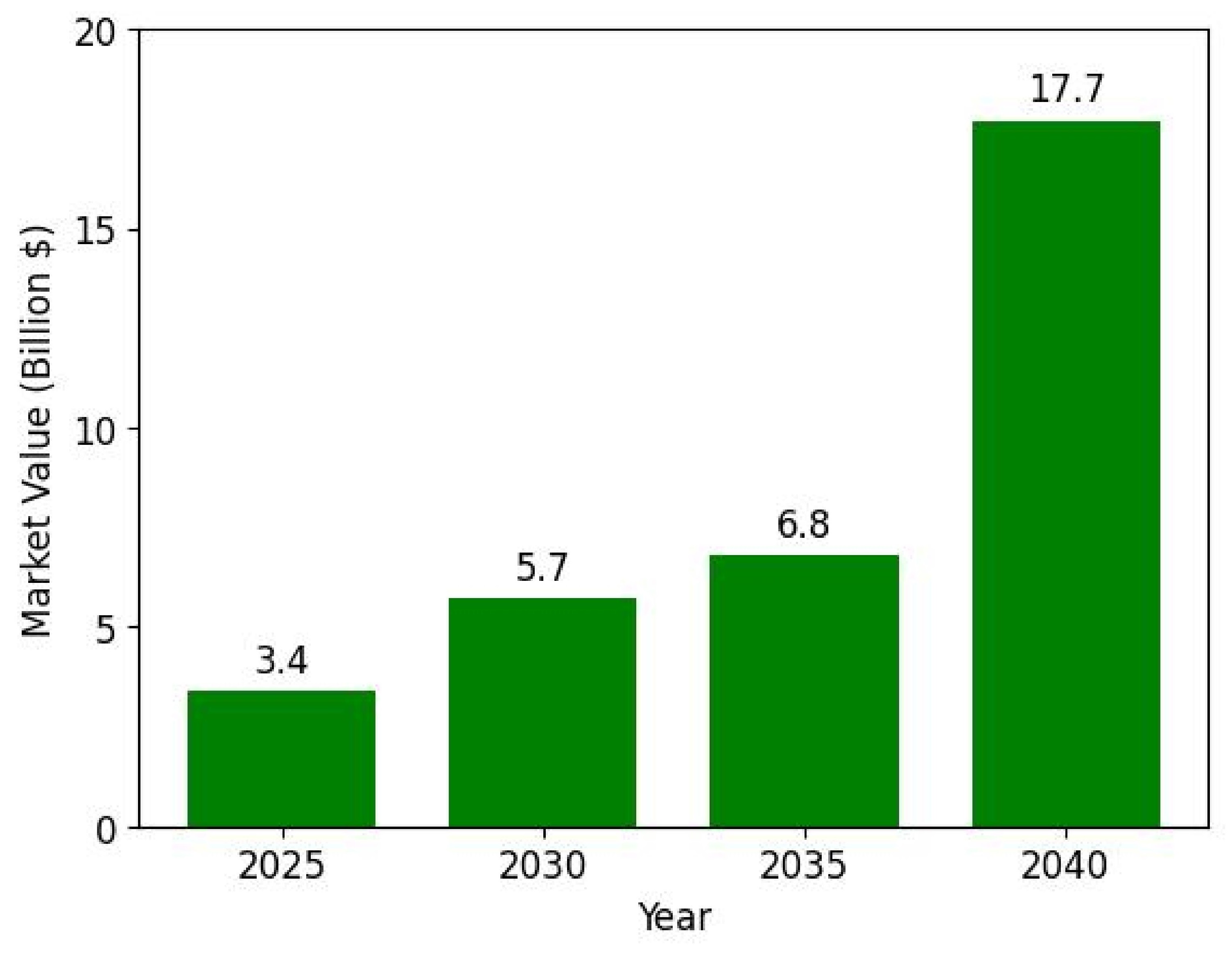}
    \caption{Estimated US passenger eVTOL market size.}
    \label{fig0}
\end{figure}

In recent times, eVTOLs have generated wide interest across industry, government, and academia because of its potential to revolutionize urban mobility. By leveraging the airspace above cities, eVTOLs provide an efficient means of  point-to-point transportation within and between cities. They can fly over congested urban areas, avoiding traffic jams and reducing travel times significantly \cite{rothfeld2021potential}. Unlike conventional ground-based vehicles, eVTOLs do not require lengthy, physical roadways and associated extensive infrastructure investments \cite{yazan2021technology}. Their vertical takeoff and landing capability will allow eVTOLs to access hard-to-reach areas, such as remote or densely populated regions, where traditional transportation options are limited. Furthermore, the transition to eVTOLs aligns with worldwide initiatives to combat climate change and decrease greenhouse gas (GHG) emissions. As fully electric aircraft, eVTOLs emit zero pollutants during flight, which paves the way for a more sustainable transportation ecosystem. Over the long run, it is anticipated that eVTOL operating expenses and trip fares will be reasonably low and comparable to that of ground-based vehicles \cite{goyal2018urban}. Given these significant potential benefits of eVTOLs, the demand for eVTOLs is projected to grow rapidly in the near future, according to several studies \cite{hasan2019urban, doo2021nasa}. According to \cite{delloitte0}, the US passenger eVTOL market size is projected to reach approximately \$17 billion by 2040. This growth, as shown in Figure \ref{fig0}, is expected to be driven by a compound annual growth rate of 11.6\% for the period 2025-2040. It indicates the substantial potential of the future market for electric vertical takeoff and landing aircraft manufacturers (EMs). Note that all monetary  amounts reported in our study are in terms of US dollars. Major industry players, including Airbus, Boeing, Bell, Embraer, Intel, Amazon, Honda, Toyota, and Uber, are partnering with EMs and channeling significant investments, surpassing the \$1 billion mark, towards the design, development, and manufacturing of eVTOLs \cite{doo2021nasa}.

\subsection{Motivations for Supply Chain Planning for eVTOL Manufacturers}

Manufacturers need to develop a supply chain planning that allows them to strategically manage the flow of goods, services, and information from procuring raw materials from suppliers to delivering finished products to customers \cite{beamon1998supply}. The more optimized the supply chain planning is, the better the manufacturers can make strategic decisions and maximize their profits. It is necessary to optimize the supply chain planning for EMs to enable them to meet customer demands on time while maximizing their operating profits, with customers including AAM operators such as air taxi service providers (e.g., Uber Elevate) and logistics and delivery companies (e.g., UPS). However, supply chain planning dedicated to EMs remains largely unexplored, warranting research and analysis in this emerging area. 

The eVTOL supply chain would exhibit several distinctive characteristics compared to traditional aviation and other industries like automotive. A key distinction is the emphasis on manufacturing high-quality eVTOLs at high volumes within a shorter time frame. This requires balancing aviation standards with increased production rates. Additionally, the allocation of emerging technology, machinery, and skilled labor for eVTOL manufacturing is crucial, introducing changes in supply chain planning due to the industry's early stage and lack of a fixed supply chain structure. Consequently, EMs face challenges in predicting and planning resources, impacting manufacturing schedules and costs. Another difference is the limited number of suppliers for eVTOL parts—such as fuselage, wing, v-tail, propeller, battery, motor, and seat—making the supply chain more sensitive to disruptions \cite{scm_motiv2}. The scarcity of critical components, like lithium-ion batteries for electric propulsion, further highlights the need for thorough supply chain planning \cite{s2022lithium}. This necessity is intensified by ongoing contract renegotiations by suppliers for major aerospace companies like Boeing and Airbus, who are under pressure to renegotiate unprofitable contracts due to slow production, high inflation, and post-pandemic challenges \cite{Flight_Global}. The eVTOL supply chain also faces uncertainties in transportation modes, delivery plans, logistics partners, and contingency plans for disruptions \cite{scm_motiv1}. These challenges complicate the timely manufacturing and delivery of pre-ordered eVTOLs, leading to significant implications such as large penalty fees for delays \cite{wang2024impact}. Given these critical issues, research on supply chain planning for EMs in the uncertain and dynamic AAM environment is essential. Addressing this research gap, our study aims to develop a model that can optimize supply chain planning for EMs. To achieve this, we first identify the challenges that EMs are expected to face in their supply chain planning.

\subsection{Challenges in Supply Chain Planning for eVTOL Manufacturers}


To profitably meet the projected future demands of the emerging eVTOL industry, EMs require a comprehensive understanding of supply chain planning challenges. These include \textit{strategic supplier selection}, \textit{procurement planning}, \textit{manufacturing scheduling} and \textit{inventory management}. Strategic supplier selection \cite{jain2009strategic} is a process of identifying optimal suppliers capable of providing necessary eVTOL parts, and procurement planning determines the number and timing of orders to ensure timely part acquisition \cite{cheaitou2015integrated}. The eVTOL industry would heavily rely on contractual partnerships with suppliers of different eVTOL parts, each with varying capacities, lead times (the time interval during which suppliers receive confirmation of an order and prepare the order for dispatch), prices, and qualities. The strategic supplier selection and procurement planning should align with manufacturing scheduling \cite{qin2020scenario}, which is a process of determining the number and timing of eVTOLs to be manufactured by EMs. Due to limited manufacturing capacity, EMs would face a challenge in this process to ensure the timely fulfillment of demand for eVTOLs. Inventory management is also essential for optimizing supply chain planning \cite{esteso2023reinforcement}. The EMs should maintain the right balance between having the necessary eVTOL parts readily available to meet the demand and avoiding excessive inventory costs. With limited inventory capacity and high inventory costs for both eVTOL parts and the eVTOLs, the inventory management would become challenging for EMs. With the evolving nature of the AAM market introducing uncertainties in prices offered by suppliers for eVTOL parts, capacities of suppliers, eVTOL manufacturing costs, and customer demand for eVTOLs, these challenges would become increasingly difficult for EMs to manage.





Another challenge in supply chain planning for EMs is \textit{transportation mode selection} for receiving eVTOL parts from suppliers and delivering eVTOLs to customers. The selected transportation mode should minimize transportation costs while ensuring timely receipt of eVTOL parts and delivery of eVTOLs. This decision-making process considers factors such as the availability of transportation modes in a given origin and destination, shipping distance, delivery speed, restrictions on driving hours imposed on transportation modes, weight of shipping eVTOL parts and eVTOLs, and associated costs for each transportation mode. However, selecting a transportation mode without considering emissions during transportation can lead to significant environmental pollution. Transportation is a major contributor to carbon emissions, which play a significant role in climate change \cite{climate_transportation}. While eVTOLs are envisioned to operate with zero emissions during flight operations, the transportation of eVTOL parts from suppliers and the delivery of eVTOLs manufactured to customers still rely on traditional transportation modes like cargo trucks, cargo trains, cargo ships, and cargo aircraft. Research emphasizes the importance of government regulations, including carbon taxes and credits, in driving emission reduction efforts and promoting a sustainable supply chain \cite{xu2022assessing, liu2022government, green_sc}. Governments have imposed regulations on companies to encourage emission reduction during transportation \cite{EPA}, requiring adoption of more sustainable practices to meet reduction targets \cite{regulations}. Consequently, EMs may incur costs for carbon taxes and credits due to emissions during transportation and would be incentivized to select transportation modes not only to minimize transportation costs but also to reduce emission costs. One more challenge which holds significant importance in the aviation industry that must be considered in the supply chain planning of EMs is the adherence to quality standards for eVTOLs. This \textit{quality control} is essential to ensure that eVTOLs meet the quality requirements set by the federal aviation administration (FAA), enabling them to obtain certifications \cite{klyde2023developing}. Meeting quality standards helps EMs gain customer trust, leading to increased pre-orders for their eVTOLs \cite{lina2022improving}.

While various studies and research have focused on AAM aircraft design (e.g., \cite{silva2018vtol}), concepts of operation (e.g., \cite{thipphavong2018urban}), market studies (e.g., \cite{hasan2018urban, goyal2018urban}), air traffic management (e.g., \cite{FAA2020}), operations planning (e.g., \cite{shihab2019schedule, shihab2020optimal, deep_dispatch}), surveillance infrastructure (e.g., \cite{dulia2024designing, open_framework_standards, PPP, mattei2024improving}), the supply chain planning for the EMs remains an area that has not been thoroughly explored yet. NASA has already emphasized the need for studies on eVTOL supply chain \cite{nasa_scm, AAM_supplychain, nasa_scm1}. However, as far as the authors of this study are aware, there is currently no article available in the existing literature that specifically focuses on supply chain planning for EMs while considering all the challenges mentioned above collectively. To bridge this research gap, our study is the first, to the best of our knowledge, to propose a \textit{stochastic sustainable supply chain optimization and planning for electric vertical takeoff and landing aircraft manufacturers} (3SCOPE) model. The objective of the 3SCOPE model is to maximize the operating profit of an EM while addressing strategic supplier selection, procurement planning, manufacturing scheduling, inventory management, transportation mode selection, and quality control of eVTOL parts and eVTOLs manufactured within the integrated supply chain planning of the EM. We consider carbon emission costs and tax credits in the computation of the operating profit to select the optimal transportation mode for each delivery. One key feature of the 3SCOPE model is its ability to handle uncertainties associated with the model parameters, such as prices offered by suppliers for eVTOL parts, capacities of suppliers, eVTOL manufacturing costs, and customer demand for eVTOLs. To account for demand uncertainty, we incorporate both endogenous demand uncertainty, which stems from the quality of eVTOL parts and eVTOLs manufactured, and exogenous demand uncertainty arising from the market growth of eVTOLs. A scenario-based two-stage stochastic programming approach is used to address these uncertainties in the 3SCOPE model. This is a mathematical optimization framework designed for the decision-makers to make decisions under uncertainties \cite{mateo2023managing}. Furthermore, we employ a rolling horizon approach to enable the 3SCOPE model to make real-time adjustments to evolving market conditions and unexpected disruptions in the supply chain of the EM. As stochastic optimization becomes more computationally expensive with increasing problem scale \cite{almeida2021decomposition}, we utilize a multi-cut Benders decomposition method to allow the EM to solve the 3SCOPE model at a larger scale. We also propose three benchmark models as potential alternatives for solving the eVTOL supply chain problem, and compare the performance of these models with the 3SCOPE model. The performance factor in our study is the operating profit of the EM generated by these models. We also run several numerical cases considering potential disruptions caused due to geopolitical and environmental factors associated with the eVTOL industry. Our observations indicate that the 3SCOPE model is capable of assisting the EM in responding to these disruptions.

\subsection{Summary of Contributions}

This paper addresses various challenges expected to be faced by the EM and offering optimal solutions to optimize supply chain planning with a maximized operating profit for the EM. Our key contributions are as follows:

\setlist[enumerate]{left=0pt}
\begin{enumerate}[label=\alph*)]

\item We introduce the 3SCOPE model, with the primary goal of maximizing the operating profit of the EM. This model serves as a comprehensive decision-making process designed to address a range of challenges within supply chain planning for the EM. 

  
\item  We consider the uncertainties in prices offered by suppliers for eVTOL parts, capacities of suppliers, manufacturing costs of eVTOL, and customer demand for eVTOLs. We implement a two-stage stochastic mixed-integer programming approach within the framework of the 3SCOPE model. To enable the 3SCOPE model to adapt to the dynamic market, we implement a rolling horizon approach, while also employing an multi-cut a Benders decomposition method to handle the large-scale supply chain problems that the EM would face.

\item  We compare the profits generated by the 3SCOPE model with three benchmark models: a deterministic model, a stochastic heuristic model, and a stochastic sequential model. Comparing the deterministic model with the 3SCOPE model helps understand the deviation of profits due to uncertainties associated with the EM's supply chain planning problem. The other two benchmark models are stochastic, offering potential alternatives to the 3SCOPE model. Our findings reveal that the 3SCOPE model consistently outperforms these stochastic models.


\item We perform numerical analyses to illustrate how the 3SCOPE model manages disruptions and adaptations across various situations, such as analyzing how the quality of eVTOLs affects customer demand, assessing supply chain disruptions due to optimal supplier sanctions, raw material scarcity, and delays in the supply chain. Additionally, we analyze the performance of eVTOL manufacturing under increased demand, similar to the automotive manufacturing case. We observe that the 3SCOPE model can handle these disruptions and offer optimal solutions accordingly, assisting the EM in making decisions to respond to these disruptions. 


  
    
    
   
\end{enumerate}
 


\subsection{Outline of Our Paper}

The subsequent sections of this article are organized as follows. Section \ref{sec2} provides an analysis of the relevant literature on supply chain challenges in both the aviation and automotive industries. It also underscores the significance of our study in contributing to the literature. Section \ref{prob_stat} describes the supply chain planning problems that we focus on solving in this study for the EM. Section \ref{sec3} introduces the 3SCOPE model and three benchmark models, providing a detailed methodology. Next, Section \ref{sec4} presents the numerical cases and discusses the results generated from each numerical case. Section \ref{insights} highlights the insights for the managers and stakeholders. Finally, in Section \ref{sec5}, the paper concludes by summarizing the key findings from the analysis and discussing potential avenues for future research.

\section{Literature Review} \label{sec2}



Supply chain planning is a critical operational need for manufacturers across all industries, including aviation and automotive industries, which is essential for meeting customer pre-orders on time while meeting quality standards, reducing operating costs, and increasing operating profits. However, the existing literature lacks studies analyzing supply chain planning for the EM and developing models to aid in its decision-making processes. The EM currently are in two main situations: 1) the eVTOL industry is still in its early stages, with the EM involved in development and testing, yet to commence commercial delivery of eVTOLs; and 2) there is no standard supply chain established for eVTOL manufacturing \cite{delivery_joby_evtol}. Consequently, the EM currently requires more time to manufacture eVTOLs and fulfill customer pre-orders, similar to the aviation industry rather than the automotive sector, even though the size of the eVTOL is not as large as a traditional commercial aircraft. However, as the demand for eVTOLs is expected to rise in the coming years, the EM should scale its manufacturing capabilities accordingly. 
This rise in demand requires the EM to exceed the manufacturing rate of the current aircraft industry and aim for production rates comparable to those in the automotive sector \cite{lineberger2021advanced}. 
Therefore, our focus is on developing solutions appropriate for addressing both aviation and automotive supply chain challenges expected to be faced by the EM.

\subsection{Expected Supply Chain Challenges Ahead for eVTOL Manufacturers}


There is a growing need to examine how strategic supplier selection and procurement planning positively affects manufacturers' supply chain planning \cite{moses2008problems}. For insights into how strategic supplier selection and procurement planning positively affect manufacturers' supply chain planning, readers are referred to the following articles: \cite{moses2008problems, tao2018digital, chiang2012empirical, mocenco2015supply, kotabe2023relationship}. These challenges include uncertainties concerning suppliers' capacity and pricing of parts in the aviation and automotive industries \cite{dweiri2016designing, hashemi2015integrated}, which the EM is also expected to face in its supply chain planning. Along with this, another main challenge for the EM in strategic supplier selection and procurement planning would be the scarcity of certain metals required by battery suppliers to manufacture batteries for eVTOLs \cite{bills2023massively}, which are one of the vital eVTOL parts. This scarcity can lead to a shortage of suppliers in the market, potentially disrupting the supply chain planning of the EM. The next challenge in supply chain planning is expected to be faced by the EM is manufacturing scheduling. Effective coordination between procurement planning and manufacturing scheduling is crucial to minimize disruptions and delays \cite{benton2020purchasing}. In aviation manufacturing, a study \cite{liang2023production} focused on reducing order-delivery time, while another study \cite{ma2023multi} aimed to minimize delays and costs in aerospace component manufacturing. Furthermore, ref. \cite{inproceedings} analyzed the challenges in aviation manufacturing schedule management, including the lack of integrated procurement processes and resource scarcity, such as labor shortages, which necessitate overtime to meet demand \cite{chigbu2021future, sharma2012production}. The EM is also expected to face significant inventory management challenges in its supply chain planning, a critical requirement for modern manufacturers \cite{khan2023relationship}. Excessive inventory ties up capital and escalates storage expenses \cite{bhattacharya2021working}, while leaner inventory levels risk stockouts and increased downtime \cite{derhami2021assessing}. To address these challenges, ref. \cite{singh2018inventory} presented a conceptual methodology for strategic management of material flow, storage, and distribution throughout the supply chain. In \cite{dillon2019study} the authors focused on minimizing inventory costs while ensuring timely production. Additionally, in \cite{masoud2016integrated} an optimization model was developed to minimize total costs, integrating production and transportation planning for these industries.

In the aviation industry, quality control is a major challenge due to the safety-critical nature of its products \cite{article, FINDLAY200218}. Manufacturers must consistently produce safe and reliable products to comply with FAA quality requirements \cite{roca2017risks}. The strategic supplier selection and procurement of aircraft parts without adequate consideration of quality poses risks \cite{q1}. To address this, NASA centralizes aircraft parts procurement and emphasizes compliance with manufacturing standards \cite{q2, q3}. The EM should integrate quality control into its supply chain planning to meet the rigorous standards set by NASA and the FAA for eVTOL parts and eVTOLs. Another significant challenge in EM's supply chain planning is the evaluation of all transportation modes for each route for delivering eVTOL parts from suppliers and the final products to customers, as well as selecting the optimal transportation mode for each delivery. Studies such as \cite{alkahtani2023modified} developed a model for minimizing the cost of delivering parts of the electric motors, which are essential parts of electric vehicles, considering only trucks as a transportation mode. However, the model was not developed to generate supply chain planning needed for the EM; it solely analyzed the optimal location for a pooling hub, where every supplier collects materials, near both the suppliers and the engine factory. In addition to this transportation challenge, the integration of sustainable logistics into EM's supply chain planning is also crucial. Despite the emission-free nature of eVTOLs during operations, their parts are transported using traditional modes such as ground, rail, air, and water, which are not emission-free. The current literature lacks models specifically addressing this challenge in the aviation and automotive industries, although generic studies do exist. For instance, \cite{barzinpour2018dual} developed a model balancing transportation costs and emissions by optimizing manufacturing locations and shipment quantities. Similarly, \cite{rad2018novel} proposed a model for green and cost-efficient transportation mode selection, considering factors like cost, transit time, and environmental impact.

\subsection{Addressing Research Gaps} \label{gaps}

Although research has been conducted addressing supply chain challenges in the aviation and automotive industries, a comprehensive study focusing on the supply chain planning of the EM has not been conducted yet. To bridge this gap, we identify the potential challenges that the EM is expected to face in its supply chain planning and introduce the 3SCOPE model to enable the EM to optimize its supply chain planning. The goal of the 3SCOPE model is to enhance sustainability in the supply chain planning of the EM while maximizing its operating profit. Many studies have developed models for supply chain planning in generic applications, which are adaptable to various industries. These models can be implemented in the aviation and automotive industries after tailoring the model parameters accordingly. However, these existing models also have gaps that we identify in this section to develop a model capable of addressing those. Our study stands out by collectively incorporating several contributions, as outlined in Table \ref{LR_COMP}, compared to other relevant studies. In the table, MCDM refers to multi-criteria decision-making approach; DO to deterministic optimization; and SO to stochastic optimization. 








Researchers have utilized various techniques to tackle the supply chain planning problem, primarily including MCDM and mathematical optimization. Each of these techniques encompasses two main categories of models: deterministic and stochastic. Deterministic models assume all input data, such as demand and manufacturing costs, is known with certainty, and consider the variables as fixed values \cite{liberti2005comparison}. While these models offer simplicity and efficiency, they fail to account for uncertainties such as demand fluctuations and supply chain disruptions. Stochastic models, however, acknowledge the inherent uncertainty in supply chain data by incorporating randomness in the model variables. This enables them to provide a more realistic depiction of supply chain behavior under various scenarios, making decisions more robust to disruptions \cite{pinsky2010introduction}. By comparing the results of stochastic models to deterministic ones, we gain insight into the extent of profit variability caused by uncertainties in the supply chain planning challenges faced by EMs, especially eVTOL manufacturers. This comparison helps us understand the impact of uncertainties on profit outcomes, emphasizing the importance of accounting for these uncertainties in the planning process \cite{koirala2022hosting}. With insight into which situations cause significant profit variation, managers can allocate resources more efficiently to the most uncertain ones. They can strategically invest time, manpower, and funds into designing supply chain planning aimed at minimizing the impact of uncertainties on profitability.

The MCDM approach has commonly been implemented to address mainly for strategic supplier selection challenge in supply chain planning problem. For example, the authors of \cite{dweiri2016designing} proposed a decision support model for strategic supplier selection in the automotive industry of Pakistan, utilizing the analytic hierarchy process (AHP), which one of the MCDM processes. In this paper, criteria such as price, quality, delivery, and service are ranked using AHP, followed by the identification and ranking of sub-criteria. But their model is not an integrated model incorporating the other challenges such as procurement planning, manufacturing scheduling, inventory management. It is a deterministic model as it did not consider uncertainties in model parameters such as price, capacities, manufacturing cost, demand etc. An example of a stochastic MCDM model can be found in \cite{foroozesh2021new}, where the authors introduced a model for addressing the strategic supplier selection problem in uncertain environments. This model employed Monte Carlo simulation to manage uncertainty; however, it is essential to note that while the model provides insights into the probabilistic behavior of the system, it does not inherently optimize decision variables or prescribe specific actions to achieve a specific objective, such as minimizing total supply chain costs or maximizing profits.

Mathematical optimization has been used to address challenges in supply chain planning. For instance, \cite{kannegiesser2014sustainable} addressed manufacturing scheduling and transportation. In this article, the authors developed a deterministic model which included variables to identify the optimal quantity and timing for manufacturing products and transportation, taking carbon emissions into account. Another deterministic model developed in \cite{rad2018novel} addressed strategic supplier selection and procurement planning, inventory management, and transportation challenges to identify the optimal timing and quantity in the manufacturing process, inventory, and transportation. On the other hand, in \cite{ghadimi2023safety}, the authors developed a deterministic model where they considered strategic supplier selection, procurement planning, manufacturing scheduling, and inventory management challenges. In another study \cite{barzinpour2018dual}, a deterministic model was developed to address manufacturing scheduling and transportation challenges. The model aimed to determine the optimal scheduling of manufacturing processes and deliver final products to the customer via different transportation modes, including penalties if the manufacturing plant lost sales.

The deterministic models struggle to handle uncertainties inherent in supply chain planning across various applications. Consequently, research has focused on addressing these uncertainties using stochastic models. For instance, ref. \cite{gruler2018combining} developed a model tackling the inventory routing problem by addressing transportation and inventory challenges with stochastic demands. The model was solved with a variable neighborhood search framework, a simheuristics approach. Additionally, the authors benchmarked their model against other greedy heuristics, including selecting inventory stock and delivery strategies that yield the highest reduction in inventory and delivery costs for the manufacturer. However, greedy stochastic approaches are not robust enough to handle uncertainties, as they prioritize immediate cost reductions or benefits in each decision step and do not consider the long-term impact of these decisions on the entire supply chain under uncertain future conditions. In contrast, stochastic programming offers a more robust solution by finding strategies that minimize expected costs or maximize expected profits across different scenarios. For example, the authors of \cite{fattahi2022data} developed a stochastic model using two-stage stochastic programming to address inventory challenges, uncertain demand, and penalties for delaying deliveries to customers. They employed a rolling horizon approach to continuously update model parameters based on changing market conditions and analyzed the model's resilience to supply chain disruptions. In another study \cite{almeida2021decomposition}, the authors developed a two-stage stochastic model aiming to address procurement planning from a given set of suppliers, transportation of parts to manufacturing plant, and inventory management of these parts. However, they did not consider the challenges related to optimal supplier selection, manufacturing scheduling, and the transportation of manufactured products to customers. This study accounted for penalties for delayed delivery, uncertainties in customer demand, supplier prices and capacities, and manufacturing costs. They also used a rolling horizon approach to update model parameters in response to market conditions and employed Benders decomposition method to handle large scale problems. Another instance of a two-stage stochastic model was presented in \cite{inproceedings}, where the authors addressed strategic supplier selection, procurement planning, inventory management, and transportation challenges under uncertain demand. Among the existing studies in literature, this study integrated most of the challenges into a stochastic model. In this study, the authors also developed a deterministic version of the model and compared it with the stochastic model to analyze how uncertainties affected the solutions in real-world scenarios compared to situations where uncertainties were not considered. However, this study did not address manufacturing scheduling challenges, such as determining the optimal timing for ordering and receiving parts from suppliers, as well as the optimal quantity and timing for manufacturing customer orders. Additionally, the study did not consider the impact of transportation emissions on transportation selection, nor did it account for uncertainties in supplier prices, capacities, and manufacturing costs. Moreover, their model was not capable of solving large-scale problems. 

The existing stochastic models for supply chain planning have been formulated to address only certain challenges. Consequently, a manufacturing plant manager must employ additional independent models for the remaining challenges \cite{kaur2020sustainable}. For instance, if the manager utilizes the model developed in \cite{inproceedings}, they can determine optimal suppliers and procurement quantities, as well as the optimal transportation mode to meet the demand. Following this, the manager must sequentially run another model for manufacturing scheduling, using the outputs of the first model as inputs. Moreover, the existing studies did not incorporate a model to address the missing challenges of the challenges anticipated for the EM, which could be run sequentially with their main model. Our main contribution is to build the integrated 3SCOPE model, which addresses all the potential challenges including strategic supplier selection, procurement planning, manufacturing scheduling, inventory management, sustainable and cost-effective transportation mode selection, and quality control of eVTOL parts and eVTOLs manufactured by the EM. Reflecting real-world challenges in the aviation industry, we consider penalties to be paid by the EM to customers if eVTOL deliveries are delayed beyond the deadline. In developing the 3SCOPE model, we incorporate uncertainties into the model parameters, including supplier prices for eVTOL parts, supplier capacities, eVTOL manufacturing costs, and customer demand for eVTOLs. To allow the 3SCOPE model to make real-time adjustments and adaptations to changing market conditions, demand fluctuations, and unforeseen disruptions, we implement a rolling horizon approach. To handle large-scale supply chain problems, we employ a multi-cut Benders decomposition method to run the 3SCOPE model. Furthermore, we construct benchmark models and compare their performance with that of the 3SCOPE model across several numerical cases. We also address supply chain planning challenges for the EM to assist them in responding to disruptions inherent in this evolving industry, such as geopolitical sanctions, raw material shortages, and supplier unavailability in the market. These combined contributions distinguish our study, setting it apart from other research in the existing literature.




\begin{spacing}{0.9}
\begin{longtable}[htbp]{lccccccccccc}
    \caption{Comparison of our study with other relevant studies.} 
    \label{LR_COMP} \\
    \hline
    \textbf{Contributions} & \cite{dweiri2016designing} &  \cite{foroozesh2021new} & \cite{kannegiesser2014sustainable} &  \cite{rad2018novel} & \cite{ghadimi2023safety} & \cite{barzinpour2018dual} & \cite{gruler2018combining} & \cite{fattahi2022data} &  \cite{almeida2021decomposition} & \cite{inproceedings} & \textbf{Our}  \\
       & &   &   &  & &  &&& &  &  \textbf{Study} \\
    \hline
    \endfirsthead
    \hline
     \textbf{Contributions} & \cite{dweiri2016designing} & \cite{foroozesh2021new} &  \cite{kannegiesser2014sustainable} &  \cite{rad2018novel} & \cite{ghadimi2023safety} &  \cite{barzinpour2018dual} & \cite{gruler2018combining} & \cite{fattahi2022data} &  \cite{almeida2021decomposition} & \cite{inproceedings} &  \textbf{Our} \\
      & &   &   &   &  && && & &  \textbf{Study} \\
    \hline
    \endhead
    \hline
    Method  & MCDM & MCDM & DO & DO & DO &  DO & SO & SO & SO & SO & SO \\ \hline
    Optimal suppliers for each part  & \checkmark & \checkmark & X & \checkmark & \checkmark & X & X & X & X & \checkmark  &  \checkmark \\ \hline
         
    Optimal number of parts to be  & X & X &  X & \checkmark & \checkmark & X & X & X & \checkmark & \checkmark & \checkmark \\ 
    
        procured from the selected suppliers  &  &  & &  &  & & & \\ 
       \hline
        
     Optimal timing for ordering parts   & X & X &  X & X & \checkmark &  X & X & X & X & X & \checkmark \\
        from suppliers and receiving them        &  &  && &  &  & & \\
      \hline
         
    Optimal quantity and timing  & X & \checkmark &  \checkmark & X & X & \checkmark &  X & X & X & X & \checkmark \\
        for manufacturing products    &  &  & &&  &  & & \\
        ordered by the customers      &  &  & &  &&  & & \\   \hline

       Optimal transportation mode for & X  & X &  \checkmark & \checkmark & X &  \checkmark &  \checkmark & X & \checkmark & \checkmark & \checkmark \\
         each delivery of parts and products &  & & & &  &  &  & & \\
      \hline
 
          
       Optimal transportation modes  & X & X &  \checkmark & \checkmark & X &  \checkmark &   X & X  & X & X & \checkmark \\
         considering carbon emissions        &  &  & & &  &  & & \\
         \hline
       Optimal quantity of parts and& X & \checkmark &  X & \checkmark & \checkmark &  X & \checkmark &   \checkmark & \checkmark & \checkmark & \checkmark \\
         products in manufacturer's inventory &  &  && & & &  & & \\
   \hline
        
       Penalty for delaying delivery & X & X &  X & X & X &  \checkmark & \checkmark & \checkmark & \checkmark & X & \checkmark \\
         of products to customer      &  &  & &  &  & & & \\ \hline
         
     Endogenous demand uncertainty & X & X &  X & X & X &  X & X & X & X & X & \checkmark \\
     arising from quality of products   &  &  & &  & & & & \\ 
        \hline
         
       Exogenous demand uncertainty  & X & X &  X  & X & X &  X & \checkmark  & \checkmark & \checkmark & \checkmark & \checkmark \\
       arising from market growth of product  &  &  & &  & & & & \\ 
    \hline
         
      Uncertainties in prices and & X & \checkmark &  X & X & X &  X & X  & X & \checkmark & X & \checkmark \\
         capacities of suppliers     &  &  & &  &  && & \\ \hline
         
      Uncertainty in manufacturing cost & X  & X &  X & X & X &  X &  X & X & \checkmark & X & \checkmark \\
       \hline
         
     Rolling horizon approach   & X & X &  X &X & X & X  & X & \checkmark & \checkmark & X & \checkmark \\
       \hline

     Scalability for handling & X & X &  X & X & X &  X & X  & X & \checkmark & X & \checkmark \\
     large-scale problems &  &  & &  & & & & \\ \hline
        


   Compared performance of   & X & X &  X & X & X &  X & \checkmark & X & X & \checkmark & \checkmark \\
       proposed model with benchmark models  &  &  & &  &&  & & \\
       
       \hline

   Resilience to disruptions in supply chain & X & X &  X & X & X &  X & X & \checkmark & X & X & \checkmark \\
    
    
        \hline
\end{longtable}

\end{spacing}

\section{Supply Chain Planning Problem for eVTOL Manufacturer} \label{prob_stat}

\subsection{Objective, Sub-problems, and Constraints}


In the EM supply chain planning problem, the main objective of the commercial EM considered is to maximize its profits over a given \textit{analysis horizon (AH)}. A \textit{rolling horizon approach (RHA)} is used to maximize the EM's profits over AH. The AH is composed of a set of time periods $T^P$. 
Successive groups of $V$ subsequent periods of AH are defined as the \textit{planning horizon (PH)}. Given the evolving nature of the AAM market, it is reasonable for the EM to set $V=2$, allowing RHA to adapt effectively for the EM to plan and optimize decisions in the dynamic eVTOL industry. Therefore, the PH includes two periods in each iteration of RHA, with the first period denoted as $F^P$ and the second period $S^P$. Each PH consists of a set of days $T$.

\begin{figure}[b!] 
    \centering
    \includegraphics[width=\textwidth]{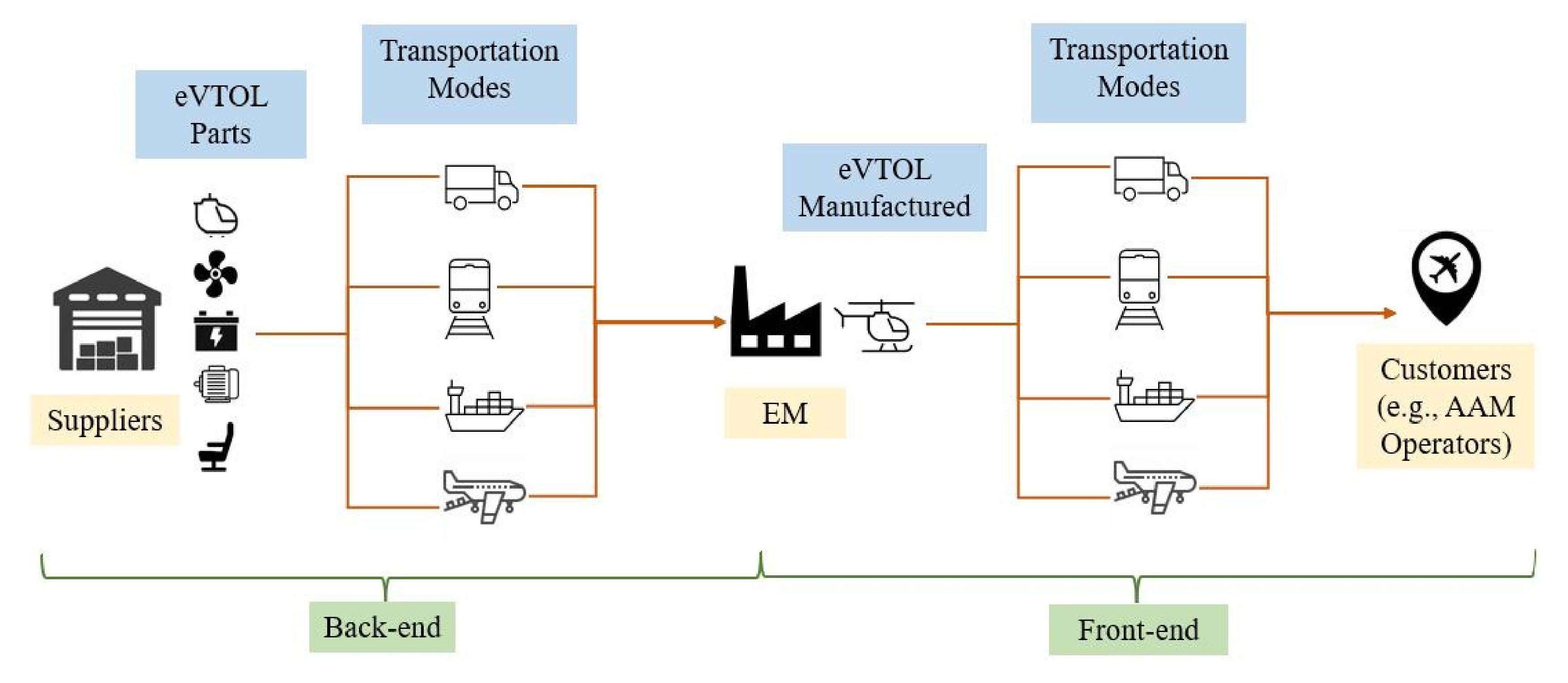}
    \caption{The back-end and front-end of the supply chain of the EM.}
    \label{fig_sc}
\end{figure}

\begin{figure}[tb!] 
    \centering
    \includegraphics[width=\textwidth]{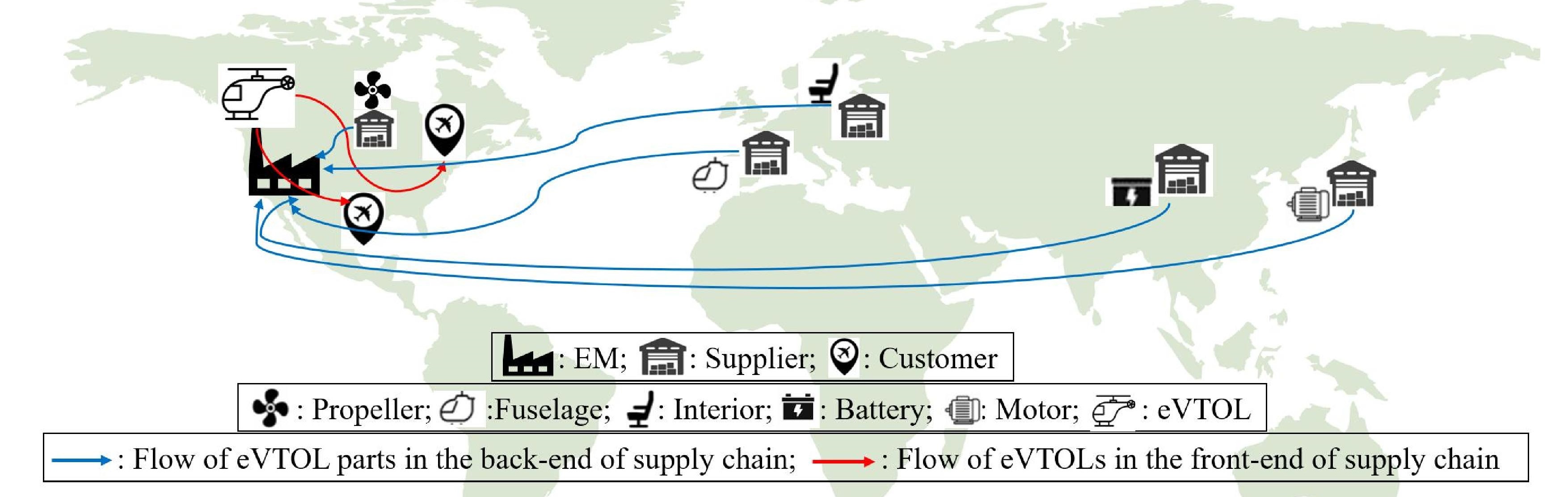}
    \caption{A supply chain network of the EM.}
    \label{fig_sc_network}
\end{figure}

The EM supply chain planning is considered to have a back-end and a front-end, as visualized in Figure \ref{fig_sc}. The back-end of the supply chain involves procuring a set of necessary eVTOL parts, denoted by $P$, required for manufacturing the eVTOLs. Through comprehensive survey of the market, the EM identifies $S_i$, the set of suitable suppliers for each specific part $i \in P$. The back-end of the supply chain also includes the delivery of these parts from suppliers to the EM, utilizing a transportation mode $k$ from a set of transportation modes $M$ available between the supplier and the EM. On the other hand, the front-end of the supply chain is concerned with 1) the manufacturing of the eVTOLs pre-ordered by a set of customers $C$; and 2) the transportation of manufactured eVTOLs via a available transportation mode $k \in M$ from the EM to customer $l \in C$. All AAM operators within a city, with potential for AAM application, are considered as a single customer. Consequently, the number of eVTOL pre-orders from each customer in $C$ represents the total number of eVTOL pre-orders originating from that city. Figure \ref{fig_sc_network} illustrates a supply chain network of the EM on a world map, showing the locations of EM, suppliers, and customers, as well as the flow of eVTOL parts and eVTOLs between them.


To integrate the two ends of the supply chain and address the challenges involved in supply chain planning discussed in Section \ref{intro}, the EM is expected to simultaneously address several sub-problems. These sub-problems encompass determining: 1) which suppliers to select for each eVTOL part (strategic supplier selection); 2) how many eVTOL parts to procure from the selected suppliers (procurement planning); 3) when to place orders with selected suppliers for these number of eVTOL parts (procurement planning); 4) when to receive the ordered eVTOL parts and in what quantities (procurement planning); 5) how many eVTOL pre-orders to accept from customers for a specific time period (manufacturing scheduling); 6) when and how many eVTOLs to manufacture to meet customer pre-orders within given deadlines (manufacturing scheduling); 7) which transportation mode to select, when to schedule deliveries, and how many eVTOLs to dispatch in each delivery of at both ends of the supply chain, considering mode availability and sustainability (transportation mode selection); 8) how many eVTOL parts and eVTOLs to hold in the inventory at both ends of the supply chain, respectively (inventory management); and 9) how many eVTOL pre-orders from customers are expected to be affected by the quality of eVTOLs (quality control). To ensure feasibility of the EM supply chain planning, several constrains, including quality constraints, customer order fulfillment constraints, manufacturing rate constraints, capacity constraints of suppliers of eVTOL parts, and inventory management constraints, have been considered.



\subsection{Customer Pre-Orders of eVTOLs} \label{customer_pre}



The EM is expected to implement a pre-order strategy for receiving eVTOL orders from customers and then manufacturing the eVTOLs according to a given deadline. This is a common order-taking strategy in the eVTOL industry, as reported in \cite{pre_order1, pre_order2}. According to this strategy, the EM supply chain planning is driven by the initial and additional number of eVTOL pre-orders received from the customers in $C$ by the EM. For each period in PH, the EM would receive the initial number of eVTOL pre-orders that need to be manufactured according to the deadline for the corresponding period. In addition to the initial number of eVTOL pre-orders, the EM would continue to accept the additional number of eVTOL pre-orders for the period until it begins, which would be added to that initial number of eVTOL pre-orders.


At the beginning of $F^P$, the EM would stop taking eVTOL pre-orders for $F^P$ and would determine the final number of eVTOL pre-orders, represented by $D^1_{l}$, required to be manufactured and delivered by a specified order-delivery deadline, denoted by $T^{d^1}_{l}$, for customer $l$ in $F^P$. At the beginning of $F^P$, the EM would also receive the initial number of eVTOL pre-orders, represented by $D^2_{l}$, required to be manufactured and delivered by a specified order-delivery deadline, denoted by $T^{d^2}_{l}$, for customer $l$ in $S^P$. However, they would not know the final number of eVTOL pre-orders for $S^P$, as throughout $F^P$, they would continue accepting eVTOL pre-orders from customers for $S^P$. These additional numbers of eVTOL pre-orders will be added to $D^2_{l}$. When $S^P$ begins, the EM would stop accepting additional eVTOL pre-orders for $S^P$, and the number of eVTOL pre-orders accepted up to that point would represent the final number of eVTOL pre-orders for $S^P$. The EM would need to manufacture eVTOLs to fulfill the pre-orders according to the deadlines associated with the final number of eVTOL pre-orders. 


The additional number of eVTOLs pre-ordered by customer $l$ for $S^P$ is divided into two distinct categories. The first category represents the additional number of eVTOL pre-orders driven by the quality of eVTOLs, denoted by $d^2_{l}$. In the EM's supply chain problem, $d^2_{l}$ is an endogenous variable, indicating that its value is driven by other variables within the model \cite{motamed2021multistage}. If the quality of eVTOLs manufactured in $F^P$ exceeds the base FAA quality requirements, the EM would receive $d^2_{l}$ eVTOL pre-orders from customer $l$ for $S^P$. This is due to the positive relationship between quality of products and demand of customers, as mentioned in these studies \cite{cummins1997price, paulley2006demand, soderlund1998customer}. Therefore, the EM would be interested to manufacture high quality of eVTOLs and select suppliers known for producing certified and qualified eVTOL parts. At the beginning of $F^P$, the EM would test and assess the quality of part $i$ received from supplier $j$ and assign a quality value $Q_{ij}$ to it. The quality of these eVTOL parts would affect the quality of eVTOLs manufactured in $F^P$, thereby influencing the value of $d^2_{l}$ received from customer $l$ for $S^P$. The second category consists of additional eVTOL pre-orders for $S^P$, which are subject to exogenous uncertainty arising from the market growth of eVTOLs. This uncertainty cannot be controlled by other variables in the model.


The EM would not accurately predict the additional number of eVTOL pre-orders for $S^P$ until it begins. On the other hand, due to the dynamic nature of the AAM market, there would be several disruptions in the supply chain, such as supplier unavailability, increased lead time of suppliers due to geopolitical conflicts and sanctions, environmental or economic factors. In response to these disruptions, the EM would need to procure more eVTOL parts in $F^P$ than are needed to manufacture eVTOLs for $F^P$. This ensures that the EM has adequate eVTOL parts and eVTOLs manufactured in inventory to satisfy the demand of customers for $S^P$. On the other hand, the inventory cost continues to accumulate as long as inventory is held within the supply chain, and minimizing inventory cost is essential for maximizing supply chain profitability. Therefore, at the beginning of $F^P$, the EM would have to predict the additional number of eVTOL pre-orders for $S^P$ as closely as possible to prevent shortages or excesses of eVTOL parts and eVTOLs in $F^P$. This uncertainty regarding additional eVTOL pre-orders for $S^P$ would add complexity to the EM supply chain planning.

\subsection{Contracting with Suppliers} \label{suplier_contract}


Based on $D^1_{l}$ and $D^2_{l}$ received from customer $l$ by the EM in PH, the EM would select suppliers for the eVTOL parts through strategic supplier selection and establish contracts with them. The contracts would specify the number of each eVTOL part $i$ required in $F^P$ and their procurement prices, though the timing for placing orders can vary. Such contracts are commonly used in the manufacturing sectors \cite{masten1984organization, dimitri2006handbook}. Upon establishing a contract with a supplier, the specified number of eVTOL parts in the contract would not be delivered all at once upon signing. Instead, the EM would proceed to place orders for the required eVTOL parts at various points throughout the agreed time frame in the contract. When an order is placed to supplier $j$ by the EM, the EM would receive eVTOL parts $i$ after the lead time $T''_{ij}$ and delivery time $T'_{ijk}$ to the EM via transportation mode $k$. The number specified in the contract would indicate that the suppliers must provide that number of eVTOL parts at the same price throughout the agreed time frame to the EM.

In the aviation industry, manufacturers traditionally procure aviation parts at discounted prices from suppliers when they place orders well in advance and in large quantities \cite{buyukdaug2020effect, alford2002effects}. However, the aviation part suppliers are currently striving to move away from discount price deals offered to the aviation manufacturers in exchange for increased volume delivered over longer periods, as mentioned in Section \ref{intro}. Therefore, the suppliers would tend to negotiate with the EM not to give the same discounted price over longer periods. In this AAM market, it is more reasonable to consider that the price would remain consistent only during $F^P$, rather than for all periods in PH. At the beginning of $S^P$, the contract needs to be updated to reflect the current market conditions. 


The AAM market can experience fluctuations due to several dynamic factors, including inflation, reduction in the number of suppliers, reduction in the production capacity of chosen suppliers, shortages of raw materials used by suppliers, change in quality requirements for eVTOLs, higher fuel costs in transportation, change in carbon credit regulations, and unavailability of suppliers due to political disruptions. These factors introduce uncertainties in availability and capacities of suppliers, prices of eVTOL parts from the suppliers, and the final number of eVTOL pre-orders received by the EM for $S^P$. Therefore, the EM supply chain planning problem would become more complex, as it is not possible for the EM to know these parameters for $S^P$ at the beginning of $F^P$ when the EM signs the contract with suppliers, specifying the number of eVTOL parts that would be delivered to the EM in $F^P$. For example, if the prices of eVTOL parts offered by suppliers are higher at the beginning of $S^P$ compared to the prices at the beginning of $F^P$, or if there is a possibility of suppliers being unavailable in the market in $S^P$, the EM should procure as many eVTOL parts as possible in advance during $F^P$. This ensures that the EM can fulfill eVTOL pre-orders within the specified deadlines in $S^P$, even in any disruptive situation. However, procuring these eVTOL parts in advance would lead to increased inventory costs across PH. On the other hand, if suppliers are available or their offered prices are lower at the beginning of $S^P$ compared to the prices at the beginning of $F^P$, the EM should procure the required eVTOL parts for $S^P$ at the beginning of $S^P$, thus avoiding unnecessary inventory costs. Since these parameters remain unknown until the end of $F^P$, these uncertainties introduce challenges in the EM supply chain planning.



%

\subsection{Allocating Resources for Manufacturing eVTOLs} \label{manufacturing_contract}

Another source of uncertainty is the eVTOL manufacturing cost for the EM. The eVTOL manufacturing cost depends on various resources required during eVTOL manufacturing and the operation of the EM's plant, including skilled labor, machinery, tools, electricity, and more. At the start of PH, the EM would allocate the resources for both $F^P$ and $S^P$ to maintain workflow stability and prevent disruptions from resource shortages throughout PH. The EM would need to know the final number of eVTOL pre-orders for the periods in PH before allocating the resources. However, the EM would not have precise information about the final number of eVTOL pre-orders for $S^P$. Therefore, the EM would face uncertainty in eVTOL manufacturing cost, as over-allocation of resources would result in excess manufacturing costs, while underestimating the amount of resources needed can lead to shortages later.

\subsection{Contracting with Logistics Companies} \label{logistic_contract}

In addition to establishing contracts with suppliers, the EM would also contract logistics companies for 1) transporting eVTOL parts from suppliers to the manufacturing plant in the back-end of the supply chain, and 2) transporting manufactured eVTOLs from the EM's plant to customer locations in the front-end of the supply chain. Outsourcing these logistics operations is a common practice due to the expenses associated with owning and maintaining a fleet of trucks, including maintenance costs and driver salaries. This becomes especially important when dealing with overseas international suppliers, where truck transportation may not be viable, and alternatives like shipping and air transport are more appropriate. Also, it is not practical to have proprietary vehicles for each transportation mode; hence, it would be better for the EM to contract with corresponding logistics companies. The logistics companies would charge the EM based on the total weight of eVTOL parts and eVTOLs transported, as well as the shipping distance \cite{adenso2023metafrontier}. Therefore, the EM would need to conduct transportation mode selection to determine the optimal mode for each origin-destination pair within both the back-end (involving suppliers and the EM) and the front-end (involving the EM and customers) of the supply chain at any given time. The optimal transportation mode for the same pair may vary over time, in alignment with other challenges in the supply chain planning at that time. Consequently, the transportation mode chosen at a given time within a period may differ from previous selections for that pair within the same period. The EM should also take into account the constraint of maximum daily driving hours imposed on drivers in its supply chain planning before contracting with logistics companies. Maximum daily driving hours refer to the legal limitations on the duration that drivers can spend operating a vehicle on a given day, which is crucial for maintaining the well-being of drivers, promoting safe transportation practices, and avoiding potential penalties and legal issues \cite{hours1}.

At the start of $F^P$, the EM should engage in contracts with the appropriate logistics companies. This ensures that these companies can plan their operations effectively. Such proactive measures prevent the EM from facing logistical issues, such as unavailability, when transporting eVTOL parts and eVTOLs, which could arise without prior contractual arrangements. Alongside total transportation costs, total emission costs should also be considered in the transportation mode selection through carbon taxes and carbon credits associated with the emissions manufactured during transportation as per the regulations, as discussed in Section \ref{intro}. Hence, when determining the optimal transportation mode at each end of the supply chain, the EM would also be incentivized to consider not only which transportation mode would result in lower transportation costs, but also lower total emission costs.

\section{Solution Methodology}\label{sec3}

\begin{figure}[htb!] 
    \centering
    \includegraphics[width=\textwidth]{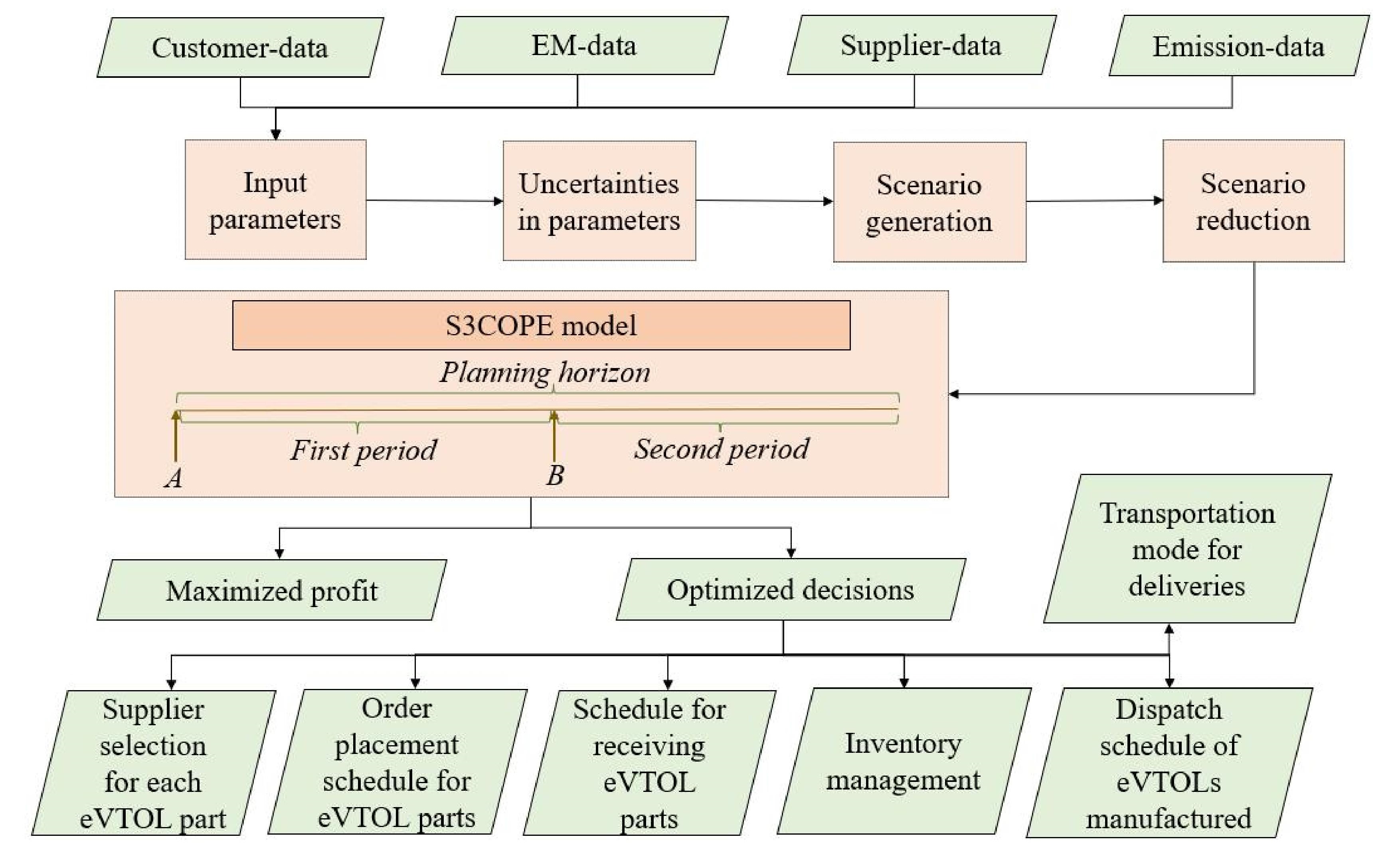}
     \caption{The EM supply chain planning framework.}
    \label{flowchart}
\end{figure}

The supply chain planning framework of the EM, as illustrated in Figure \ref{flowchart}, is proposed to solve the EM supply chain planning problem. The cornerstone of this framework is the 3SCOPE model, a stochastic mixed-integer programming model. The objective of this model is to maximize the operating profit of the EM by providing a decision-making process to address the various challenges of the supply chain planning, including: strategic supplier selection and procurement planning; manufacturing scheduling; inventory management; transportation mode selection considering sustainability; and quality control.

\subsection{Input Data}

The initial step towards solving the EM supply chain planning problem involves collecting data on input parameters for the 3SCOPE model. This data includes information about customers, the EM, suppliers, and emissions. The customer-data encompasses a variety of crucial details, including  initial number of eVTOL pre-orders and their associated deadlines; their geographical locations; additional number of eVTOL pre-orders subject arising from the market growth of eVTOLs; the available transportation modes for delivering eVTOLs manufactured at the EM's facility to each customer's destination; the shipping distances for each transportation mode; and the corresponding time required for delivery. The EM's data includes manufacturing capacity, manufacturing cost per eVTOL, maximum daily manufacturing rate, and eVTOL selling price. The supplier-data includes details such as their geographic locations, prices of eVTOL parts lead times, capacities, and qualities of eVTOL parts. The emission-data accounts for carbon emissions associated with different transportation modes and shipping distances, as well as carbon credits related to these emissions as specified in regulations. 



\subsection{Rolling Horizon Approach} \label{rollingHorizon_approach}

In the EM supply chain planning problem, many of the input variables, including the availability and capacities of suppliers, prices of eVTOL parts, customer demand for eVTOLs, availability of manufacturing resources, and transportation mode availability at both ends of the supply chain, can change throughout AH in response to the various dynamic market factors discussed in Section \ref{prob_stat}.
Therefore, we employ RHA to enable the EM to respond effectively to dynamic changes in the AAM market and support optimal decision-making throughout AH. The RHA is outlined in Algorithm \ref{alg:rolling_horizon}. 
After each iteration of this approach, the 3SCOPE model finds solutions for both $F^P$ and $S^P$ of the PH, but it only saves the solutions for $F^P$. Then, the 3SCOPE model updates the model parameters using the newly available information for $S^P$ and the subsequent periods in $T^P$, and proceeds to the next iteration. This enables the 3SCOPE model to consistently generate optimal solutions by utilizing RHA.

\begin{spacing}{1}
\begin{algorithm} [htb!]
\caption{Rolling Horizon Approach}
\label{alg:rolling_horizon}
\begin{algorithmic}[1]

\State \textbf{Initialization:}
\begin{itemize}
    \item Initialize a set of next periods, $S^P$, as an empty set. 
    \item Initialize a set of subsets, $O$, as an empty set, with each subset representing the optimal solutions for individual periods within $T^P$. 
\end{itemize}

\For{$h$ from 1 to $|T^P|$} \Comment{Loop through each iteration $h$ of RHA.}
    \State Set $F^P = h$. \Comment{Set the first period of iteration $h$.}
    \If{$h < |T^P|$} \Comment{If iteration $h$ is not the last iteration.}
        \For{$i$ from $F^P + 1$ to $F^P + V - 1$} \Comment{Iterate through subsequent periods in iteration $h$ following $F^P$.}
            \State Add $i$ to the set $S^P$. \Comment{Add subsequent periods to $S^P$.}
        \EndFor
        \If{the last period in $S^P > |T^P|$} \Comment{Adjust the last period in iteration $h$ if it exceeds AH.}
            \State Set the last period in $S^P$ to $|T^P|$.
        \EndIf 
        \State Set $P^H = A^H(F^P \cup S^P)$. \Comment{Set PH.}
        \State Update the input parameters of the 3SCOPE model for PH. \Comment{Update input parameters.}
        \State Run the 3SCOPE model for PH. \Comment{Run the 3SCOPE model for PH in iteration $h$.}
        \State Obtain optimal solutions for $F^P$ and for the periods in $S^P$. \Comment{Obtain optimal solutions for PH in iteration $h$.}
        \State Add optimal solutions for $F^P$ to $O_h \in O$. \Comment{Add optimal solutions to the set of optimal solutions.}
    \Else                       \Comment{If iteration $h$ is the last iteration.}
        \State Run the 3SCOPE model only for $F^P$. \Comment{Run the 3SCOPE model only for the last period in $T^P$.}
        \State Add optimal solutions for $F^P$ to $O_h \in O$. \Comment{Add optimal solutions to $O$.}
    \EndIf
\EndFor

\For{each iteration $h$ in RHA}
    \If{Lead times of suppliers increase in iteration $h$}
        \For{each supplier $j$ in the set $S$}
            \For{each part $i$ in the set $P$}
                \State Observe $T''_{ij}$.
                \For{each transportation mode $k$ in the set $M$}
                    \State Observe $T'_{ijk}$.
                \EndFor
            \EndFor
        \EndFor
        \For{each customer $l$ in the set $C$}
            \State Determine the number of days needed to manufacture eVTOLs for customer $l$.
            \For{each transportation mode $k$ in the set $M$}
                \State Observe $T'_{lk}$.
            \EndFor
        \EndFor
        \State Renegotiate customer deadlines and penalty cost for the delay.
        \State Increase $|T|$ based on all observed information.
    \EndIf
    \State Continue with the remaining steps of RHA.
\EndFor

\State \textbf{Output:} $O$. \Comment{Output the set of optimal solutions.}
\end{algorithmic}
\end{algorithm}

\end{spacing}

\indent 

\begin{figure}[htbp] 
    \centering
    \includegraphics[width=\textwidth]{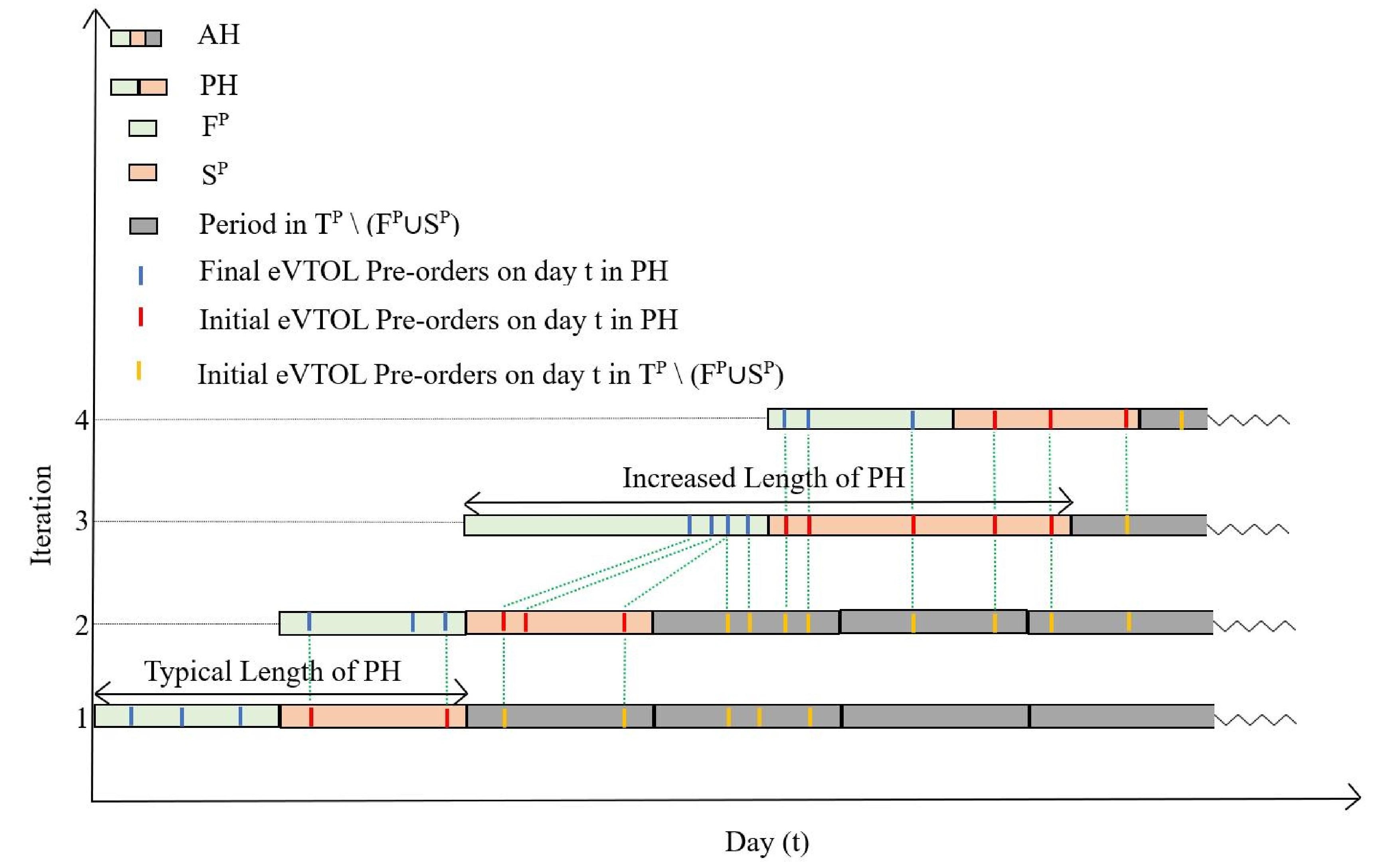}
     \caption{The rolling horizon approach.}
    \label{RHA}
\end{figure}

A demonstration of RHA, as depicted in Figure \ref{RHA}, illustrates how the approach progresses through the periods in $T^P$ during each iteration across the days in AH. Each iteration encompasses a typical length of PH, involving both $F^P$ and $S^P$. The green dotted lines illustrate the transition of the initial and final numbers of eVTOL pre-orders from one iteration to the next. The initial number of eVTOL pre-orders for $S^P$ in an iteration becomes the final number of eVTOL pre-orders in the next iteration, as the additional eVTOL pre-orders are added. For example, at the beginning of the first iteration, the EM has three final eVTOL pre-orders for $F^P$ and two initial eVTOL pre-orders for $S^P$. While there are eVTOL pre-orders from customers for other periods in $T^P$, they are not included in PH in this iteration. During the first iteration, the EM receives an additional order for $S^P$. Subsequently, at the beginning of the second iteration, the newly available information is updated for all the periods. In this iteration, the $S^P$ from the first iteration transitions into $F^P$ for the second iteration, including the three final eVTOL orders for this period. The third period in $T^P$ becomes $S^P$ in this iteration. The RHA iterates similarly; however, in any iteration, the EM would find that due to disruptions in the AAM market, the lead times of all suppliers for any eVTOL part have increased. As a result, the required eVTOL parts would not be available to the EM for manufacturing eVTOLs and fulfilling customer pre-orders during that period. In this situation, at the beginning of the iteration, the EM should increase the length of PH, renegotiate with the customers, and delay the deadlines over the extended PH, as illustrated in Figure \ref{RHA}. The EM would need to pay penalty fees to the customers for delaying the deadlines in this situation, which depend on the conditions of the deal made with customers before signing the agreements. However, if the situation reverts to its typical market conditions, the length of PH should be decreased as before. 

\subsection{Scenario-based Two-stage Stochastic Programming} \label{model_approach}

A key feature of the 3SCOPE model is addressing uncertainties associated with the prices offered by suppliers for eVTOL parts, capacities of suppliers, eVTOL manufacturing costs, and customer demand for eVTOLs. Utilizing the collected data and taking into account the inherent uncertainties, a set of discrete scenarios are subsequently generated to develop a scenario-based two-stage stochastic programming approach for constructing the 3SCOPE model.

\subsubsection{First-stage Decisions} \label{fist_stage}

The \textit{first-stage} decisions, also referred to as here-and-now decisions, involve proactive decisions made before the uncertainties over PH become known. The first-stage decisions are made at the beginning of PH, marked as point \textit{A} in Figure \ref{flowchart}. This stage involves several key actions for the EM: 1) accepting pre-orders from customers and finalizing contracts with them; 2) selecting and finalizing contracts with suppliers for order placement and eVTOL part delivery; 3) selecting and finalizing contracts with logistics companies responsible for transporting eVTOL parts and eVTOLs at both ends of the supply chain; and 4) allocating resources for eVTOL manufacturing in the EM's plant. The 3SCOPE model formulates first-stage decisions encompassing several variables, including 1) number of eVTOLs to be manufactured in $F^P$ to satisfy the eVTOL pre-orders for $F^P$; 2) number of eVTOLs to be over-manufactured in $F^P$ to satisfy the eVTOL pre-orders for $S^P$; 3) timing of manufacturing these eVTOLs in $F^P$; 4) number of eVTOL parts to be ordered from suppliers in $F^P$; 5) number of eVTOL parts to be received from suppliers in $F^P$; 6) timing of the placing orders and receiving eVTOL parts in $F^P$; 7) number of eVTOL parts to be held in the EM's inventory in $F^P$; and 8) number of eVTOLs to be held in the EM's inventory in $F^P$. 
As a result of these first-stage decisions, the EM bears the total cost associated with its supply chain in $F^P$. This encompasses procurement cost, transportation cost at both ends of the supply chain, emission cost at both ends of the supply chain, manufacturing cost, inventory cost for eVTOL parts and manufactured eVTOLs, and penalty cost for delays in meeting customer orders. Collectively, these costs are referred to as first-stage total cost.

\subsubsection{Second-stage Decisions}

The \textit{second-stage} decisions, also referred to as wait-and-see decisions, represent reactive decisions made in recourse or response to compensate for the decisions made in the first-stage. The second-stage decisions allow the EM to make adjustments at point B, as depicted in the figure, in response to unfolding uncertainties in their supply chain planning. When the uncertainties in variables are revealed to the EM at the beginning of $S^P$, the EM would formulate second-stage decisions for $S^P$, which include 1) manufacturing schedule of the EM in $S^P$ to manufacture eVTOLs and satisfy the eVTOL pre-orders for this period; 2) timing of manufacturing these eVTOLs in $F^P$; 3) number of eVTOL parts to be ordered from suppliers in $S^P$; 4) number of eVTOL parts to be received from suppliers in $S^P$; 5) timing of the placing orders and receiving eVTOL parts in $S^P$; 6) number of eVTOL parts to be held in the EM's inventory in $S^P$; and 7) number of eVTOLs to be held in the EM's inventory in $S^P$. These decisions incur a second-stage total cost for $S^P$. The uncertainties are modeled through a set of discrete scenarios, as discussed next.




\subsubsection{Scenario Generation} \label{sc_gen}



To generate scenarios, we follow the scenario generation process presented in \cite{mohammadi2014scenario, bornapour2019efficient}. Initially, a normally-distributed probability density function (PDF) is defined for each uncertain parameter. The uncertainty of each parameter is modeled by incorporating an error value (EV), which represents the potential deviation from a expected value (EX) given for that parameter. This EV is added to EX of the parameter, and it is assumed to follow the PDF with a given standard deviation (SD). The value of uncertain parameter in a scenario would then be EX + EV $\times$ SD.

\begin{figure}[t] 
    \centering
    \includegraphics[width=0.7\textwidth]{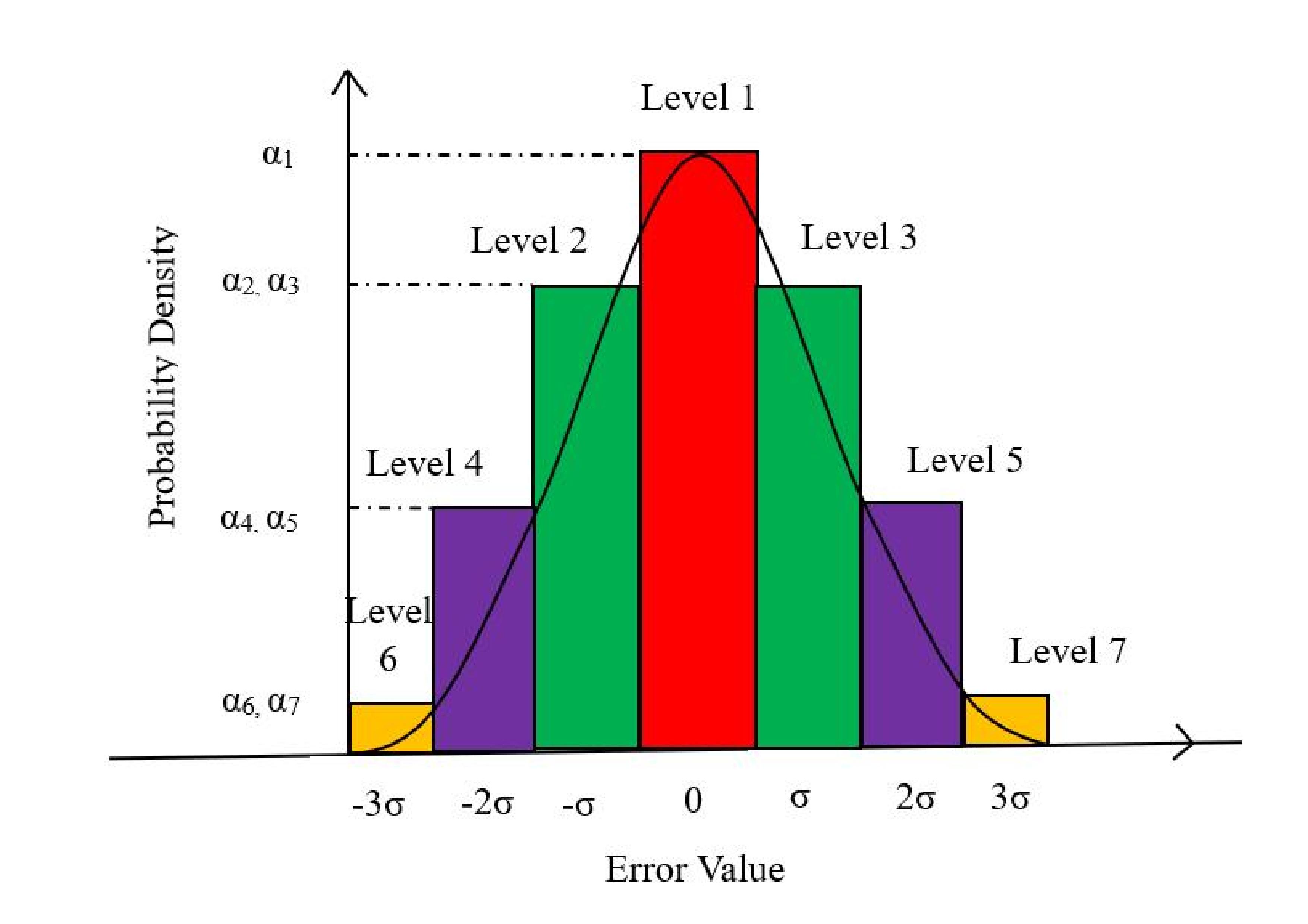}
    \caption{Discretization of the PDF defined for EV of each uncertain variable.}
    \label{fig:pdf}
\end{figure}

\begin{figure}[tb!] 
    \centering
    \includegraphics[width=0.7\textwidth]{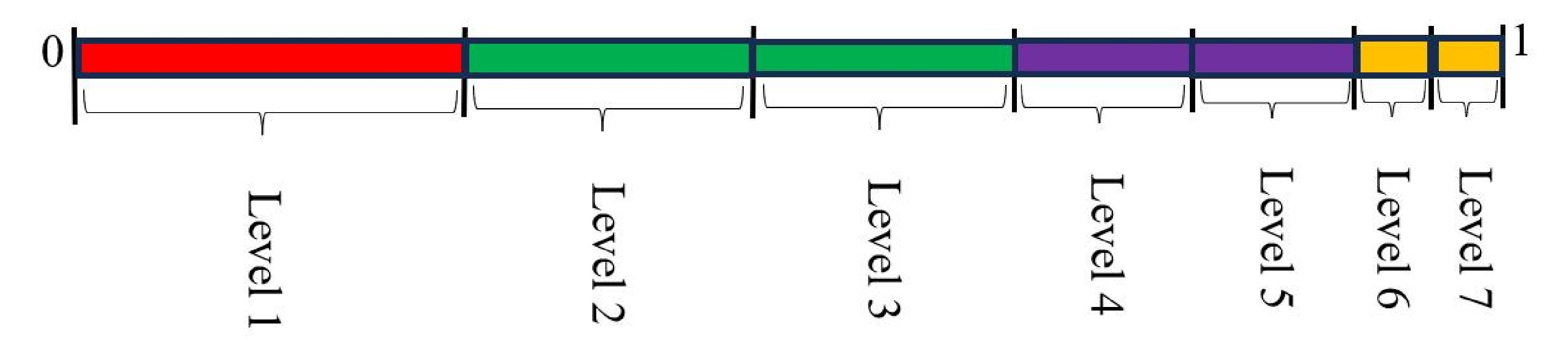}
    \caption{RWM for the normalized probabilities of the probability levels generated for each uncertain variable.}
    \label{fig:RWM}
\end{figure}

The PDF is discretized into seven probability levels, as depicted in Figure \ref{fig:pdf}. To obtain normalized probabilities for the probability levels, we employ a \textit{Roulette wheel mechanism (RWM)}, as illustrated in Figure \ref{fig:RWM}. In RWM, the difference between two adjacent probability levels equals the SD, which is associated with the PDF of the respective parameter. To generate scenarios in each iteration of RHA, initially, a random number within a range of [0, 1] is generated using RWM. Depending on the value of the random number drawn in a scenario, it determines the tier on the roulette wheel, representing a specific EV for the associated parameter in that scenario. This EV is then utilized to calculate the value of that uncertain parameter and its probability of occurrence in that scenario. This procedure is repeated for all uncertain parameters in that scenario. Finally, to compute the probability of the given scenario, all individual probabilities for the parameter values are multiplied together. Similarly, this process is iteratively repeated to generate all scenarios in this study.

\subsubsection{Scenario Reduction} \label{sc_reduc}

The scenario reduction technique involves approximating the scenario set using a smaller yet representative number of scenarios, providing a reasonably accurate approximation of the uncertainties. This approach is chosen to maintain a balance between computation time and solution quality. 
In this study, the primary strategy behind scenario reduction is twofold: removing similar scenarios and scenarios with very low probabilities. This scenario reduction strategy was also implemented in \cite{ahmadi2016novel, ahmadi2013risk}. The strategy involves excluding the most unlikely and repeated scenarios from the initial set. Consequently, the reduced set of scenarios, denoted by $N$, can be treated and solved as a deterministic problem, with each scenario associated with a specific probability.


\subsection{3SCOPE Model Formulation}

The 3SCOPE model is developed according to the parameters, indices, and decision variables listed in Table \ref{notation}.

\begin{spacing}{1}
\begin{longtable}{@{}p{3.5cm}p{12cm}@{}}
\caption{Parameters, indices, and decision variables in the SAND model.} \label{notation} \\
\toprule
\textbf{Parameters} & \textbf{Definition} \\ \midrule
\endfirsthead
$T^P$ & Set of periods in AH. \\
$V$ & Number of subsequent periods in PH: $V=2$. \\
$F^P$ & First period in PH: $F^P = P^H(1)$. \\
$S^P$ & Second period in PH except $F^P$: $S^P = P^H \backslash F^P$. \\
$P$ & Set of parts of an eVTOL. \\
$S_i$ & Set of suppliers for part $i \in P$. \\
$M$ & Set of transportation modes. \\
$C$ & Set of customers in any given period. \\
$T$ & Set of days in PH. \\
$N$ & Set of scenarios. \\
$\pi^s$ & Probability of scenario $s$. \\
$x_{ij}$ & First-stage procurement price of eVTOL part $i \in P$ when the eVTOL part is procured from supplier $j \in S_i$ at the beginning of $F^P$. \\
$x^s_{ij}$ & Second-stage procurement price of eVTOL part $i \in P$ when the eVTOL part is procured from supplier $j \in S_i$ at the beginning of $S^P$ in scenario $s \in N$.\\
$y_{ij}$ &  First-stage capacity of supplier $j \in S_i$ for eVTOL part $i \in P$ in $F^P$. \\
$y^s_{ij}$ &  Second-stage capacity of supplier $j \in S_i$ for eVTOL part $i \in P$ in scenario $s \in N$ in $S^P$. \\
$\lambda_{ijk}$ & Shipping distance between supplier $j \in S_i$ for eVTOL part $i \in P$ and the EM via transportation mode $k \in M$ in the back-end of the supply chain. \\
$T'_{ijk}$ & Transportation time for delivering eVTOL part $i \in P$ from supplier $j \in S_i$ to the EM via transportation mode $k \in M$ in the back-end of the supply chain. \\
$\gamma_{kt}$ & Average freight cost per ton-km transportation mode $k \in M$ on day $t \in T$. \\
$W_i$ & Weight of eVTOL part $i\in P$. \\
$\nu_k$ & Speed of transportation mode $k \in M$.\\
$\mu_k$ & Maximum daily driving hours for transportation mode $k \in M$. \\
$G_{ijk}$ & Emission of GHG per ton-mile generated for delivering eVTOL part $i \in P$ from supplier $j \in S_i$ to the EM via transportation mode $k \in M$ in the back-end of the supply chain. \\
$C^g_{t}$ & Cost per ton of GHG on a given day $t \in T$. \\
$\lambda_{lk}$ & Shipping distance between the EM and customer $l \in C$ via transportation mode $k \in M$ in the front-end of the supply chain. \\
$T'_{lk}$ & Transportation time for delivering eVTOLs from the EM to customer $l \in C$ via transportation mode $k \in M$ in the front-end of the supply chain. \\
$W^e$ &  Weight of an eVTOL. \\
$G_{lk} $ & Emission of GHG generated for delivering eVTOLs from the EM to customer $l \in C$ via transportation mode $k \in M$ in the front-end of the supply chain. \\
$L$ &  First-stage cost incurred by the EM for manufacturing an eVTOL. \\
$L^s$ &  Second-stage cost incurred by the EM for manufacturing an eVTOL in scenario $s \in N$. \\
$H^a_{i}$ & Holding cost of eVTOL part $i \in P$. \\
$H^b$ & Holding cost of an eVTOL on day. \\
$T^{d^1}_{l}$ &  Order-delivery deadline for customer $l \in C$ in $F^P$. \\
$T^{d^2}_{l}$ &  Order-delivery deadline for customer $l \in C$ in $S^P$. \\
$T^{d^\Delta}_{l}$ &  New order-delivery deadline for customer $l \in C$ in $F^P$ if PH is extended due to disruption. \\
$\Delta_l$ &  Penalty cost incurred by the EM for delaying the deadline for customer $l \in C$. \\
$\Upsilon^1_l$ &  Penalty cost incurred by the EM for delaying the deadline for customer $l \in C$ by one day from $T^{d^1}_{l}$. \\

$\eta$ & Selling price of an eVTOL. \\
$Q^r_i$ &  Minimum quality requirement for eVTOL part $i \in P$. \\ 
$Q_{ij}$ & Quality value of eVTOL part $i \in P$ of supplier $j \in S_i$ in $F^P$. \\ 
$ \Xi_{ij}$ & Normalized quality value of supplier $j \in S_i$ of part $i \in P$.\\ 
$\rho_{ij}$ &  Probability value associated with $ \Xi^s_{ij}$. \\
$\epsilon_{i}$ & Entropy of part $i \in P$. \\ 
$\psi_{i}$ & Degree of divergence of part $i \in P$. \\ 
$\kappa$ & Quality sensitivity parameter. \\
$\omega_{i}$ & Entropy weight of part $i \in P$. \\
$D^1_{l}$ & Final number of eVTOLs pre-ordered by customer $l \in C$ for $F^P$. \\
$D^2_{l}$ & Initial number of eVTOLs pre-ordered by customer $l \in C$ for $S^P$. \\
$D^s_{l}$ & Additional number of eVTOLs pre-ordered by customer $l \in C$ for $S^P$ in scenario $s$, subject to exogenous uncertainty arising from the market growth of eVTOLs. \\
$Q^b$ &  Minimum quality requirement for eVTOLs manufactured. \\ 
$R$ & Maximum number of eVTOLs that can be manufactured by the EM on any day. \\
$T''_{ij}$ & Lead time of supplier $j \in S_i$ of eVTOL part $i \in P$. \\
$\theta_i$ & Number of eVTOL parts $i \in P$ needed to manufacture an eVTOL. \\
$\alpha^1_{i}$ & Number of eVTOL parts $i \in P$ in the EM's inventory at the beginning of $F^P$. \\
$\sigma^\alpha_{i}$ & Daily capacity of the EM's inventory for holding eVTOL parts $i \in P$. \\
$\sigma^\beta$ &  Daily capacity of the EM's inventory for holding eVTOLs manufactured. \\
$\Upsilon^2$ & Minimum lead time of supplier across all eVTOL parts in $P$. \\
$\Upsilon^3$ & Minimum transportation time in the back-end of the supply chain across all transportation modes in $M$. \\
$\Upsilon^4$ & Minimum  transportation time for customer $l \in C$ in the front-end of the supply chain across all transportation modes in $M$. \\
$\Upsilon^5$ & Time required to manufacture $D^1_{l}$ eVTOLs for customer $l \in C$. \\
$\Upsilon^6$ & Minimum delay for customer $l \in C$. \\
$\Upsilon^7$ & Maximum delay for customer $l \in C$. \\
$\Upsilon^8$ & Maximum lead time of supplier across all eVTOL parts in $P$. \\
$\Upsilon^9$ & Maximum transportation time in the back-end of the supply chain across all transportation modes in $M$. \\
$\Upsilon^{10}$ & Maximum transportation time for customer $l \in C$ in the front-end of the supply chain across all transportation modes in $M$. \\

\midrule

\textbf{Indices} &  \\ \midrule
$i$ & $i$-th part in $P$. \\
$j$ & $j$-th supplier in $S_i$. \\
$k$ & $k$-th transportation mode in $M$. \\
$l$ & $l$-th customer in $C_p$. \\
$p$ & $p$-th period of PH. \\
$t$, $t_1$, $t_2$ & $t$-th day, $t_1$-th day, $t_2$-th day in $T$, respectively. \\
$s$ & $s$-th scenario in $N$.\\

\midrule

\textbf{Decision Variables} &  \\  \midrule

$f^1$ & Continuous variable representing first-stage operating profit to be generated in $F^P$. \\

$f^2$ & Continuous variable representing second-stage operating profit to be generated in $S^P$. \\

$\delta$ & Continuous variable representing first-stage revenue to be generated in $F^P$.  \\

$\delta^s$ & Continuous variable representing second-stage revenue to be generated in $S^P$ in scenario $s \in N$. \\

$\Omega$ & Continuous variable representing first-stage total cost to be incurred in $F^P$. \\

$\Omega^s$ & Continuous variable representing second-stage total cost to be incurred in $S^P$ in scenario $s \in N$. \\

$\Omega^\Upsilon$ & Continuous variable representing first-stage total cost to be incurred in $F^P$ if PH is extended due to disruption. \\

$v_{t i j k}$ & Non-negative integer variable representing first-stage number of eVTOL parts $i \in P$ to be ordered to supplier $j \in S_i$ at the beginning of a given day $t \in T$ in $F^P$  when transportation mode $k \in M$ is chosen in the back-end of the supply chain.  \\

$v^s_{t i j k}$ & Non-negative integer variable representing second-stage number of eVTOL parts $i \in P$ to be ordered to supplier $j \in S_i$ at the beginning of a given day $t \in T$ in $S^P$ in scenario $s \in N$ when transportation mode $k \in M$ is chosen in the back-end of the supply chain. \\

$z_{ t i j k}$ &  Non-negative integer variable representing first-stage number of eVTOL parts $i \in P$ to be received from supplier $j \in S_i$ via transportation mode $k \in M$ at the beginning of a given day $t \in T$ in $F^P$. \\

$z^s_{ t i j k}$ & Non-negative integer variable representing second-stage number of eVTOL parts $i \in P$ to be received from supplier $j \in S_i$ via transportation mode $k \in M$ at the beginning of a given day $t \in T$ in $S^P$ in scenario $s \in N$. \\

$r_{tlk} $ & Non-negative integer variable representing first-stage number of eVTOLs to be manufactured at the beginning of a given day $t \in T$ in $F^P$ for customer $l \in C$ when the mode in the front-end of the supply chain is $k \in M$.\\

$\varphi_{l}$ & Non-negative integer variable representing first-stage number of eVTOLs to be over-manufactured in $F^P$ for customer $l \in C$.\\

$r^s_{tlk} $ &  Non-negative integer variable representing second-stage number of eVTOLs to be manufactured at the beginning of a given day $t \in T$ in $S^P$ for customer $l \in C$ in scenario $s \in N$ when the mode in the front-end of the supply chain is $k \in M$.\\

$\alpha^2_{i}$ & Number of eVTOL parts $i \in P$ to be held in the EM's inventory at the beginning of $S^P$. \\

$\alpha_{ti}$ & Non-negative integer variable representing first-stage number of eVTOL parts $i \in P$ to be held in the EM's inventory at the beginning of a given day $t \in T$ in $F^P$. \\

$\alpha^s_{ti}$ & Non-negative integer variable representing second-stage number of eVTOL parts $i \in P$ to be held in the EM's inventory at the beginning of a given day $t \in T$ in $S^P$ in scenario $s \in N$. \\

$\beta_{tlk}$ & Non-negative integer variable representing first-stage number of eVTOLs to be held in the EM's inventory for customer $l \in C$ at the beginning of a given day $t \in T$ in $F^P$ when the delivery transportation mode in the front-end of the supply chain is $k \in M$.\\

$\beta^s_{tlk}$ &  Non-negative integer variable representing second-stage number of eVTOLs to be held in the EM's inventory for customer $l \in C$ at the beginning of a given day $t \in T$ in $S^P$ in scenario $s \in N$ when the delivery transportation mode in the front-end of the supply chain is $k \in M$. \\

$q$ & Continuous variable representing average quality value of eVTOLs across the eVTOL-batch manufactured in $F^P$. \\

$d^2_{l}$ & Non-negative integer endogenous variable representing additional number of eVTOLs pre-ordered by customer $l \in C$, driven by $q$, for $S^P$. \\


\bottomrule
\end{longtable} 
\end{spacing}

\subsubsection{Objective of 3SCOPE Model}

The objective of the 3SCOPE model is to maximize the operating profit of the EM by optimizing supply chain planning decisions over AH, implementing RHA. The objective function for each iteration of RHA is represented by Eq. \ref{obj}, where $f^1$ represents the first-stage operating profit to be generated in $F^P$, and $f^2$ represents the second-stage operating profit to be generated in $S^P$. These $f^1$ and $f^2$ are calculated based on the revenues generated and the costs incurred over the supply chain, following Eqs. \ref{profit1} and \ref{profit2}, respectively. $\delta$ denotes the first-stage revenue to be generated in $F^P$, and $\delta^s$ refers to the second-stage revenue to be generated in $S^P$ in scenario $s$. These revenues are generated from the sale of eVTOLs manufactured for customers, as detailed in Eqs. \ref{delta1} and \ref{delta2}. Here, the variable $\eta$ represents the selling price of an eVTOL, and $D^s_{l}$ denotes the additional number of eVTOLs pre-ordered by customer $l$ for $S^P$, driven by market growth in scenario $s$. In Eq. \ref{profit1}, $\Omega$ represents the first-stage total cost to be incurred in $F^P$. $\Omega$ comprises eight distinct total cost components that are anticipated to be incurred after the contracts with customers, suppliers, logistics companies and other entities at the beginning of $F^P$. These components, as presented in Eq. \ref{omega}, are as follows: 1) total eVTOL parts procurement cost, 2) total transportation cost in the back-end of the supply chain, 3) total transportation cost in the front-end of the supply chain, 4) total emission cost in the back-end of the supply chain, 5) total emission cost in the front-end of the supply chain, 6) total eVTOL manufacturing cost, 7) total inventory cost for holding eVTOL parts,  and 8) total inventory cost for holding manufactured eVTOLs. In Eq. \ref{profit2}, $\Omega^s$ refers to the second-stage total cost to be incurred in $S^P$, and $\pi^s$ denotes the probability of scenario $s \in N$, created to capture the uncertainties revealed in the second stage. $\Omega^s$, as defined in Eq. \ref{omega_s}, encompasses the same eight distinct total cost components as those outlined for $\Omega$. These costs are expected to be incurred in $S^P$ in scenario $s$ when the uncertainties are anticipated to be revealed.

\begin{spacing}{0.6}

\begin{equation} \label{obj}
\text{maximize:}  f^1 +  f^2 
\end{equation}

\begin{equation} \label{profit1}
 f^1 = \delta - \Omega
\end{equation}

\begin{equation} \label{profit2}
 f^2 = \sum_{s \in  N} {\pi^s (\delta^s - \Omega^s)} 
\end{equation}

\begin{equation} \label{delta1}
\delta =   \sum_{l \in C} \eta D^1_{l}
\end{equation}

\begin{equation} \label{delta2}
\delta^s =  \sum_{l \in C} \eta (d^2_{l} + D^2_{l} + D^s_{l}),  \quad \forall s \in N
\end{equation}

\begin{equation} \label{omega}
\begin{split}
    \Omega &= \sum_{t=1}^{\frac{|T|}{2}} \sum_{i \in P} \sum_{j \in S_i} \sum_{k \in M} x_{ij} z_{tijk}  +  \sum_{t=1}^{\frac{|T|}{2}} \sum_{i \in P} \sum_{j \in S_i} \sum_{k \in M} \gamma_{kt} W_i \lambda_{ijk} z_{tijk}  +  \sum_{l \in C} \sum_{k \in M} \sum_{t=1}^{\frac{|T|}{2}} \gamma_{kt} W^e \lambda_{lk} r_{tlk}  \\
       & \quad +  \sum_{t=1}^{\frac{|T|}{2}} \sum_{i \in P} \sum_{j \in S_i} \sum_{k \in M} C^t_{t} G_{ijk} W_i z_{tijk} +
     \sum_{l \in C} \sum_{k \in M} \sum_{t=1}^{\frac{|T|}{2}} C^t_{t} G_{lk} W^e r_{tlk} + \sum_{l \in C} \sum_{k \in M} \sum_{t=1}^{\frac{|T|}{2}}  L r_{tlk}   \\
    & \quad  + \sum_{t =1}^{\frac{|T|}{2} - 1} \sum_{i \in P} \alpha_{ti} H^a_{i} +  \sum_{t=1}^{T^{d^1}_{l}-T'_{lk}} \sum_{l \in C} \sum_{k \in M} \beta_{tlk} H^b
\end{split}
\end{equation}

\begin{equation} \label{omega_s}
\begin{split}
    \Omega^s &= \sum_{t= \frac{|T|}{2}+1}^{|T|} \sum_{i \in P} \sum_{j \in S_i} \sum_{k \in M} x^s_{ij} z^s_{tijk}  + \sum_{t= \frac{|T|}{2}+1}^{|T|} \sum_{i \in P} \sum_{j \in S_i} \sum_{k \in M} \gamma_{kt} W_i \lambda_{ijk} z^s_{tijk} +  \sum_{l \in C} \sum_{k \in M} \sum_{t= \frac{|T|}{2}+1}^{|T|} \gamma_{kt} W^e \lambda_{lk} r^s_{tlk}  \\
    & \quad  + \sum_{t= \frac{|T|}{2}+1}^{|T|} \sum_{i \in P} \sum_{j \in S_i} \sum_{k \in M} C^t_{t} G_{ijk} W_i z^s_{tijk} +  \sum_{l \in C} \sum_{k \in M} \sum_{t= \frac{|T|}{2}+1}^{|T|} C^t_{t} G_{lk} W^e r^s_{tlk}  + \sum_{l \in C} \sum_{k \in M} \sum_{t= \frac{|T|}{2}+1}^{|T|} L^s r^s_{tlk}   \\
    & \quad + \sum_{t = \frac{|T|}{2}+1}^{|T| - 1} \sum_{i \in P} \alpha^s_{ti} H^a_{i}  + \sum_{t= \frac{|T|}{2}+1}^{T^{d^2}_{l}-T'_{lk}} \sum_{l \in C} \sum_{k \in M} \beta^s_{tlk} H^b, \quad \forall s \in  N
\end{split}
\end{equation}

\end{spacing}
\indent

Next, we outline the calculations for the eight distinct total cost components. In the first stage of the 3SCOPE model, when calculating the total eVTOL parts procurement cost (the first component in Eq. \ref{omega}), we take into account the procurement price denoted by $x_{ij}$ of eVTOL part $i$ procured from supplier $j \in S_i$ at the beginning of $F^P$. Additionally, we consider $z_{tijk}$, representing the number of eVTOL parts $i$ to be received from supplier $j$ via transportation mode $k$ at the beginning of a given day $t$ in $F^P$. We then analyze the total transportation costs at both ends of the supply chain, represented by the second and third cost components. As discussed in Section \ref{logistic_contract}, the logistics companies would charge the EM based on the the shipping distance and the total weight of eVTOL parts and eVTOLs transported. The second cost component in Eq. \ref{omega} is the total transportation cost to be incurred in the back-end of the supply chain in $F^P$, which is driven by the decision variable, $z_{tijk}$, and three parameters. The parameters are 1) the shipping distance (km) between supplier $j$ for part $i$ and the EM via transportation mode $k$, denoted by $\lambda_{ijk}$; 2) the average freight cost per ton-km transportation mode $k$ on a given day $t$, represented by $\gamma_{kt}$; and 3) the weight of part $i$, denoted by $W_i$, to determine the total weight of $z_{tijk}$.  By considering the shipping distance $\lambda_{ijk}$, we calculate $T'_{ijk}$ using Eq. \ref{transtime1}. $T'_{ijk}$ refers to the transportation time for delivering eVTOL part $i$ from supplier $j$ to the EM via transportation mode $k$ in the back-end of the supply chain. To calculate $T'_{ijk}$, we consider the constraints imposed by the speed of transportation mode $k$, denoted by $\nu_k$, and the maximum daily driving hours for that mode, denoted by $\mu_k$. The third cost component in Eq. \ref{omega} is the total transportation cost to be incurred in the front-end of the supply chain in $F^P$. The variable $r_{tlk}$ represents the number of eVTOLs to be manufactured at the beginning of a given day $t$ in $F^P$ for customer $l$ when transportation mode $k$ is chosen. The logistics company of transportation mode $k$ would charge the EM based on the weight of $r_{tlk}$ eVTOLs that need to be transported to deliver to customer $l$. The shipping distance between the EM and customer $l$ via transportation mode $k$ is denoted by $\lambda_{lk}$, and the weight of an eVTOL is denoted by $W^e$. Similar to the calculation of transportation time in the back-end of the supply chain, $T'_{lk}$, the time required to transport eVTOLs from the EM to customer $l$ via transportation mode $k$, is calculated using Eq. \ref{transtime2}.

\begin{spacing}{0.6}
    
\begin{equation}\label{transtime1}
T'_{ijk}= \frac{\lambda_{ijk}}{\nu_k \mu_k}, \quad \forall i \in P, \forall j \in S_i, \forall k \in M
\end{equation}

\begin{equation}\label{transtime2}
T'_{lk}= \frac{\lambda_{lk}}{\nu_k \mu_k}, \quad \forall l \in  C, \forall k \in M
\end{equation}

\end{spacing}
\indent 




Furthermore, to account for the associated total emission costs during the transportation, we consider GHG emission for the back-end of the supply chain. In the fourth cost component in Eq. \ref{omega}, the cost per ton of GHG emission on a given day $t$ is referred to as $C^g_{t}$, while $G_{ijk}$ indicates the tons of GHG emissions generated during the transportation of each ton of eVTOL part $i$ from supplier $j$ via transportation mode $k$ over the shipping distance $\lambda_{ijk}$ in the back-end of the supply chain. Similarly, the fifth cost component in Eq. \ref{omega} refers to the GHG emission in the front-end of the supply chain. Next, the sixth cost component in Eq. \ref{omega} refers to the total eVTOL manufacturing cost to be incurred in $F^P$, where $L$ represents the cost incurred by the EM for manufacturing an eVTOL. Another significant cost that affects the profitability is inventory cost, which accounts for the expenses to be incurred in holding inventory over a given period. We consider inventory costs at both ends of the supply chain in $F^P$: 1) total inventory cost for holding eVTOL parts received from the suppliers in the back-end of the supply chain (the seventh cost component in Eq. \ref{omega}), and 2) total inventory cost for holding eVTOLs manufactured in the front-end of the supply chain before delivering them to the customers (the eighth cost component in Eq. \ref{omega}). In the seventh cost component, the first-stage variable $\alpha_{ti}$ denotes the number of eVTOL parts $i$ in inventory at the beginning of a given day $t$ in $F^P$. $H^a_i$ refers to the holding cost of eVTOL part $i$ in the inventory, reflecting the cost associated with holding one unit of that particular part. Moving on to the eighth cost component, $H^b$ refers to the cost for holding an eVTOL in the inventory. $\beta_{tlk}$ denotes the first-stage number of eVTOLs manufactured in inventory held for customer $l$ at the beginning of a given day $t$ in $F^P$, when transportation mode $k$ is chosen in the front-end of the supply chain. The order-delivery deadline for customer $l$ in $F^P$ is represented by ${T^{d^1}_{l}}$. The eVTOLs manufactured should be dispatched from the inventory for delivery by $(T^{d^1}_{l} - T'_{lk})^{th}$ day, considering it would need $T'_{lk}$ days for transportation to reach customer $l$ by ${T^{d^1}_{l}}$ if transportation mode $k$ is chosen. The inventory cost for holding manufactured eVTOLs would be incurred for the days that the eVTOLs remain in inventory after manufacturing, prior to their dispatch for delivery. 

Similarly, the description holds true for the cost components presented in Eq. \ref{omega_s} corresponding to the second-stage of the 3SCOPE model. The second-stage variables are as follows: $x^s_{ij}$ represents the procurement price of eVTOL part $i$ procured from supplier $j$ at the beginning of $S^P$ in scenario $s$, $z^s_{tijk}$ the number of eVTOL parts $i$ to be received from supplier $j$ via transportation mode $k$ at the beginning of a given day $t$ in $S^P$ within scenario $s$, $r^s_{tlk}$ the number of eVTOLs to be manufactured at the beginning of a given day $t$ in $S^P$ for customer $l$ when transportation mode $k$ is selected in scenario $s$, $L^s$ the cost incurred by the EM for manufacturing an eVTOL in $S^P$ in scenario $s$, $\alpha^s_{ti}$ the number of eVTOL parts $i$ to be held in inventory at the beginning of a given day $t$ in $S^P$ in scenario $s$, and $\beta^s_{tlk}$ the number of eVTOLs manufactured to be held in inventory allocated for customer $l$ at the beginning of a given day $t$ in $S^P$ in scenario $s$, when transportation mode $k$ is chosen in the front-end of the supply chain.


\subsubsection{Constraints of 3SCOPE model}

The 3SCOPE model addresses practical limitations imposed by real-world conditions on the variables influencing the objective function. This is achieved by incorporating various constraints at both stages of its formulation, which are discussed next.

\makeatletter
\renewcommand\paragraph{\@startsection{paragraph}{5}{\z@}%
                                     {-3.25ex\@plus -1ex \@minus -.2ex}%
                                     {1.5ex \@plus .2ex}%
                                     {\normalfont\normalsize\bfseries}}
\renewcommand\theparagraph{\arabic{paragraph}} 

\renewcommand\theparagraph{\arabic{section}.\arabic{paragraph}}

\paragraph{Quality Constraints} \label{quality_cons}



In the aviation industry, it is crucial to consider quality constraints to ensure that the eVTOLs manufactured meet the minimum quality requirements set by the FAA, as discussed in Section \ref{sec2}. Within the 3SCOPE model, we incorporate a set of quality constraints to achieve these standards. The first quality constraint is represented by Eq. \ref{eqc1}, which ensures that the quality value for eVTOL part $i$ of supplier $j$, denoted as $Q_{ij}$, meets the minimum quality requirement $Q^r_i$ for eVTOL part $i$. The EM would contact suppliers known for producing high-quality eVTOL parts and would assess the quality of these eVTOL parts before procurement. The quality of these eVTOL parts would eventually impact the final quality of the entire eVTOL batch manufactured in $F^P$. When evaluating the quality of different parts used in the eVTOL manufacturing process, \textit{Entropy Weights Method (EWM)} can be implemented to determine the weights assigned to each part. The EWM is a decision-making technique widely used to assign weights to different parts of a product based on their relative importance \cite{EWM, EWM1}. These weights would reflect the contribution of each part to the quality of a manufactured eVTOL. By considering the entropy and divergence of values across the different eVTOL parts, EWM identifies eVTOL parts with lower entropy, which receive higher weights indicating their higher importance in contributing to the quality of the final product \cite{EWM2, EWM3}. In EWM, the first step, as presented in Eq. \ref{EWMeq1}, is to normalize the quality values by dividing them by the maximum quality value among all eVTOL parts, $max(Q_{ij})$. The normalization ensures that the quality values are scaled within the range of 0 to 1. The other steps are presented in Eqs. \ref{EWMeq2} to \ref{EWMeq5}, following the steps mentioned in \cite{EWM, EWM1}. In these equations, $\Xi_{ij}$ represents the normalized quality value of eVTOL part $i$ from supplier $j$, $\rho_{ij}$ the probability value associated with each normalized quality value $\Xi_{ij}$, $\epsilon_{i}$ the entropy value for eVTOL part $i$, and $\psi_{i}$ the degree of divergence of eVTOL part $i$, and $\omega_{i}$ the relative importance of the quality of eVTOL part $i$ in the manufactured eVTOLs.

\begin{spacing}{0.5}
\begin{equation} \label{eqc1}
Q_{ij} \geq Q^r_i, \quad   \forall i \in  P,  \forall j \in S_i
\end{equation}

\begin{equation} \label{EWMeq1}
 \Xi_{ij} = \frac{Q_{ij}}{max(Q_{ij})}, \quad   \forall i \in  P,  \forall j \in S_i
\end{equation}

\begin{equation}\label{EWMeq2}
\rho_{ij} = \frac{\Xi_{ij}}{\sum_{{j \in S_i}} \Xi_{ij}}, \quad  \forall i \in  P,  \forall j \in  S_i 
\end{equation}

\begin{equation} \label{EWMeq3}
\epsilon_{i} = -\frac{1}{\log(|P|)} \sum_{{j \in S_i}} (\rho_{ij} \log(\rho_{ij})), \quad  \forall i \in  P
\end{equation}

\begin{equation} \label{EWMeq4}
\psi_{i} = |1-\epsilon_{i}|, \quad  \forall i \in  P
\end{equation}

\begin{equation} \label{EWMeq5}
\omega_{i} = \frac{\psi_{i}}{\sum_{i \in P} \psi_{i}}, \quad  \forall i \in  P
\end{equation}

\end{spacing}
\indent 


Using the weights derived from EWM, we select \textit{Weighted Averaging (WA)} approach to determine $q$, which represents the average quality value of eVTOLs across the entire eVTOL-batch manufactured in $F^P$. The WA is a commonly used approach to determine the average quality value of products and services in various fields \cite{boran2009multi, w1, w}. We compute $q$ following Eq. \ref{WA}, where in the denominator, the total number of eVTOL parts to be received and used in manufacturing eVTOLs during $F^P$ are weighted according to their respective weights. The numerator represents the total weighted quality values of those eVTOL parts to be received from suppliers in $F^P$. The second quality constraint is the minimum FAA quality requirement constraint for the eVTOL-batch manufactured, as presented in Eq. \ref{eqc2}. This constraint ensures that $q$ cannot be less than the base quality value, denoted as $Q^b$, which is considered to be the minimum FAA quality requirement for eVTOLs.


\begin{equation} \label{WA}
q = \frac{\sum_{t=1}^{\frac{|T|}{2}} \sum_{i \in P} ( \omega_{i} \sum_{j \in S_i} (Q_{ij} \sum_{k \in M} z_{ t i j k}))}{ \sum_{t=1}^{\frac{|T|}{2}} \sum_{i \in P} ( \omega_{i} \sum_{j \in S_i} \sum_{k \in M} z_{ t i j k})}
\end{equation}

\begin{equation} \label{eqc2}
q > Q^b
\end{equation}


\paragraph{Quality and Customer Order Constraint}  \label{quality_o}


The decision variable $q$ influences the additional number of eVTOLs pre-ordered by customer $l$ for $S^P$. If $q$ exceeds $Q^b$, customer pre-orders for $S^P$ are expected to increase, as discussed in Section \ref{customer_pre}. We formulate first-stage quality and customer order constraint following Eq. \ref{eq:demand2}, which incorporates the relationship between $q$ and $D^1_{l}$. Such a relationship can be established between quality of a product and demand of customers for that product, as mentioned in Section \ref{customer_pre}. $d^2_{l}$ represents the portion of the additional number of eVTOLs pre-ordered by customer $l$ for $S^P$ that is affected by $q$. The quality and customer order constraint considers an adjustment term $(q - Q^b)$. This adjustment term is multiplied by $D^1_{l}$ and quality sensitivity parameter $\kappa$, which determines the sensitivity of $d^2_{l}$ to changes in $q$. The constraint captures the change in customer pre-orders by considering how the deviation of $q$ from $Q^b$ influences customer preferences. If $q$ is higher than $Q^b$, the adjustment term will be positive, indicating that the EM would expect to receive $d^2_{l}$ from customer $l$ for $S^P$. However, if $q$ is the same as $Q^b$, there is no additional eVTOL pre-orders from the customers driven by $q$.

\begin{equation} \label{eq:demand2}
d^2_{l} =  (q - Q^b) \kappa D^1_{l}, \quad  \forall l \in C
\end{equation}

\paragraph{Manufacturing Rate and Customer Order Fulfillment Constraints}

To meet the order requirements of customer $l$ within the specified deadline $T^{d^1}_{l}$, Eq. \ref{manu_rate_cons} ensures that the total number of eVTOLs manufactured for customer $l$ from day 1 to day $(T^{d^1}_{l}-T'_{lk})$ in $F^P$ must be equal to at least $D^1_{l}$ when transportation mode $k$ is chosen. In addition to manufacturing eVTOLs to fulfill all the eVTOL orders of $F^P$, the EM would also manufacture $\varphi_{l}$ eVTOLs in $F^P$, which is computed following Eq. \ref{manu_rate_cons1}. $\varphi_{l}$ represents the first-stage number of eVTOLs to be over-manufactured for customer $l$. These additional eVTOLs are produced in advance to fulfill eVTOL orders for $S^P$, accounting for uncertainties in the EM's supply chain, as mentioned in Section \ref{prob_stat}. The second-stage order fulfillment constraint in Eq. \ref{manu_rate_cons2} ensures that the total number of eVTOLs manufactured for customer $l$ from day $\frac{|T|}{2}+1$ to day $(T^{d^2}_{l}-T'_{lk})$ in $S^P$ must be equal to the remaining number of eVTOLs pre-ordered $(d^2_{l} + D^2_{l} + D^s_{l} - \varphi_l)$ by that customer for $S^P$ in scenario $s$. The manufacturing capacity constraints are presented in Eq. \ref{manu_rate_cons3} for the first-stage  and in Eq. \ref{manu_rate_cons4} for the second-stage.

\begin{spacing}{0.5}
\begin{equation} \label{manu_rate_cons}
  \sum_{{k \in M}} \sum_{{t=1}}^{T^{d^1}_{l}-T'_{lk}}   r_{tlk}  \geq   D^1_{l} , \quad  \forall l \in C
\end{equation}

\begin{equation}  \label{manu_rate_cons1}
\varphi_l =  \sum_{{k \in M}} \sum_{{t=1}}^{\frac{|T|}{2}} r_{tlk} - D^1_{l}, \quad  \forall l \in C
\end{equation}

\begin{equation} \label{manu_rate_cons2}
  \sum_{{k \in M}} \sum_{{t=\frac{|T|}{2}+1}}^{T^{d^2}_{l}-T'_{lk}}  r^s_{tlk} = d^2_{l} + D^2_{l} + D^s_{l} - \varphi_l, \quad \forall l \in C,  \forall s \in N
\end{equation}

\begin{equation} \label{manu_rate_cons3}
 \sum_{l \in C}  \sum_{{k \in M}}  r_{tlk} \leq   R, \quad \forall t \in \{1, 2, ..., \frac{|T|}{2}\}
\end{equation}

\begin{equation} \label{manu_rate_cons4}
 \sum_{l \in C}  \sum_{{k \in M}}  r^s_{tlk}  \leq   R, \quad \forall t \in \{\frac{|T|}{2}+1, \frac{|T|}{2}+2, ..., |T|\},  \forall s \in N
\end{equation}

\end{spacing}
\indent





\paragraph{Inventory Management Constraints}


The first-stage inventory part-flow constraint for eVTOL part $i$ in Eq. \ref{inven_manage1} ensures a proper balance in the inventory level in the back-end of the supply chain. $\alpha_{ti}$ represents the number of eVTOL parts $i$ to be held in the EM's inventory at the beginning of a given day $t$ in $F^P$. $\alpha_{ti}$ is determined by 1) the initial inventory for eVTOL part $i$ at the beginning of $F^P$ ($\alpha^1_i$); 2) the total number of eVTOL parts $i$ received from suppliers from the beginning of $F^P$ up to day $t$ ($\sum_{{t_2=1}}^{t} \sum_{{j \in S_i}} \sum_{{k \in M}} z_{t_2ijk}$); and 3) the total number of eVTOL parts $i$ used in the manufacturing of eVTOLs from the beginning of $F^P$ up to day $t$ ($\sum_{{t_2=1}}^{t} \sum_{{l \in C}} \sum_{{k \in M}} \theta_i r_{t_2lk}$). The parameter $\theta_i$ represents the number of eVTOL parts $i$ required to manufacture an eVTOL. At $t = \frac{|T|}{2}$, $\alpha_{ti}$ in Eq. \ref{inven_manage2} represents the number of eVTOL parts $i$ to be held in the EM's inventory at the beginning of $S^P$, denoted by $\alpha^2_i$. The second-stage inventory part-flow constraint is presented in Eq. \ref{inven_manage3} for eVTOL part $i$, where $\alpha^s_{ti}$ represents the number of eVTOL parts $i$ to be held in the EM's inventory at the beginning of a given day $t$ in $S^P$ in scenario $s$. Next, the first-stage inventory availability constraint in Eq. \ref{inven_manage4} states that in $F^P$, $\alpha_{ti}$ must be greater than or equal to the number of eVTOL parts $i$ required for the manufacturing of eVTOLs on the next day. There should be a sufficient number of eVTOL parts $i$ available in the inventory so that the EM can assemble the eVTOLs required for the next day's operations. Similarly, the second-stage inventory availability constraint is presented in Eq. \ref{inven_manage5}. Then the first-stage inventory capacity constraint in Eq. \ref{inven_manage6} ensures that $\alpha_{ti}$ does not exceed the EM's inventory capacity $\sigma^\alpha_{i}$ for eVTOL part $i$. Similarly, the second-stage inventory capacity constraint is presented in Eq. \ref{inven_manage7}, which ensures that $\alpha^s_{ti}$  in scenario $s$ does not exceed the EM's inventory capacity $\sigma^\alpha_{i}$ for eVTOL part $i$.

\begin{spacing}{0.5}

\begin{align} \label{inven_manage1}
\alpha_{ti} &= \alpha^1_i + \sum_{{t_2=1}}^{t} \sum_{{j \in S_i}} \sum_{{k \in M}} z_{t_2ijk} - \sum_{{t_2=1}}^{t} \sum_{{l \in C}} \sum_{{k \in M}} \theta_i r_{t_2lk}, \quad \forall t \in \{1, 2, ..., \frac{|T|}{2}\}, \forall i \in P
\end{align}

\begin{align} \label{inven_manage2}
\alpha^2_i = \alpha_{ti}, \quad \forall t = \frac{|T|}{2}, \forall i \in P
\end{align}

\begin{align} \label{inven_manage3}
\alpha^s_{ti} &= \alpha^2_i + \sum_{{t_2=1}}^{t} \sum_{{j \in S_i}} \sum_{{k \in M}}  z^s_{t_2ijk} - \sum_{{t_2=1}}^{t} \sum_{{l \in C}} \sum_{{k \in M}} \theta_i r^s_{t_2lk}, \quad \forall t \in \{\frac{|T|}{2}+1, \frac{|T|}{2}+2, ..., |T|\}, \forall i \in P, \forall s \in N
\end{align}

\begin{equation} \label{inven_manage4}
\alpha_{ti} \geq \sum_{{l \in C}} \sum_{{k \in M}}  \theta_i  r_{\left(t+1 \right)lk}, \quad  \forall t \in \{1, 2, ..., \frac{|T|}{2}\}, \forall i \in P
\end{equation}

\begin{equation} \label{inven_manage5}
\alpha^s_{ti} \geq \sum_{{l \in C}} \sum_{{k \in M}}  \theta_i   r^s_{\left(t+1 \right)lk}, \quad  \forall t \in \{\frac{|T|}{2}+1, \frac{|T|}{2}+2, ..., |T|\}, \forall i \in P, \forall s \in  N
\end{equation}

\begin{equation} \label{inven_manage6}
 \alpha_{ti} \leq \sigma^\alpha_{i}, \quad  \forall  t \in \{1, 2, ..., |T|\}, \forall i \in P
\end{equation}

\begin{equation} \label{inven_manage7}
 \alpha^s_{ti} \leq \sigma^\alpha_{i}, \quad  \forall  t \in \{1, 2, ..., |T|\}, \forall i \in P, \forall s \in N
\end{equation}

\end{spacing}

\indent




The eVTOLs manufactured by the EM at the beginning of a given day $t$ in $F^P$ for customer $l$ are held in the front-end inventory of the supply chain. When the chosen transportation mode for customer $l$ in the front-end of the supply chain is $k$, and the delivery time is $T'_{lk}$, the eVTOLs remain in the inventory from day $t$ until day $(T^d_{l}-T'_{lk})$. The first-stage inventory constraint for holding eVTOLs manufactured is given in Eq. \ref{inven_manage8}, where $\beta_{tlk}$ represents the number of eVTOLs to be held in the EM's inventory for customer $l$ at the beginning of a given day $t$ in $F^P$ when the delivery transportation mode in the front-end of the supply chain is $k$. Similarly, Eq. \ref{inven_manage9} is for the second-stage inventory constraint for holding eVTOLs manufactured, where $\beta^s_{tlk}$ represents the number of eVTOLs to be held in the EM's inventory for customer $l$ at the beginning of a given day $t$ in $S^P$ in scenario $s$ when the delivery transportation mode in the front-end of the supply chain is $k$. Then Eq. \ref{inven_manage10} ensures that $\beta_{tlk}$ does not exceed the capacity $\sigma^\beta$ of the inventory for eVTOLs manufactured in $S^P$. Similarly, Eq. \ref{inven_manage11} ensures that $\beta^s_{tlk}$ does not exceed the capacity $\sigma^\beta$ of the inventory for eVTOLs manufactured in $S^P$ in scenario $s$.

\begin{spacing}{0.6}
    
\begin{equation} \label{inven_manage8}
\beta_{tlk} = (T^{d^1}_{l}-T'_{lk}-t)  r_{tlk}, \quad  \forall l \in C, \forall k \in M, \forall t \in \{1, 2, ..., (T^{d^1}_{l}-T'_{lk})\}
\end{equation}

\begin{equation} \label{inven_manage9}
\beta^s_{tlk} = (T^{d^2}_{l}-T'_{lk}-t) r^s_{tlk} , \quad \forall l \in C, \forall k \in M, \forall t \in \{\frac{|T|}{2}+1, \frac{|T|}{2}+2, ..., (T^{d^2}_{l}-T'_{lk})\}, \forall s \in N
\end{equation}

\begin{equation} \label{inven_manage10}
\sum_{{l \in C}} \sum_{{k \in M}}  \beta_{tlk} \leq  \sigma^\beta , \quad  \forall t \in \{1, 2, ..., (T^{d^1}_{l}-T'_{lk})\}, \forall s \in N
\end{equation}

\begin{equation} \label{inven_manage11}
\sum_{{l \in C}} \sum_{{k \in M}}   \beta^s_{tlk} \leq  \sigma^\beta , \quad  \forall t \in \{\frac{|T|}{2}+1, \frac{|T|}{2}+2, ..., (T^{d^2}_{l}-T'_{lk})\}, \forall s \in N
\end{equation}

\end{spacing}

\indent

\paragraph{eVTOL Parts Procurement Constraints} \label{critical}

The number of eVTOL parts $i$ ordered from supplier $j$ on a particular day $t_1$ is expected to be received after a certain duration on day $t$. This duration is calculated by $t_1 = t - T''_{ij} - T'_{ijk}$. The EM would need to order the eVTOL parts $i$, considering both the lead time ($T''_{ij}$) required by supplier $j$ and the transportation time ($T'_{ijk}$) for delivery from supplier $j$ to the EM, in order to receive the eVTOL parts $i$ on time. Eq. \ref{pcons1} ensures first-stage procurement constraint, where $v_{tijk}$ represents the number of eVTOL parts $i$ to be ordered to supplier $j$ at the beginning of a given day $t$ in $F^P$ when transportation mode $k$ is chosen in the back-end of the supply chain. Similarly, Eq. \ref{pcons2} ensures second-stage procurement constraint, where the $v^s_{tijk}$ represents the number of eVTOL parts $i$ to be ordered to supplier $j$ at the beginning of a given day $t$ in $S^P$ in scenario $s$ when transportation mode $k$ is chosen in the back-end of the supply chain. Eq. \ref{pcons3} restricts the total number of eVTOL parts $i$ to be ordered on day $t$ in $F^P$ from supplier $j$ for part $i$ so that it does not exceed the capacity of that supplier to accept the EM's order. This ensures that the supplier can meet the demand of the EM for the specified part without exceeding its first-stage capacity $y_{ij}$. Eq. \ref{v^s} ensures that supplier $j$ can meet the demand of the EM for eVTOL part $i$ without exceeding its second-stage capacity $y^s_{ij}$ in $S^P$ in scenario $s$.

\begin{spacing}{0.5}

\begin{equation} \label{pcons1}
v_{t_1ijk} = z_{tijk}, \quad  t_1 = t - T''_{ij} - T'_{ijk}, t_1 \geq 1, \forall t \in \{1, 2, ..., \frac{|T|}{2}\}, \forall i \in P, \forall j \in S_i, \forall k \in M
\end{equation}

\begin{equation} \label{pcons2}
v^s_{t_1ijk} = z^s_{tijk}, \quad  t_1 = t - T''_{ij} - T'_{ijk}, t_1 \geq 1, \forall i \in P, \forall j \in S_i, \forall k \in M, \forall t \in \{\frac{|T|}{2}+1, \frac{|T|}{2}+2, ..., |T|\}, \forall s \in N
\end{equation}

\begin{equation} \label{pcons3}
\sum_{{t =1}}^{\frac{|T|}{2}} \sum_{{k \in M}} v_{tijk} \leq y_{ij}, \quad  \forall i \in P, \forall j \in S_i
\end{equation}

\begin{equation} \label{v^s} 
\sum_{t={\frac{|T|}{2}}+1}^{|T|} \sum_{{k \in M}} v^s_{tijkh} \leq  y^s_{ij}, \quad  \forall i \in P, \forall j \in S_i,  \forall s \in N
\end{equation}

\end{spacing}
\indent


\subsubsection{Planning Horizon Extension Due to Disruptions}

In response to disruptions in the supply chain, when the lead times of available suppliers in the AAM market increase during any iteration of RHA, the EM should extend the length of PH, renegotiate with customers, and revise deadlines accordingly over the extended PH, as discussed in Section \ref{rollingHorizon_approach}. The EM would incur penalty cost for extending the deadlines. To incorporate this penalty cost into $\Omega$, Eq. \ref{omega} is updated as shown in Eq. \ref{penalty_cons1}. Here $\Omega^\Upsilon$ represents the updated first-stage total cost for the iteration in which the length of PH is extended. $\Upsilon_l$ denotes the penalty cost incurred by the EM for customer $l$. Prior to initiating renegotiation with a customer to set a new deadline, the EM should calculate both the minimum and maximum possible delays for customer $l$. With this range of delays, the customer can then select a new deadline based on their own scheduling constraints and availability to receive the manufactured eVTOLs. There are four factors considered to determine the minimum delay for the customer. The first factor involves identifying the lead time of supplier $j$, which is the minimum among the lead times of all available suppliers in $S_i$ for eVTOL part $i$, while also satisfying the quality constraints. Consequently, a minimum lead time is selected for each eVTOL part $i$, and among these lead times, the maximum one is represented by $\Upsilon^2$. The second factor, as denoted by $\Upsilon^3$ in Eq. \ref{penalty_cons2}, determines the minimum transportation time in the back-end of the supply chain across all transportation modes in $M$. The third factor, as represented by $\Upsilon^4_l$ in Eq. \ref{penalty_cons3}, identifies the minimum transportation time for customer $l$ in the front-end of the supply chain across all transportation modes in $M$. The fourth factor, as denoted by $\Upsilon^5_l$ in Eq. \ref{penalty_cons4}, refers to the time required to manufacture $D^1_{l}$ eVTOLs for customer $l$. The minimum delay $\Upsilon^6_l$ for customer $l$ can be calculated using Eq. \ref{penalty_cons5}. Similarly, the maximum delay for customer $l$, as denoted by $\Upsilon^{7}_l$, can be computed following Eq. \ref{penalty_cons6}. In Eq. \ref{penalty_cons7}, $\Upsilon^8$ and $\Upsilon^9$ refer to the maximum lead time of supplier across all eVTOL parts and the maximum transportation time in the back-end of the supply chain across all transportation modes, respectively. $\Upsilon^{10}_l$ in Eq. \ref{penalty_cons7} represents the maximum transportation time for customer $l$ in the front-end of the supply chain across all transportation modes. Based on the possible range of delays ($\Upsilon^6_l$ - $\Upsilon^{7}_l$), the EM would renegotiate with customer $l$ and set the new deadline at the beginning of the iteration in which the length of PH is extended. The new deadline for customer $l$ is denoted by $T^{d^\Upsilon}_{l}$. If both the minimum and maximum delays are 0 for customer $l$, it implies that the EM would receive the eVTOL parts, the deadline of that customer is not affected in this situation, and the EM would not need to contact the customer for renegotiation. For customer $l$, $\Upsilon_l$ in Eq. \ref{eq149} refers to the penalty cost incurred by the EM for delaying the deadline for customer $l$ by one day from the original $T^{d^1}_{l}$. In iterations where disruptions occur requiring an extension of PH, the 3SCOPE model extends the length of the PH, adopts $\Omega^\Upsilon$ instead of $\Omega$, and updates the parameter $T^{d^1}_{l}$ in Eqs. \ref{obj} to \ref{v^s} with $T^{d^\Upsilon}_{l}$. If the situation reverts to typical market conditions, the 3SCOPE model utilizes the typical length of PH, $\Omega$, and $T^{d^1}_{l}$ again.

\begin{spacing}{0.4}

\begin{equation} \label{penalty_cons1}
\Omega^\Upsilon = \Omega +  \sum_{{l \in C}} \Upsilon_l
\end{equation}

\begin{equation} \label{penalty_cons2}
\Upsilon^3 = \min_{k \in M} \left\{ T'_{ijk} \; \middle| \; \Upsilon^2 = \max_{i \in P} \left( \min_{j \in S_i} T''_{ij} \right) \right\}
\end{equation}

\begin{equation} \label{penalty_cons3}
\Upsilon^4_l = \min_{k \in M} T'_{lk}, \quad \forall l \in C
\end{equation}

\begin{equation} \label{penalty_cons4}
\Upsilon^5_l = \frac{D^1_{l}}{R}, \quad \forall l \in C
\end{equation}

\begin{equation} \label{penalty_cons5}
\Upsilon^6_l = \begin{cases} \Upsilon^2 + \Upsilon^3 + \Upsilon^4_l + \Upsilon^5_l - T^{d^1}_{l}, & \text{if } T^{d^1}_{l} < \Upsilon^2 + \Upsilon^3 + \Upsilon^4_l + \Upsilon^5_l \\ 0, & \text{otherwise} \end{cases}, \quad \forall l \in C
\end{equation}

\begin{equation} \label{penalty_cons6}
\Upsilon^{7}_l = \begin{cases} \Upsilon^8 + \Upsilon^9 + \Upsilon^{10}_l + \Upsilon^{5}_l - T^{d^1}_{l}, & \text{if } T^{d^1}_{l} < \Upsilon^8 + \Upsilon^9 + \Upsilon^{10}_l + \Upsilon^{5}_l \\ 0, & \text{otherwise} \end{cases}, \quad \forall l \in C 
\end{equation}

\begin{equation} \label{penalty_cons7}
\Upsilon^9 = \max_{k \in M} \left\{ T'_{ijk} \; \middle| \; \Upsilon^8 = \max_{i \in P} \left( \max_{j \in S_i} T''_{ij} \right) \right\}
\end{equation}

\begin{equation} \label{penalty_cons8}
\Upsilon^{10}_l = \max_{k \in M} T'_{lk}, \quad \forall l \in C
\end{equation}

\begin{equation} \label{eq149}
\Upsilon_l =  \Upsilon^1_l |T^{d^\Upsilon}_{l} - T^{d^1}_{l}|, \quad \forall l \in C
\end{equation}

\end{spacing}





\subsection{Multi-cut Benders Decomposition}

Stochastic programming problems typically require significant computational resources, particularly as they transition from small to large-scale problems, resulting in increased complexity and running time \cite{almeida2021decomposition}. To enhance the scalability of our 3SCOPE model, enabling it to handle large-scale problems effectively, we implement a relaxation technique followed by a multi-cut Benders decomposition method, as described in \cite{you2013multicut}. The method involves breaking down the 3SCOPE model into a master problem and relaxed slave problems, where recourse decisions are made. With first-stage integer variables in $F^P$ and second-stage continuous variables in $S^P$, the 3SCOPE model can be effectively approached using the L-Shaped Method \cite{laporte1993integer}, which is a scenario-based decomposition technique based on Benders decomposition \cite{benders2005partitioning}.

\subsection{Benchmark Models} \label{bench}

In addition to the 3SCOPE model, which is built on two-stage stochastic programming, we also consider three benchmark models. We compare the profits generated from the 3SCOPE model with the operating profits generated from these three benchmark models. Following \cite{inproceedings}, where the deterministic version of their stochastic model was compared, we also conduct this comparison to reflect the impact of uncertainty on supply chain performance. Therefore, our first benchmark model is a deterministic version of the 3SCOPE model. This alternative option of the 3SCOPE model is implemented directly for the periods in AH, using the same input values considered for $F^P$ in RHA to run the 3SCOPE model. Next, we develop a stochastic heuristic model, also developed in \cite{gruler2018combining}, which serves as a greedy heuristic version of the 3SCOPE and is the second benchmark model in our study. To develop the stochastic heuristic model, we utilize a heuristic for selecting suppliers with the lowest prices, as presented in Algorithm \ref{algo2}. Lastly, we consider a stochastic sequential approach where the stochastic model is run first to address the challenges as considered in \cite{inproceedings}, followed by the second model for manufacturing scheduling, as discussed in Section \ref{gaps}. Then, we compare the profits generated from these sequential models with the profits generated from the integrated 3SCOPE model.

\begin{spacing}{0.8}
\begin{algorithm}[tb]
\caption{Stochastic Heuristic Model}
\label{algo2}

\begin{algorithmic}[1]
    \For{each iteration $h \in$ RHA}
        \For{each scenario $s$ in $S$}
            \For{each eVTOL part $i$}
                \State Select supplier $j \in S_i$ with the lowest procurement price for eVTOL part $i$ for manufacturing $D^1_{l}$ and $D^2_{l}$.
                \indent \indent \indent eVTOLs for $F^P$ and $S^P$, respectively.
                \State Determine the number of eVTOL parts $i$ to be procured from supplier $j \in S_i$.
                \State Select transportation mode for each delivery which incurs the lowest transportation cost and emission 
                    \indent \indent \indent  cost in the back-end of the supply chain.
                \State Determine inventory cost for holding eVTOL parts.
            \EndFor
            \State Determine the cost for manufacturing $D^1_{l}$ and $D^2_{l}$ eVTOLs for $F^P$ and $S^P$, respectively.
            \State Determine inventory cost for holding eVTOLs manufactured.
            \State Select transportation modes for transporting eVTOLs manufactured to customer $l$, ensuring they arrive by the deadline.
            \State Consider the transportation mode that incurs the lowest transportation cost and emission cost in the back-end of the 
              \indent \indent supply chain.
        \EndFor
        \State Save the profit generated for $F^P$.
    \EndFor
    \State \textbf{return} Optimized decisions and profit.
\end{algorithmic}

\end{algorithm}

\end{spacing}


\subsection{Computational Environment}

The 3SCOPE model described in Section \ref{model_approach} and the benchmark models described in Section \ref{bench} are implemented in Python 3.10.0 using the Gurobi-Python API. Using Gurobi 10.0.1x64 sovler, we run the models on a computer equipped with an Intel® Core™ i9-10980XE processor (3.00 GHz, 36 cores) and 128 GiB of memory.

\section{Results}\label{sec4}


This section includes the experimental setup for running the 3SCOPE model and the analysis of various numerical cases for the EM. It also presents the results generated from the 3SCOPE model and benchmark models across different numerical cases.



\subsection{Experimental Setup} \label{data}

This study employs a rolling horizon approach, as discussed in Section \ref{rollingHorizon_approach}, where $T$ = \{day 1, day 2, ..., day 360\}, and AH ranges from day1 to day900. Each PH is divided into two periods, each consisting of 180 days. The RHA progresses through the periods in $T^P$ = \{Period 1, Period 2, Period 3, Period 4, Period 5\}, with Period 1 starting from the year 2026. The start and end days for each period over AH are listed in Table \ref{period_years}. In the experimental setup for running the 3SCOPE model, we consider a conceptual design for intercity travel, which is capable of carrying four passengers and one pilot, as presented in \cite{akash2021design}. The authors identified and provided weight estimates for various components of the eVTOL. From this study, we select a subset of eVTOL parts and also include an additional eVTOL part. The set of eVTOL parts are denoted by $P = [a, b, c, d, e]$, where $a$ stands for fuselage, wing and v-tail, $b$ for propeller, $c$ for battery, $d$ for motor, and $e$ for seat. The weight of eVTOL part $i$, $W_i$, is as follows: $W_a$ = 0.426 tons, $W_b$ = 0.0033 tons, $W_c$ = 0.128 tons, and $W_d$ = 0.02 tons \cite{akash2021design}. Additionally, for the eVTOL part $e$, which was not covered in \cite{akash2021design}, we consider an estimated weight of $W_e$ = 0.012 tons based on sources \cite{seat, seat1}. To produce a single eVTOL, the values of $\theta_i$ for eVTOL parts $i$ are $\theta_a$ = 1, $\theta_b$ = 6, $\theta_c$ = 4, $\theta_d$ = 6, and $\theta_e$ = 5. This results in a baseline eVTOL weight of $W^e$ = 1.1128 tons, incorporating a miscellaneous weight of 0.01 tons. 


\begin{table}[ht]
\begin{spacing}{1}
\centering
\caption{Period start and end days over AH.}
\label{period_years}
\begin{tabular}{lccccc}
\toprule
\textbf{} & \textbf{Period 1} & \textbf{Period 2} & \textbf{Period 3} & \textbf{Period 4} & \textbf{Period 5} \\ \midrule
\textbf{Start Day} & 1 & 181 & 361 & 541 & 721 \\
\textbf{End Day} & 180 & 360 & 540 & 720 & 900 \\
\textbf{Year} & 2026 ($1^{st}$ Half) & 2026 ($2^{nd}$ Half) & 2027 ($1^{nd}$ Half) & 2027 ($2^{nd}$ Half)  & 2028 ($1^{st}$ Half)  \\ 
\bottomrule
\end{tabular}
\end{spacing}
\end{table}




Since the market for eVTOL is not yet established, we conduct an analysis of the market for existing suppliers who provide parts for automotives and traditional aircraft and select five different suppliers from both within and outside the US for each eVTOL part. The set of suppliers selected for each eVTOL part is denoted by $[s1, s2, s3, s4, s5]$. Given the scarcity of available data regarding the lead time of suppliers, we randomly sample values for $T''_{ij}$ from a uniform distribution within the 1 to 5-day range. A similar range of lead time previously employed in \cite{cachon2000supply}. In this study, we have sourced the values of $x_{ij}$, $y_{ij}$, and $L$ from various sources, including \cite{airframe, weigon, p_e_1}. As discussed in Sections \ref{sec2} and \ref{quality_cons}, aviation products demand the strictest adherence to quality requirements. According to the priorities of the major parts of the eVTOLs manufactured, we assume $Q^r_i$ for eVTOL parts $i$ to be as follows: $Q^r_a$ = 7.5, $Q^r_b$ = 8, $Q^r_c$ = 8.5, $Q^r_d$ = 8, and $Q^r_e$ = 7.5. We utilize a rating scale ranging from 0 to 10 to quantify the quality of the eVTOL parts provided by the suppliers. Here, 0 signifies the lowest quality, while a rating of 10 represents the highest quality among suppliers involved in the manufacturing of eVTOL parts. We assume $Q^b$ to be 8 out of 10, which indicates that any $q$ value below 8 would not meet the minimum quality standard required by the FAA, resulting in a failure to obtain certification from the FAA. We then randomly sample values for $Q_{ij}$ from a uniform distribution ranging between 6 and 10. We set $\kappa$ to 0.2, indicating a 20\% increase in eVTOL pre-orders from a customer if the quality of eVTOLs in a given period changes by 1 point on the scale of 10 from the quality of eVTOL-batch manufactured in the previous period. The first-stage $Q_{ij}$, $x_{ij}$, and $y_{ij}$ values in the first iteration of RHA are listed in Table \ref{suppliers}. Here, supplier capacity is measured by the number of eVTOL parts available to the supplier for providing to the EM. In the same iteration, $L$ is considered to be \$410k. At the beginning of each iteration, the EM would update these values to reflect the current information available in the market.

\begin{table}[b!]
\begin{spacing}{1}
\centering
\caption{First-Stage supplier data for eVTOL parts considered in the first iteration of the rolling horizon approach.} \label{suppliers}
\begin{tabular}{cccccc|ccccc|ccccccc}
\toprule
\textbf{eVTOL Parts} & \multicolumn{5}{c|}{$a$} & \multicolumn{5}{c|}{$b$} & \multicolumn{5}{c}{$c$} \\
\midrule
\textbf{Supplier} & $s1$ & $s2$ & $s3$ & $s4$ & $s5$ & $s1$ & $s2$ & $s3$ & $s4$ & $s5$ & $s1$ & $s2$ & $s3$ & $s4$ & $s5$ \\
\midrule
\textbf{Quality} & 8.7 & 9.0 & 7.2 & 8.6 & 8.7 & 8.7 & 8.8 & 8.6 & 8.8 & 9.1 & 8.9 & 8.5 & 9.2 & 8.5 & 9.2 \\
\textbf{Price (\$1000)} & 30 & 29.5 & 28.5 & 30.3 & 30.7 & 10 & 9.8 & 11 & 9.5 & 10.5 & 22.2 & 22 & 22.5 & 22.9 & 21.7 \\
\textbf{Capacity} & 10 & 15 & 12 & 15 & 8 & 115 & 100 & 120 & 110 & 130 & 90 & 95 & 100 & 85 & 80 \\
\midrule
\textbf{eVTOL Parts} & \multicolumn{5}{c|}{$d$} & \multicolumn{5}{c|}{$e$} \\
\midrule
\textbf{Supplier} & $s1$ & $s2$ & $s3$ & $s4$ & $s5$ & $s1$ & $s2$ & $s3$ & $s4$ & $s5$ \\
\midrule
\textbf{Quality} & 8.2 & 8.4 & 8.1 & 9.3 & 8.8 & 8.1 & 9.4 & 8.4 & 7.5 & 7.0 \\
\textbf{Price (\$1000)} & 11.2 & 11 & 10.7 & 12 & 11.8 & 9.5 & 10.5 & 9.8 & 9.5 & 8.0 \\
\textbf{Capacity} & 135 & 70 & 80 & 95 & 75 & 80 & 85 & 90 & 95 & 100 \\
\bottomrule
\end{tabular}
\end{spacing}
\end{table}


A demand comparison focused on urban areas within the US was done by the NASA in \cite{hasan2019urban}. Within this analysis, air taxi was found to be viable for 5-seater eVTOLs across a sample of ten urban areas. In our study, we focus on the four cities from different states identified in the NASA study: Dallas ($c1$), New York ($c2$), Chicago ($c3$), and San Francisco ($c4$). These cities constitute the customer set $C = [c1, c2, c3, c4]$ in our study, selected for their significant potential and projected demand for eVTOLs. The EM in our study is assumed to capture 40\% of the demand for eVTOLs in each of these cities. We set the values of $D^2_{l}$ for each customer $l$ by following the demand projection of eVTOLs, as given in market reports \cite{market1, market3}. We explore a range of deadline combinations for customers to demonstrate the flexibility of our model in accommodating diverse sets of deadlines. The values of $D^1_{l}$, $D^2_{l}$, $T^{d^1}_{l}$, and $T^{d^2}_{l}$ considered in this study for customer $l$ are listed in Table \ref{demand_deadline} for all iterations of the RHA. We have addressed uncertainties in our model parameters, including $x^s_{ij}$, $y^s_{ij}$, $L^s$, and $D^s_{l}$, and generated 30 scenarios for running the 3SCOPE model according to the process outlined in Section \ref{sc_gen}. At the beginning of each iteration, the values of these parameters across the scenarios are updated. For $x^s_{ij}$, $y^s_{ij}$, and $L^s$, we consider their EVs to be the values of $x_{ij}$, $y_{ij}$, and $L$ in the same iteration, respectively. The SDs are assumed to be 10\% of their respective EVs. As for $D^s_{l}$, a left-skewed normal distribution is assumed to be followed, with an EV and SD of 4 and 2 in the initial iteration of RHA, respectively. In subsequent iterations of RHA, the EV for $D^s_{l}$ follows the demand projection of eVTOL, as outlined in market reports \cite{market1, market3}.

\begin{table}[t]
\begin{spacing}{1}
\centering
\caption{eVTOL pre-orders and deadline data for customers for all iterations of the rolling horizon approach.} \label{demand_deadline}
\begin{tabular}{c|c|cccc|cccc|cccc|cccc}
\toprule
& & \multicolumn{4}{c|}{$\mathbf{D^1_{l}}$} & \multicolumn{4}{c|}{$\mathbf{D^2_{l}}$} & \multicolumn{4}{c|}{$\mathbf{T^{d^1}_{l}}$ } & \multicolumn{4}{c}{$\mathbf{T^{d^2}_{l}}$} \\ 
\textbf{Iteration} & \textbf{Periods}  & \textbf{c1} & \textbf{c2} & \textbf{c3} & \textbf{c4} & \textbf{c1} & \textbf{c2} & \textbf{c3} & \textbf{c4} & \textbf{c1} & \textbf{c2} & \textbf{c3} & \textbf{c4} & \textbf{c1} & \textbf{c2} & \textbf{c3} & \textbf{c4} \\
\midrule
1  & 1, 2 & 6 & 5 & 3 & 3 & 3 & 2 & 1 & 0 & 70 & 90 & 100 & 150 & 300 & 280 & 260 & 320 \\
2  & 2, 3 & 6 & 5 & 3 & 3 & 3 & 2 & 1 & 1& 300 & 280 & 260 & 320 & 490 & 450 & 490 & 510 \\
3 &  3, 4 & 8 & 7 & 4 & 4 & 3 & 2 & 2 & 1 & 490 & 450 & 490 & 510 & 700 & 600 & 655 & 635 \\
4 &  4, 5 & 8 & 7 & 4 & 4 & 6 & 4 & 3 & 2 & 700 & 600 & 655 & 635 & 815 & 850 & 888 & 830 \\
5 &  5 & 11 & 10 & 5 & 5 & - & - & - & - & 815 & 850 & 888 & 830 & - & - & - & - \\
\bottomrule
\end{tabular}
\end{spacing}
\end{table}

Sources such as \cite{lilium, joby, weigon} offer insights into the EM details, including its eVTOL manufacturing targets, profit margins, and capacities. However, due to the emerging nature of the eVTOL industry, actual data of the EM is currently limited. Hence, concerning the maximum inventory capacity of the EM, we assume the capacities of the EM for the eVTOL parts in $P$ to be 15, 60, 50, 60, and 50 units per day, respectively. Additionally, we assumed $R$ to be 5 units per day in our study. The average price of an eVTOL ranges between \$1 million and \$4 million, as indicated by \cite{sell1, sell2, sell3}. In our study, we assume that the EM would sell its five-seater eVTOLs at a price of \$1.5 million. Another EM specific data is holding cost for inventory which encompasses expenses such as damaged products, storage space, depreciation, and insurance \cite{odedairo2020system}. For a company specializing in the design, manufacturing, and supply of electronic products, the holding cost is capped at a maximum of 32.9\% of the current inventory value, as indicated in \cite{azzi2014inventory}. In accordance with this, we consider 32.9\% of the average value of an eVTOL part $i$ procured from supplier $j$, with the corresponding prices provided in Table \ref{suppliers}, to calculate the holding cost $H^a_i$ for each eVTOL part $i$ in the inventory. Subsequently, the holding cost $H^b$ for holding an eVTOL manufactured is estimated as 32.9\% of cumulative value of the eVTOL parts required for its manufacturing. In this study, we assume a delay-penalty policy applied to customer orders, as presented in Table \ref{tab:delay_penalty}. It depicts a progressive penalty policy based on the number of days a delivery is delayed by the EM for a customer pre-order. Initially, a delay of 1-15 days incurs a penalty of 0.05\% of the contract value of eVTOLs per customer per day, with incremental increases of penalty for each subsequent 15-day period.

\begin{table}[t]
\begin{spacing}{1}
    \centering
        \caption{Delay-penalty policy.}
    \begin{tabular}{ccccccc}
        \hline
        \textbf{Delay (days)} & 1-15 & 16-30 & 31-45 & 46-60 & 61-75 & 76-90\\
        \hline
        \textbf{Penalty (\% of contract value)} & 0.05 & 0.1 & 0.3 & 0.5 & 0.7 & 1 \\
        \hline
    \end{tabular}
    \label{tab:delay_penalty}
\end{spacing}
\end{table}


In this study, we assume the location of the EM to be in California, where Joby Aviation's headquarter is located, and we list suppliers from various locations within and outside the US. We consider four transportation modes: truck, rail, ship, and aircraft. To determine the optimal transportation mode in both of the ends of the supply chain among these four options, we create origin-destination pairs by pairing suppliers with EMs and EMs with customers. We then gather data on available transportation modes for each origin-destination pair for the four modes, along with their corresponding shipping distances and emissions. We use a calculator provided in \cite{eco-calculator} along with the Google Maps API for these calculations. An example of a pair between supplier $s3$ of part $b$ and the EM is presented in Table \ref{mode}. The supplier $s3$ of part $b$ is an international supplier located outside the USA with no ground connection. Therefore, truck and rail modes are unavailable due to the absence of a feasible route between the supplier and the EM for these modes of transportation. For the available modes—--ship and aircraft—--each mode is further analyzed through three sub-ends of the supply chain. When the transportation mode is ship, the total shipping distance between the supplier and the EM is approximately 10,736 km, resulting in a total emission of 0.13 tons of GHG per ton of eVTOL parts transported. We categorize the total shipping distance and emissions into three distinct segments as follows: 1) the first segment involves transporting eVTOL parts via truck from the supplier's location (referred to as origin) to seaport 1, which is near the origin; 2) the second segment involves selecting the ship mode to transport eVTOL parts from seaport 1 to seaport 2, located near the EM's location; and 3) the third segment involves selecting truck mode to transport the eVTOL parts from seaport 2 to the destination, which is the EM's location. This detailed breakdown underscores the multi-mode transportation, shedding light on the distinct contributions of each segment to the total shipping distance and total emission. The locations of the customers where the EM should deliver the eVTOLs manufactured are listed in Table \ref{LOcation_transportation}. In this table, the transportation times needed to deliver the eVTOLs via different modes from the EM's plant to the customers are also presented. Similar to calculating the transportation time for the supplier to the EM pairs, we utilize a calculator provided in \cite{eco-calculator} along with the Google Maps API to compute the transportation times for the EM and customer pairs. In transportation, it is essential to factor in both the speed and maximum daily driving hours to calculate the delivery time based on the data of shipping distance. To conduct this calculation, we utilize the speed and daily driving hours limits data presented in Table \ref{speed}, which is sourced from \cite{allowed_speed}. Using historical data on the average freight cost per ton-km from 1990 to 2021, sourced from \cite{unit_cost_mode}, we projected the transportation cost per ton-km for the years 2026 to 2028. Considering the AH spanning from 2026 to 2028, we provide the projected average freight costs for four transportation modes in Table \ref{tab_cost}. To achieve emission-reduction goals and mitigate GHG emissions, The Cap-and-Trade Program states that the cost per ton of GHG emissions is \$14-\$18 in year 2023 \cite{emission_cost}. Additionally, prices are projected to rise by 5\% plus inflation annually \cite{emission_cost1}. In our study, we start with a value of \$16.5 for the year 2026 as the initial reference point. This value is then increased by 5\% each subsequent year.

\begin{table}[t]
\begin{spacing}{1}
\centering
\caption{Transportation and emission data for supplier $s3$ of part $b$ and EM pair.} \label{mode}
\begin{tabular}{c|c|ccccc|ccc}
\toprule
 \textbf{Mode} & \textbf{Segment}&  \textbf{From}  & \textbf{To}   & \textbf{Multi} & \textbf{Distance } & \textbf{Emission}  & \textbf{Total } & \textbf{Total } &\\
 
 &  &  & & \textbf{Mode}  &  \textbf{ (km)}    & \textbf{(ton)} &\textbf{Distance } & \textbf{Emission}  \\
  &  &  &  &    &   & &\textbf{ (km)} & \textbf{(ton)}  \\
\midrule
Truck & - &- & - & - & - & - & - & - \\
\midrule
Rail & -&- & - & - & - & - & - & - \\
\midrule
  & 1  & {Origin} & {Seaport 1} & Truck & 31.68 & 0.01 &  \\
  Ship  &2&{Seaport 1} & {Seaport 2} & Ship & 10528.22 & 0.08 &  10,736.94 & 0.13\\
    &3&{Seaport 2} & {Destination} & Truck & 177.04 &  0.04& &  \\
\midrule
 &1&  {Origin} & {Airport 1} & Truck & 34.57 & 0.01 & \\
 Aircraft  &2& {Airport 1} & {Airport 2} & Aircraft & 10612.09 & 8.27 & 10,753.44 & 8.31 \\
    &3&{Airport 2} & {Destination} & Truck & 106.78 & 0.03 & & \\

\bottomrule
\end{tabular}
\end{spacing}
\end{table}

\begin{table}[t]
\begin{spacing}{1}
\centering
\caption{Transportation times for eVTOL delivery via different modes from the EM's plant to the customers.}
\label{LOcation_transportation}
\begin{tabular}{cccccc}
\toprule
\textbf{Customer} & \textbf{Destination} & \multicolumn{4}{c}{\textbf{Transportation Mode}} \\ \cmidrule(lr){3-6} 
         &  & Truck   & Rail   & Ship  & Air  \\ \midrule
c1       & Dallas Central Business District Vertiport & 4       & 2      & 14    & 1    \\
c2       & John F Kennedy International Airport & 6       & 2      & 14    & 1    \\
c3       & Vertiport Chicago & - & - & 20    & 1    \\
c4       & San Francisco International Airport & 1       & 1      & 1     & 1    \\ \bottomrule
\end{tabular}
\end{spacing}
\end{table}

\begin{table}[t]
\begin{spacing}{1}
\centering
\caption{Speed limit and maximum daily driving hours for different cargo transportation modes.}\label{speed}
\begin{tabular}{ccccc}
\toprule
\textbf{Mode} & \textbf{Truck} & \textbf{Rail} & \textbf{Ship} & \textbf{Air} \\
\midrule
\textbf{Speed limit (km/h)} & 80 & 120 & 30 & 900 \\
\textbf{Maximum daily driving hours} & 11 & 24 & 24 & 14 \\
\bottomrule
\end{tabular}
\end{spacing}
\end{table}

\begin{table}[t]
\begin{spacing}{1}
\centering
\caption{Projected average freight costs for four transportation modes in U.S. dollars per ton-km.} \label{tab_cost}
\begin{tabular}{l S[table-format=3.3] S[table-format=2.4] S[table-format=1.4] S[table-format=1.4]}
\toprule
\textbf{Year} & {\textbf{Aircraft}} & {\textbf{Truck}} & {\textbf{Rail}} & {\textbf{Ship}} \\
\midrule
2026 & 83.429 & 12.5521 & 2.9621 & 2.0041 \\
2027 & 85.211 & 12.8216 & 3.0180 & 2.0479 \\
2028 & 86.993 & 13.0911 & 3.0739 & 2.0917 \\
\bottomrule
\end{tabular}
\end{spacing}
\end{table}

\subsection{Numerical Analysis}

\subsubsection{Numerical Case 1: Base Case}
\begin{figure} 
\begin{subfigure}{1\textwidth}
  \centering
  \includegraphics[width=1\linewidth]{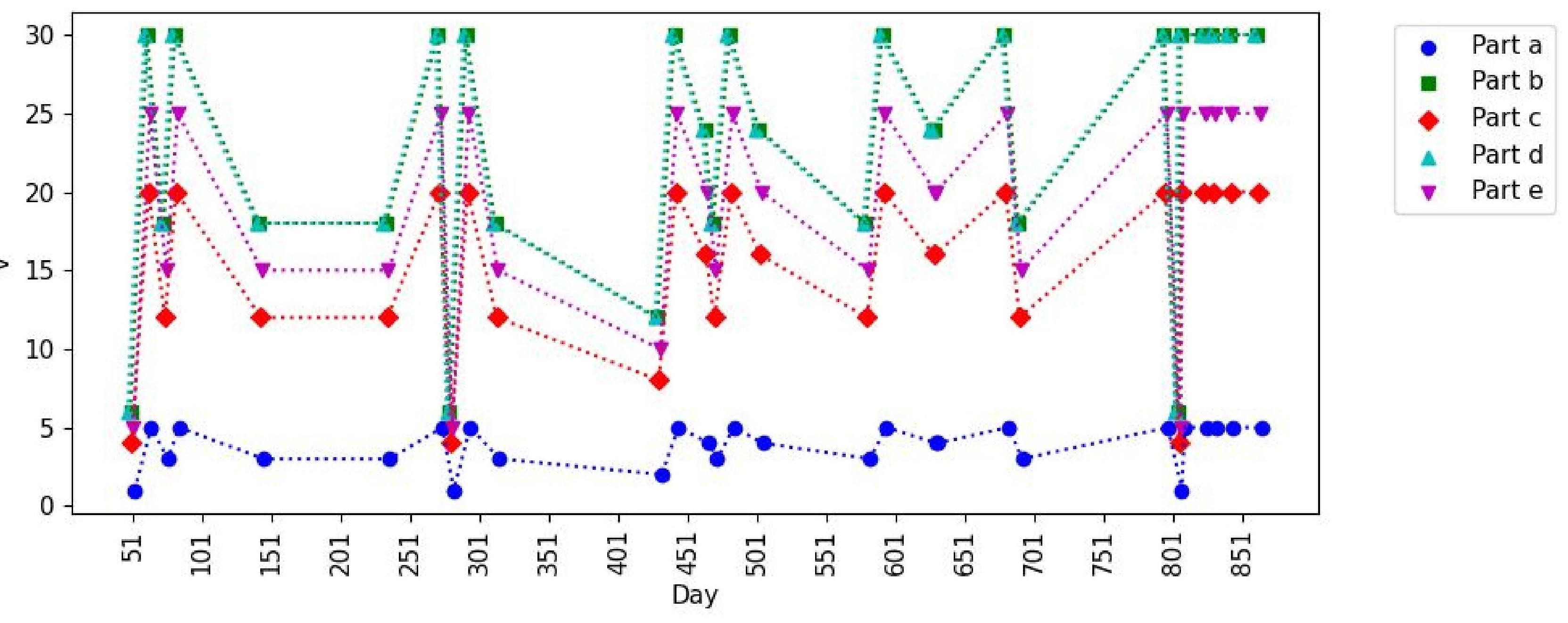}
  \caption{Optimal number of eVTOL parts ordered across AH}
  \label{fig1a}
\end{subfigure}
\begin{subfigure}{1\textwidth}
  \centering
  \includegraphics[width=1\linewidth]{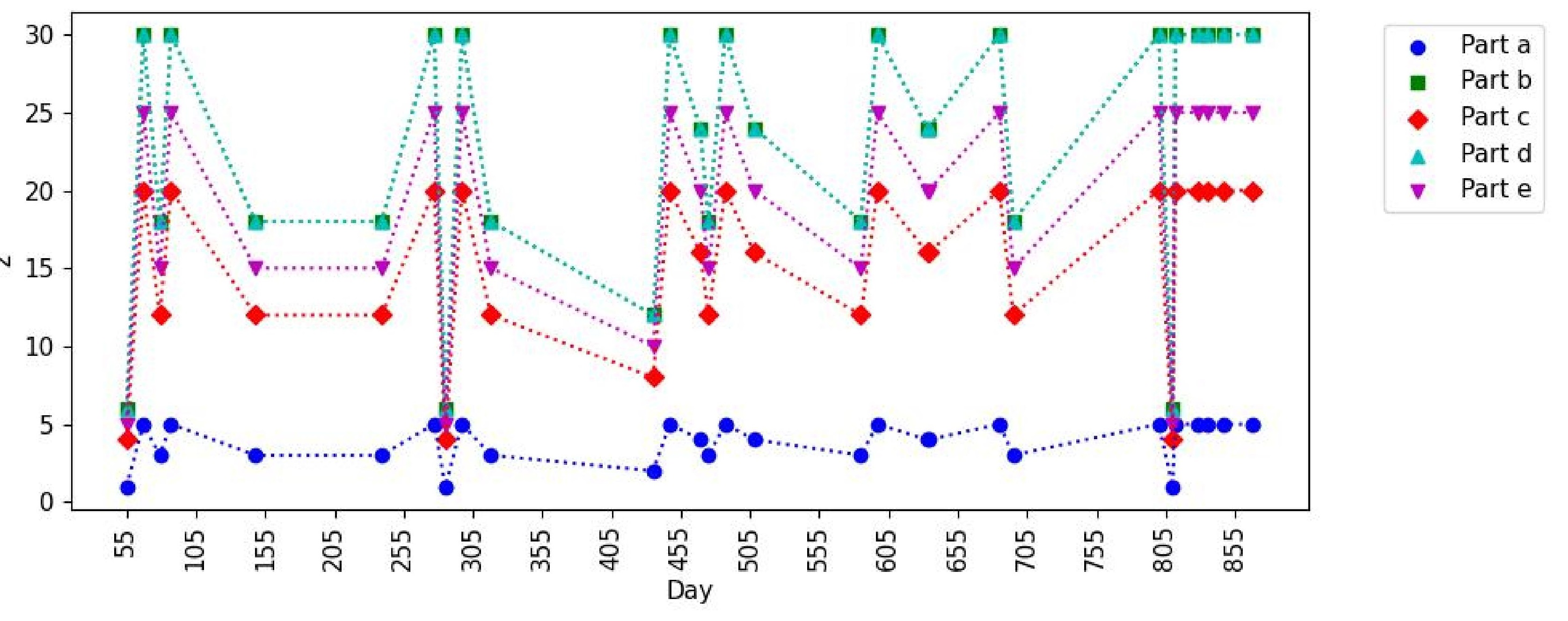}  
  \caption{Optimal number of eVTOL parts received across AH}
  \label{fig1b}
\end{subfigure}
\begin{subfigure}{1\textwidth}
  \centering
  \includegraphics[width=1\linewidth]{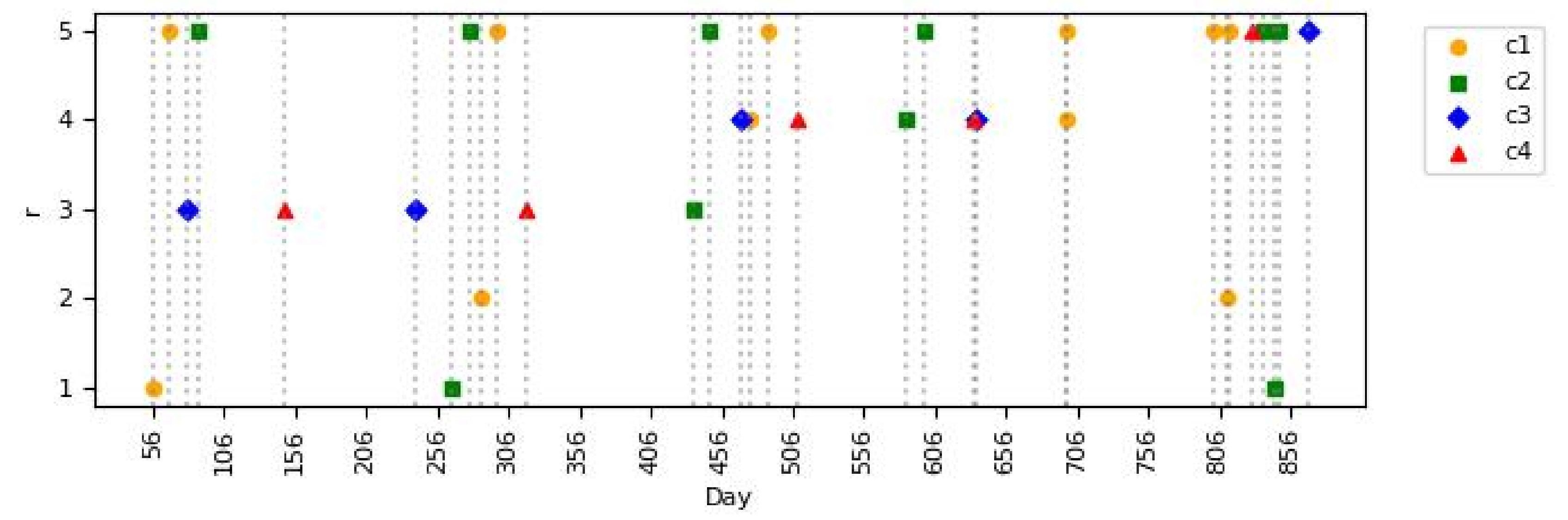}  
  \caption{Optimal number of eVTOLs manufactured across AH}
  \label{fig1c}
\end{subfigure}
\caption{Optimal number of eVTOL parts ordered, optimal number of eVTOL parts  received, and optimal number of eVTOLs manufactured across AH.}
\label{fig1}
\end{figure}

\begin{figure}[h!tb]
\begin{subfigure}{0.5\textwidth}
  \centering
  \includegraphics[width=1\linewidth]{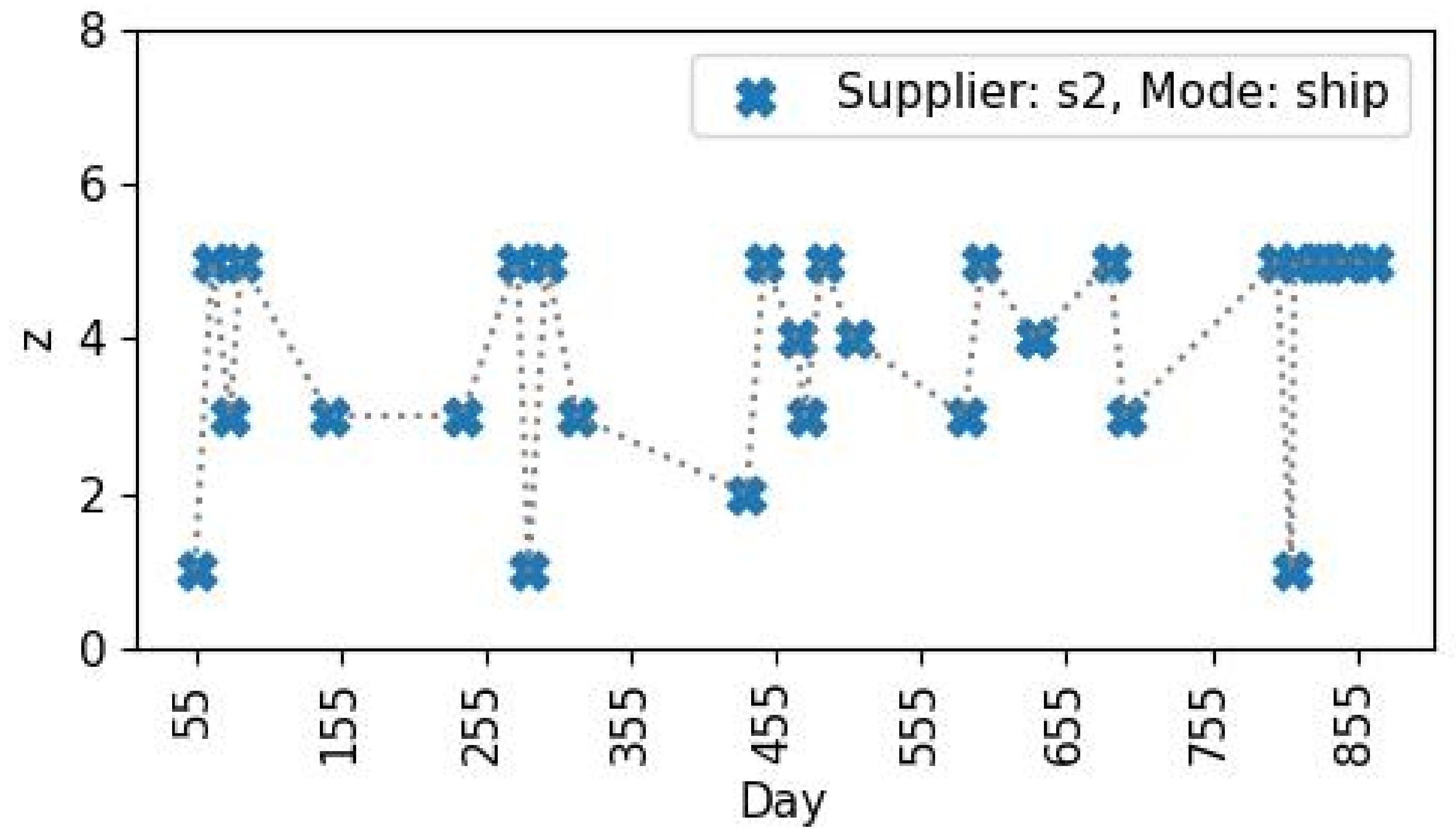}
  \caption{Optimal number of eVTOL parts $\boldsymbol{a}$ received}
  \label{figSM1a}
\end{subfigure}
\begin{subfigure}{0.5\textwidth}
  \centering
  \includegraphics[width=1\linewidth]{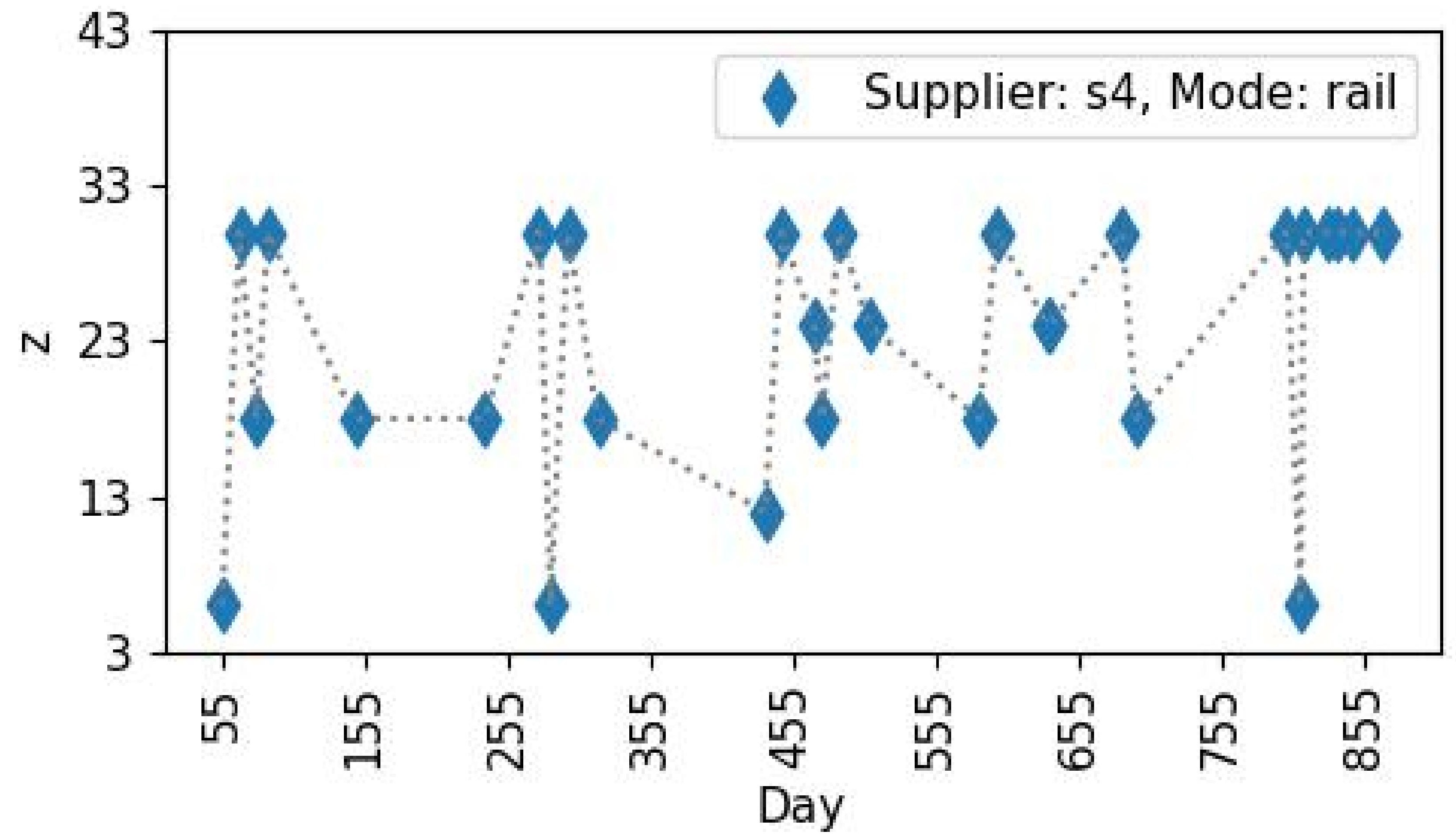}  
  \caption{Optimal number of eVTOL parts $\boldsymbol{b}$ received}
  \label{figSM1b}
\end{subfigure}
\begin{subfigure}{0.5\textwidth}
  \centering
  \includegraphics[width=1\linewidth]{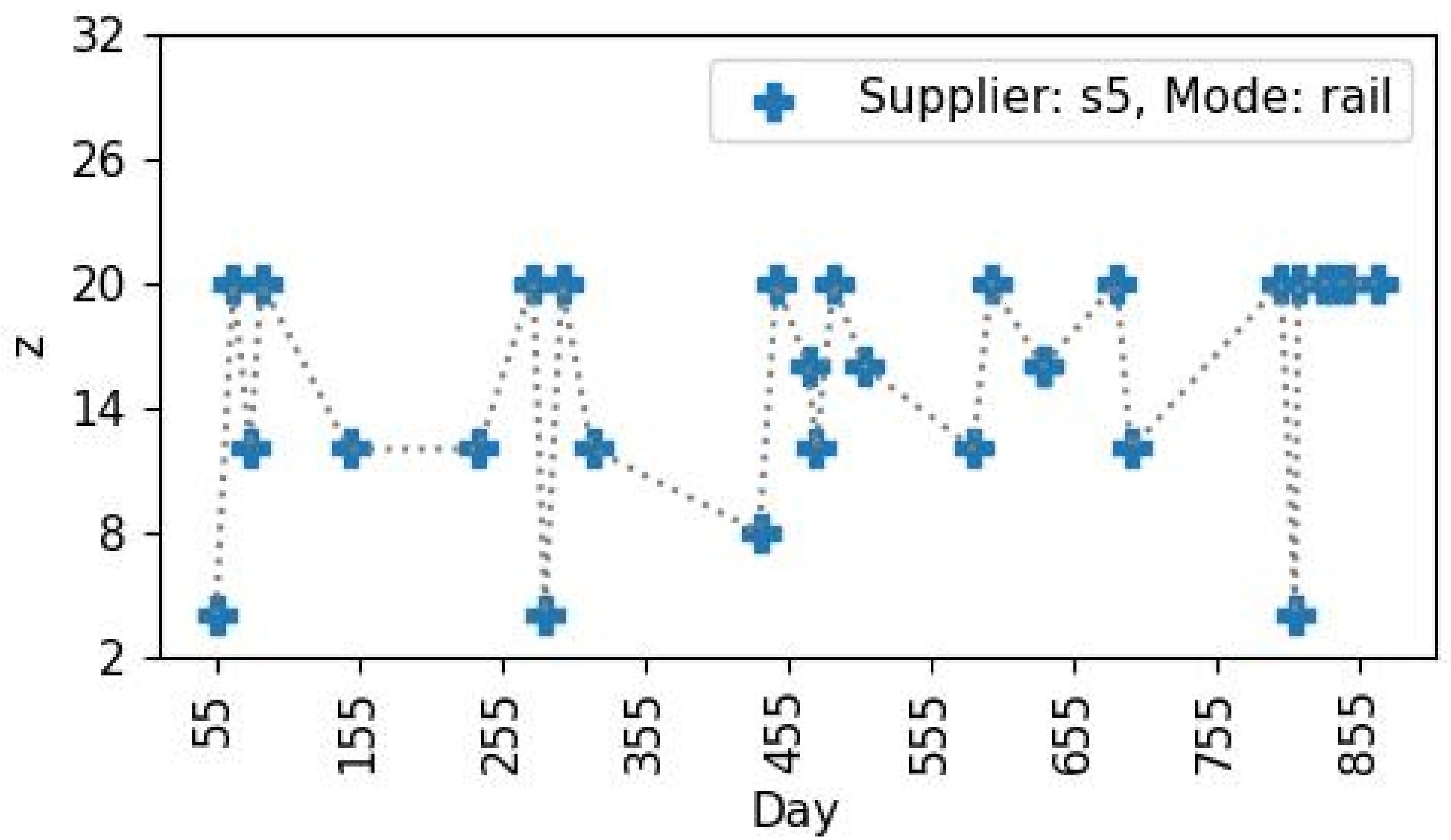}  
  \caption{Optimal number of eVTOL parts $\boldsymbol{c}$ received}
  \label{figSM1c}
\end{subfigure}
\begin{subfigure}{0.5\textwidth}
  \centering
  \includegraphics[width=1\linewidth]{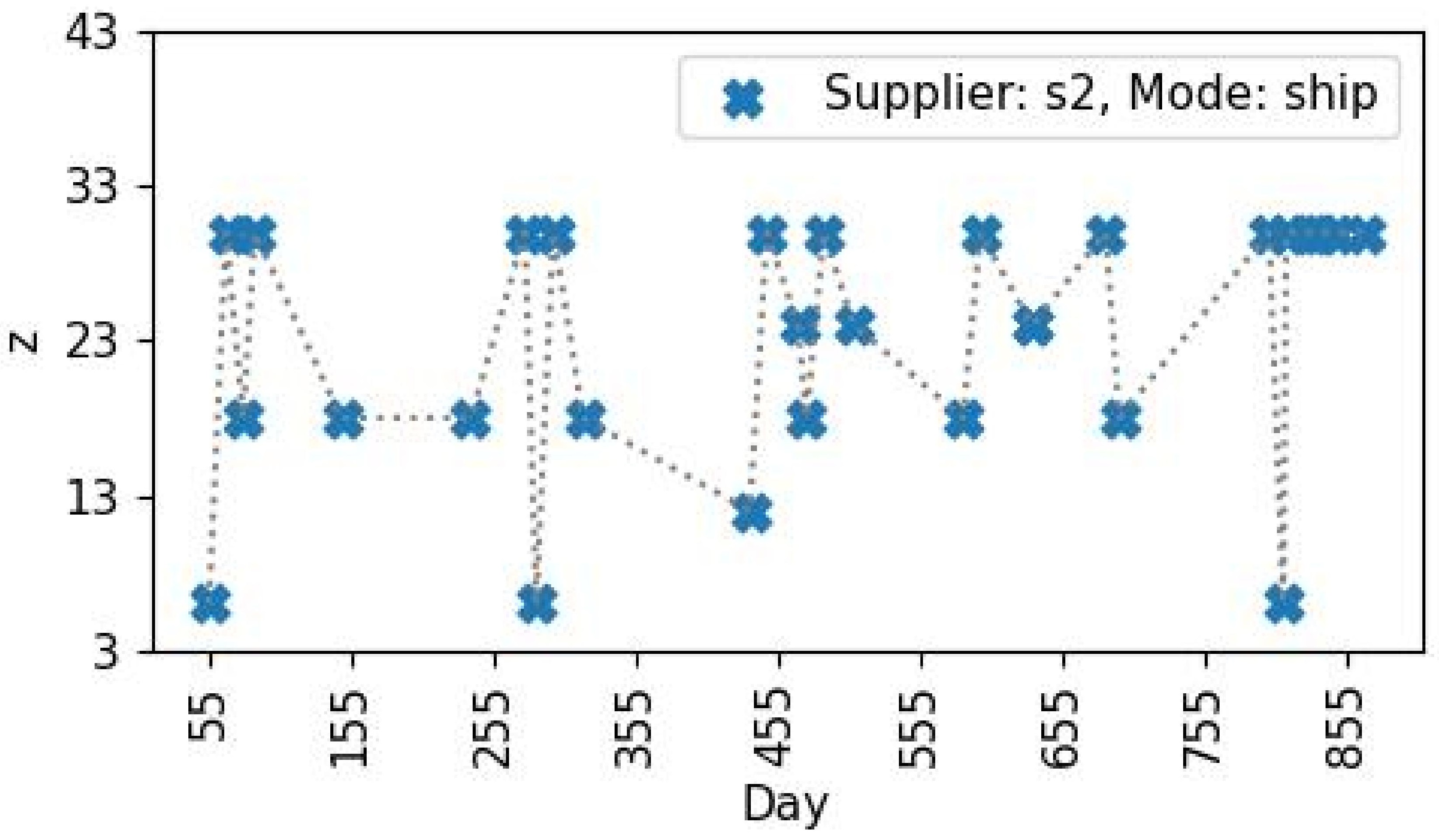}  
  \caption{Optimal number of eVTOL parts $\boldsymbol{d}$ received}
  \label{figSM1d}
\end{subfigure}

\centering
\begin{subfigure}{0.5\textwidth}
  \centering
  \includegraphics[width=1\linewidth]{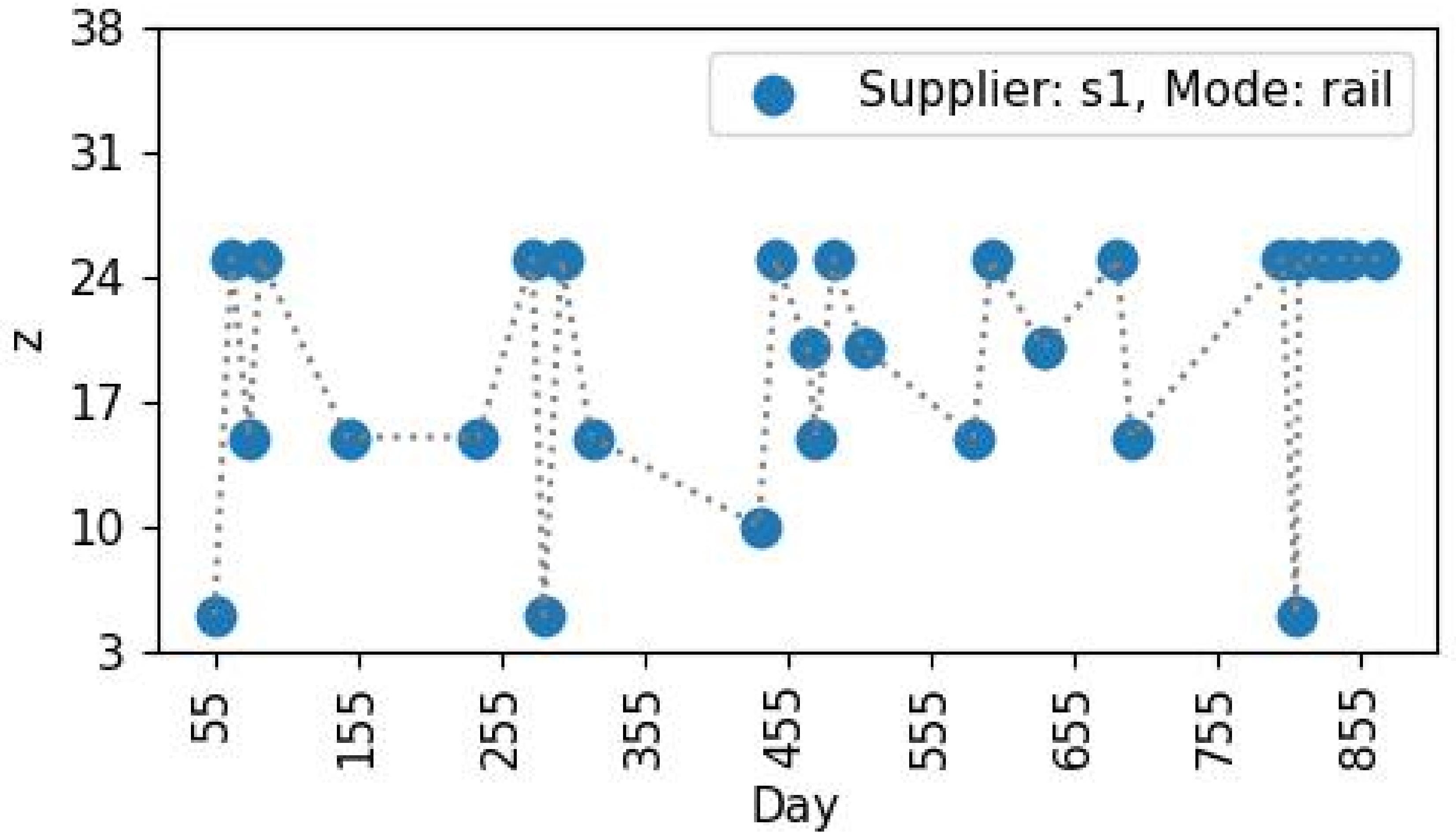}  
  \caption{Optimal number of eVTOL parts $\boldsymbol{e}$ received}
  \label{figSM1e}
\end{subfigure}
\caption{Optimal number of eVTOL parts received from the selected suppliers via optimal transportation mode across AH.} 
\label{figSM1}
\end{figure}

Following the Section \ref{data}, we study our first numerical case, which is considered to be the base case to compare with rest of the numerical cases. The results obtained from the 3SCOPE model for 900 days in AH show the optimal day and the optimal quantity of eVTOL parts to be ordered from suppliers in AH, as depicted in Figure \ref{fig1a}. Additionally, the results illustrate the optimal day and the optimal quantity of eVTOL parts to be received by the EM from the suppliers, as shown in Figure \ref{fig1b}. The optimal number of eVTOL parts is represented by five distinct colored markers, with each marker corresponding to a specific eVTOL part. To provide a visual representation of individual lines for each eVTOL part, all markers for each part are connected by dotted lines over AH. A detailed version of Figure \ref{fig1b} is provided in Figure \ref{figSM1}, illustrating the optimal supplier for each eVTOL part and the optimal transportation mode for delivering the optimal number of eVTOL parts by that supplier to the EM on the optimal day. Next, the optimal day and the optimal number of eVTOLs to be manufactured by the EM for different customers are displayed in Figure \ref{fig1c}. Here the optimal number of eVTOLs manufactured is presented with four different colored markers, where each colored marker is for each customer. Vertical dotted lines are used for visual representation, highlighting the optimal day for eVTOL manufacturing for each customer.

The results generated by the 3SCOPE model enable the EM to make several decisions and generate its optimal supply chain planning. For example, the 3SCOPE model recommends the EM to manufacture one eVTOL for customer $c1$ on the day 56 in AH. As given before, the EM would need $\theta$ number of eVTOL parts to make one eVTOL, so to start manufacturing one eVTOL for customer $c1$ on day 56, the EM should receive the required eVTOL parts by day 55. As given in Figure \ref{figSM1}, the 3SCOPE model provides that one eVTOL part $a$ from supplier $s2$ via mode ship, six eVTOL parts $b$ from supplier $s4$ via mode rail, four eVTOL parts $c$ from supplier $s5$ via mode rail, six eVTOL parts $d$ from supplier $s2$ via mode ship, and five eVTOL parts $e$ from supplier $s1$ via mode rail should be received by the EM on the day 55 to meet the manufacturing requirement at the beginning of day 56. The 3SCOPE model suggests that the EM should initiate an order with supplier $s2$ for eVTOL part $a$ on day 52, accounting for a two-day lead time of supplier $s2$ and an additional one-day transportation time via ship for the delivery of eVTOL parts $a$ from supplier $s2$ to the EM. The difference between the optimal order placement day and the optimal day for receiving the eVTOL parts reflects both the supplier's lead time and the transportation time needed to deliver the eVTOL parts to the EM's facility. Similarly, the EM should initiate orders with supplier $s4$ for eVTOL part $b$ on $49^{th}$ day, allowing for a four-day lead time and one-day transportation time via rail. Orders with supplier $s5$ for eVTOL part $c$ should be placed on day 50, factoring in a three-day lead time and two-day transportation time via rail. Orders with supplier $s2$ for eVTOL part $d$ should begin on Day 47, accounting for a four-day lead time and three-day transportation time via ship. Finally, orders with supplier $s1$ for eVTOL part $e$ should start on day 51, considering a two-day lead time and two-day transportation time via rail.

\begin{figure}[tb] 
    \centering
    \includegraphics[width=\textwidth]{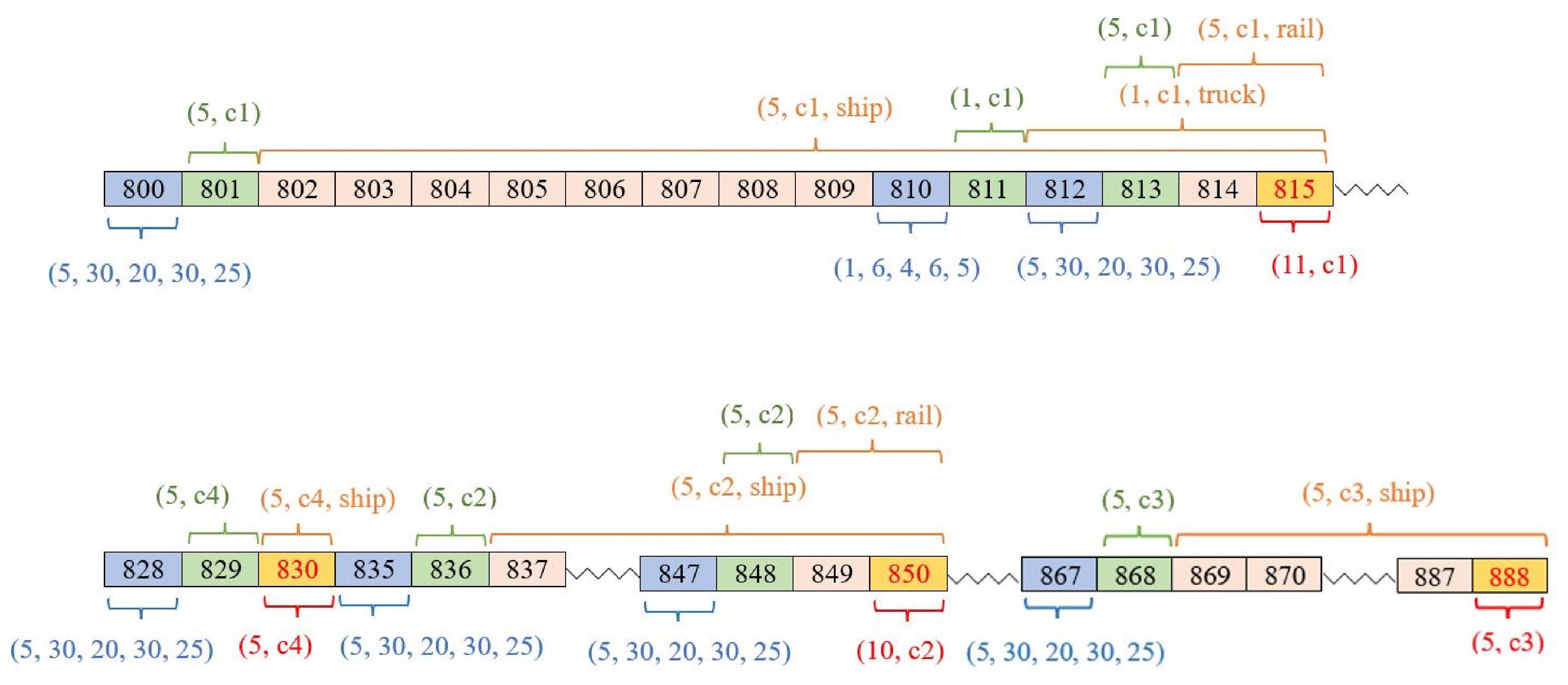}
    \caption{Optimal schedule of the front-end of the supply chain generated in the base case for Period 5, with blocks ordered by ascending days.}
    \label{schedule}
\end{figure}

An optimal schedule for the front-end of the supply chain of Period 5, which spans from day 721 to day 900, is displayed in Figure \ref{schedule}. Each blue block signifies the day when the EM should receive eVTOL parts from the suppliers and hold them in the EM's inventory at the beginning of a given day; each green block the day when the EM should manufacture the eVTOLs for the customers; orange blocks the time frame for the eVTOLs to be delivered to the customer via selected transportation mode; and each yellow block with the red font indicates the customer's deadline (by the end of this day, the customer is expected to receive the eVTOLs manufactured). The numbers written beneath the blue line represent eVTOL parts (a, b, c, d, e) received at the start of a given day. Above the green line, the first segment represents the number of eVTOLs to be manufactured at the beginning of a given day. The second segment indicates the customer for whom these eVTOLs are manufactured on that day, considering that the EM can manufacture a maximum of $R=5$ eVTOLs each day. Above the orange line, the text is divided into three parts. The first part indicates the number of eVTOLs transported to the customer, as mentioned in the second part. The third part specifies the transportation mode, with the length of the line representing the required days for delivery via that mode to the customer. Below the red line, the text is divided into two parts. The first part denotes the final number of eVTOLs pre-ordered by the customer mentioned in the second part. To illustrate how the schedule provides insights to the EM, an example is given for customer $c1$. As outlined in Table \ref{demand_deadline}, $D^1_l$ for customer $l=1$ is 11 eVTOLs with a deadline set for the end of day 815. As shown in Figure \ref{schedule}, the scheduling details for customer $c1$ are as follows: at the beginning of day 800, the EM is expected to receive eVTOL parts a, b, c, d, and e, in quantities of 5, 30, 20, 30, and 25, respectively, which should be held in the inventory. On day 801, the EM should manufacture five eVTOLs using these parts, which are then held in the inventory. At the beginning of day 802, the eVTOLs manufactured are dispatched to customer $c1$ via ship, which takes 14 days to reach the customer's location. Subsequently, the EM should receive eVTOL parts again on day 810 and day 812 to manufacture one eVTOL on day 811 and five eVTOLs on day 813, which are then delivered to customer $c1$ via truck (four days delivery time) and rail (two days delivery time), respectively. Consequently, by the end of day 815, the EM would successfully meet the deadline of customer $c1$, supplying a total of 11 eVTOLs. This schedule provides similar insights for other customers.

The total operating profits generated in AH in the base case from the 3SCOPE model and the benchmark models are presented in Table \ref{comparison1}. The total operating profits generated from the 3SCOPE model differ by approximately 3.35\% compared to the deterministic model due to the uncertainty involved in the supply chain planning of the EM. Next, when comparing the results generated by two other stochastic benchmark models, we observe that the 3SCOPE model outperforms the others. The 3SCOPE model yields profits approximately 5.72\% higher than the stochastic sequential model. This suggests that developing the integrated model is more efficient than running the independent models sequentially. The stochastic heuristic model also exhibits lower profitability when compared to the 3SCOPE model. The lower performance of the stochastic heuristic model is attributed to their simplified decision-making processes, which make locally optimal choices at each step. Focusing solely on lower supplier prices, without considering other factors in overall supply chain optimization, leads to increased inventory and transportation costs. This ultimately results in lower overall profit compared to the 3SCOPE model. A comparison of the runtime in seconds is provided in Table \ref{tab:runtime_comp} for two methods employed to solve the 3SCOPE model across varying numbers of scenarios. The runtime data for the method without multi-cut Benders decomposition is only available for up to 10 scenarios, highlighting its limitation in handling larger scenario sets with the given computer specifications. In contrast, the runtime data for the method with multi-cut Benders decomposition consistently performs well across a broader range of scenarios, indicating its scalability and efficiency. The runtime data indicates that the 3SCOPE model without the multi-cut Benders decomposition method performs better when addressing small-scale problems with an uncertain environment that can be modeled with fewer scenarios. Conversely, the 3SCOPE model with the multi-cut Benders decomposition method is more effective for solving large-scale problems involving a highly uncertain environment that needs to be modeled with a larger number of scenarios.

\begin{table} [tb]
\begin{spacing}{1}
  \centering
  \caption{Comparison of total operating profits (million \$) generated in AH from different models (case 1: base case; case 2A: quality sensitivity analysis with $\kappa = 0.15$; case 2B: quality sensitivity analysis with $\kappa = 0.25$; case 3A: analysis of eVTOL supply chain disruptions due to optimal supplier sanctions; case 3B: analysis of eVTOL supply chain disruptions due to raw material scarcity of eVTOL parts; case 3C: analysis of eVTOL supply chain disruptions due to delay in supplier lead times; and case 4: analysis of eVTOL manufacturing to resemble automotive manufacturing.} \label{comparison1}
  \begin{tabular}{cccccc}
    \toprule
    \textbf{Numerical Case} & \textbf{Deterministic Model} & \textbf{3SCOPE model} & \textbf{Stochastic Sequential Model} & \textbf{Stochastic Heuristic Model} \\
    \midrule
    Base & 51.69 & 50.01 & 47.30 & 46.22 \\
    2A & 50.68 & 49.04 & 47.78 & 45.55 \\
    2B & 52.71 & 51.01 & 49.45 & 49.42 \\
    3A & 52.01 & 49.44 & 45.77 & 42.23 \\
    3B &  44.98 & 41.56 & 36.14 & 34.28 \\
    3C & 38.96 & 34.80 & 29.94 & 27.08 \\
    4 & 45.10 & 42.90 & 40.61 & 38.46 \\
    \bottomrule
  \end{tabular}
\end{spacing}
\end{table}

\begin{table}[t]
\begin{spacing}{1}
\centering
\caption{Comparison of runtime in seconds for two methods applied to solve the 3SCOPE model.}
\label{tab:runtime_comp}
\begin{tabular}{@{}ccccccccc@{}}
\toprule
\textbf{Number of Scenarios} & 3  & 10 & 25 & 40 & 55 & 70 & 85 & 100 \\
\midrule
\textbf{Without multi-cut Benders decomposition}  & 23.02 & 4340.32 & - & - & - & - & - & - \\
\textbf{With multi-cut Benders decomposition}    & 25.04 & 52.10 & 150.59 & 231.91 & 262.06 & 364.26 & 535.46 & 1054.87 \\
\bottomrule
\end{tabular}
\end{spacing}
\end{table}

\subsubsection{Numerical Case 2: Analyzing the Impact of Quality Sensitivity Parameter}

In this numerical case, we perform an analysis by varying $\kappa$, which is assumed to be 0.2 in the base case. As discussed in Section \ref{quality_o}, $\kappa$ affects $d^2_{l}$ for customer $l$. The values of $d^2_{l}$ for different $\kappa$ values are listed in Table \ref{sens_kappa}. $d^2_{l}$ varies compared to the numbers listed in Table \ref{demand_deadline} for each $\kappa$ over the five periods in AH. The result indicates that as $\kappa$ increases, the effect of $q$ becomes more pronounced on $d^2_{l}$ for customer $l$. The amount of increase in $d^2_{l}$ is reflected in the values of the final number of eVTOL pre-orders. The total operating profits generated by the EM from different models for $\kappa$ values of 0.15 and 0.25, corresponding to numerical cases 2A and 2B, respectively, are listed in Tables \ref{comparison1}. Similar to the base case, the 3SCOPE model generates profits differing by approximately 3.34\% from the deterministic model due to supply chain planning uncertainty. The 3SCOPE model outperforms other stochastic models and proves more efficient than the stochastic sequential and heuristic models. We observe that the total operating profits generated in the base case are approximately 2\% higher than the EM's total operating profits in numerical case 2A. Furthermore, the EM's total operating profits are expected to be approximately 2.3\% higher in numerical case 2B compared to the base case. These insights demonstrate that the EM's total operating profits would increase with greater sensitivity in customer demand and higher quality of eVTOLs.



\begin{table}[b!]
\begin{spacing}{1}
\centering
\caption{Endogenous additional number of eVTOLs pre-ordered by the customers across different $\kappa$ values.}
\label{sens_kappa}
\begin{tabular}{c|cccc|cccc|cccc}
\toprule
\textbf{Period} & \multicolumn{4}{c|}{\textbf{$\mathbf{\kappa = 0.15}$}} & \multicolumn{4}{c|}{\textbf{$\mathbf{\kappa = 0.2}$}} & \multicolumn{4}{c}{\textbf{$\mathbf{\kappa = 0.25}$}} \\
\textbf{} & \textbf{c1} & \textbf{c2} & \textbf{c3} & \textbf{c4} & \textbf{c1} & \textbf{c2} & \textbf{c3} & \textbf{c4} & \textbf{c1} & \textbf{c2} & \textbf{c3} & \textbf{c4} \\
\midrule
1 & 0 & 0 & 0 & 0 & 0 & 0 & 0 & 0 & 0 & 0 & 0 & 0 \\
2 & 0 & 0 & 0 & 0 & 0 & 0 & 0 & 0 & 1 & 0 & 0 & 0 \\
3 & 0 & 0 & 0& 0 & 1 & 1 & 1 & 1 & 2 & 1 & 1 & 1 \\
4 & 1 & 1 & 1 & 0 & 2 & 1 & 1 & 1 & 2 & 2 & 1 & 1 \\
5 & 2& 1 & 1 & 1 & 3 & 2& 1 & 1 & 4 & 3 & 2 & 2 \\
\bottomrule
\end{tabular}

\end{spacing}
\end{table}

\subsubsection{Numerical Case 3: Analyzing eVTOL Supply Chain Disruptions}

\paragraph{Case 3A: Optimal Supplier Sanctioned}

\begin{figure}[tb] 
    \centering
    \includegraphics[width=0.5\textwidth]{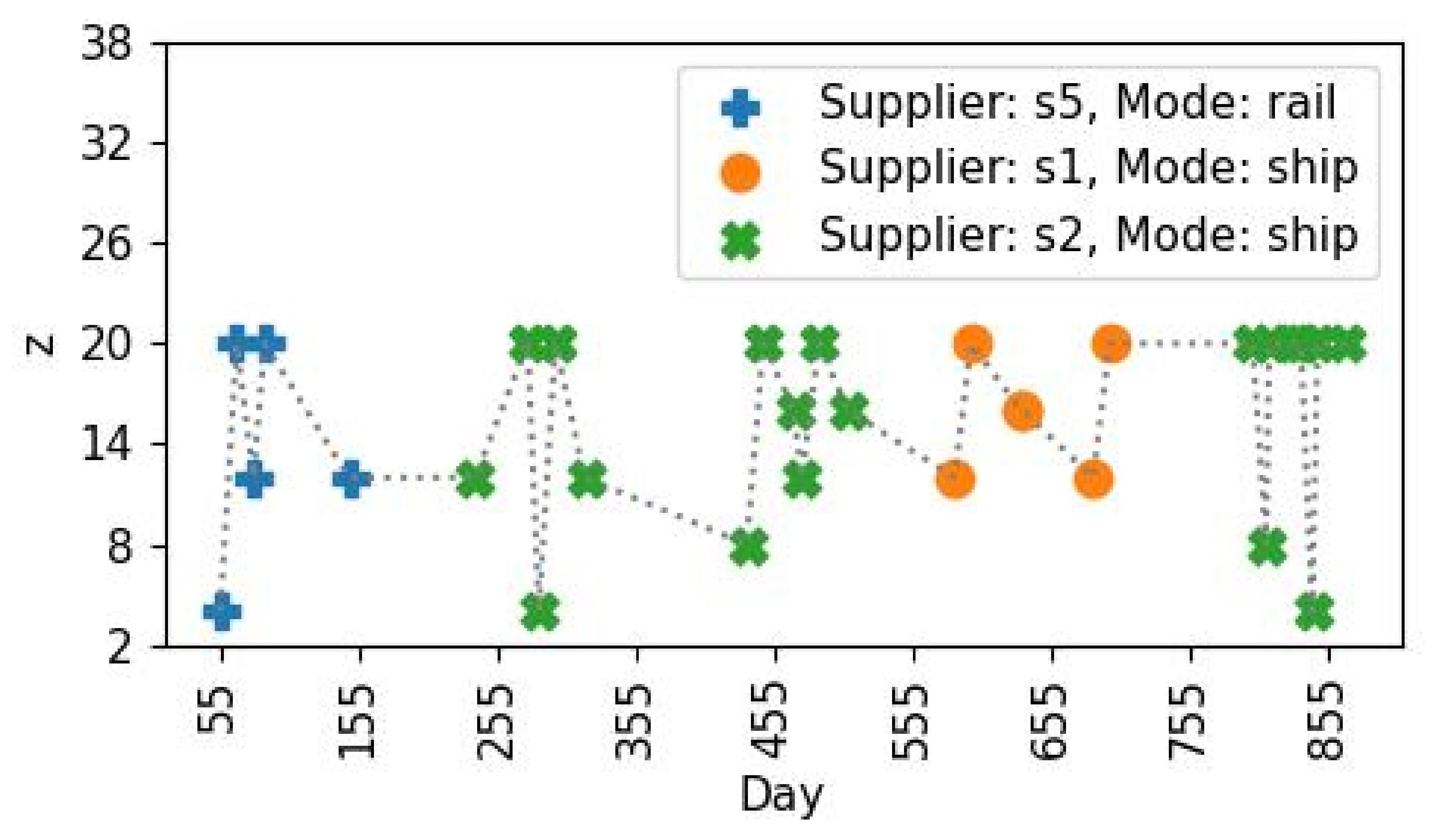}
    \caption{Optimal number of eVTOL parts $\boldsymbol{c}$ received when supplier $s5$ is sanctioned due to supply chain disruptions.}
    \label{partc_z_s5_remove}
\end{figure}

In the EM supply chain planning problem, it is crucial to consider the geopolitical stability of a supplier's region. This is because some regions may have less geopolitical stability than others. Political conflicts and sanctions between the countries where the EM and the supplier are located can disrupt the delivery of eVTOL parts during such sanctions. Conflicts can also arise within the U.S. states, leading to difficulties in delivering eVTOL products, requiring the EM to remove a supplier from the list if such sanctions or conflicts occur at any point of AH. In this case, we conduct an analysis to assess how sanctions affect the EM's supply chain. To do so, we select eVTOL part $c$ (battery), which is a critical part for eVTOL manufacturing. As shown in Figure \ref{figSM1c}, supplier $s5$ is chosen for eVTOL part $c$ throughout AH in the base case. However, we assume that, starting from day 181 of AH, the EM would face a conflict to select supplier $s5$. This conflict is unforeseen at the beginning of AH, and the EM would remove this supplier from their list at the beginning of the second iteration of RHA. The 3SCOPE model, as presented in Figure \ref{partc_z_s5_remove}, indicates that supplier $s5$ is chosen for $F^P$. However, for the subsequent periods, suppliers $s1$ and $s2$ are selected, which demonstrates that the 3SCOPE model can adapt to the removal of a supplier due to a conflict that arises in the middle of AH.

\begin{figure}[b!] 
    \centering
    \includegraphics[width=0.9\textwidth]{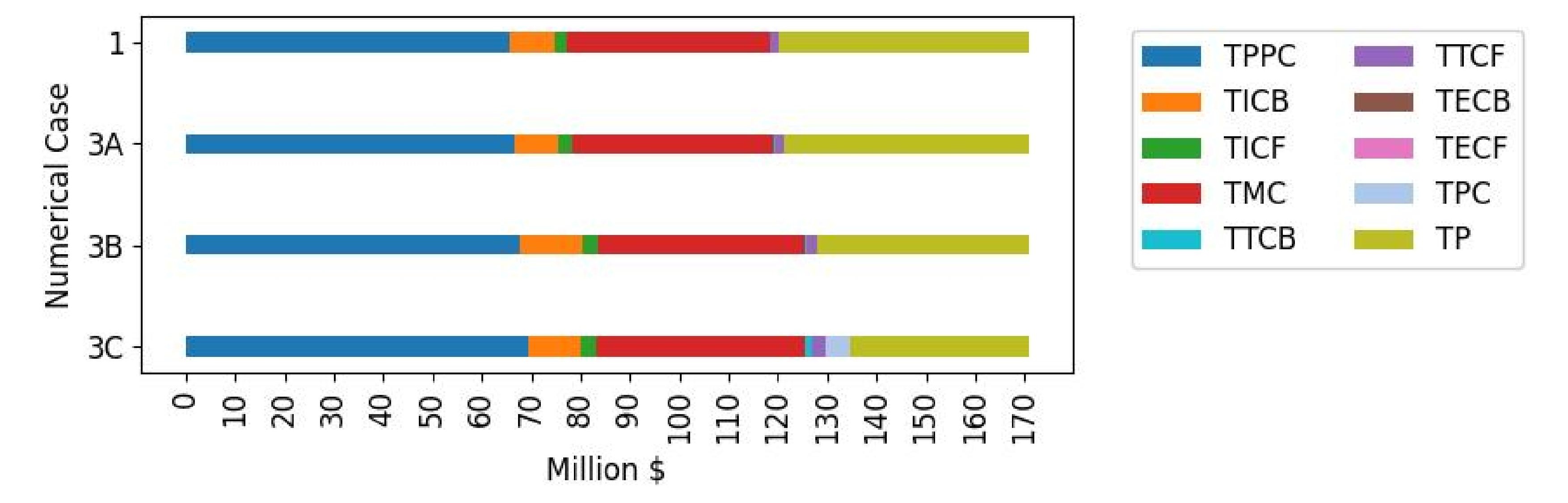}
    \caption{Breakdown of costs and profits for various numerical cases.}
    \label{cost_profit_fig}
\end{figure}

As listed in Table \ref{comparison1}, the total operating profits generated by the 3SCOPE model are approximately 5.2\% lower than the results generated by the deterministic model in case 3A. This suggests that the supply chain planning for the EM becomes more uncertain in case 3A compared to the base case. Comparing among the stochastic models, the 3SCOPE model generates higher profits than the stochastic benchmark models in this case, providing insights into its greater capability of performing under uncertainties compared to other stochastic benchmark models. Additionally, the profits in the table show a decrease of approximately 1.76\% compared to the total operating profits generated in the base case. Figure \ref{cost_profit_fig} shows the breakdown of total costs incurred in AH and total operating profits generated in AH from the 3SCOPE model. In this figure, TPPC refers to the total eVTOL part procurement cost, TICB the total inventory cost for holding eVTOL parts, TICF the total inventory cost for holding manufactured eVTOLs, TMC the total eVTOL manufacturing cost, TTCB the total transportation cost in the back-end of the supply chain, TTCF the total transportation cost in the front-end of the supply chain, TECB the total emission cost in the back-end of the supply chain, TECF the total emission cost in the front-end of the supply chain, TPC the total penalty cost, and TP the total operating profit. This figure depicts that the decrease in TP, despite the same revenue across all cases, is attributed to the increase in TPPC. Due to the higher prices offered by suppliers $s1$ and $s2$ compared to $s5$, as listed in Table \ref{suppliers}, TPPC in numerical case 3A is higher than in the base case.

\paragraph{Case 3B: Raw Material Scarcity}

When raw materials required for manufacturing eVTOL parts become scarce, suppliers may be unable to produce these parts, leading to their unavailability in the market. Consequently, the EM would face challenges in selecting suppliers for eVTOL parts due to disruptions in the AAM market. For example, we consider the scarcity of lithium, a key material used in manufacturing batteries for eVTOLs. This scarcity would result in suppliers being unable to produce and provide a sufficient number of batteries required for eVTOL manufacturing. In this example, we assume that this disruption occurs during Period 3 of AH. At the beginning of the second iteration of RHA, if the EM observes that the disruption would be in $S^P$ for eVTOL parts, the EM would update the input parameters for $S^P$ and run the 3SCOPE model for the second iteration. The total operating profits generated by the 3SCOPE model, as given in Table \ref{comparison1}, show an approximate 8.2\% decline compared to the deterministic model in case 3B. This observation suggests an increase in supply chain uncertainty for the EM in numerical case 3B relative to the base case. Upon comparison with other stochastic models, the 3SCOPE model demonstrates higher profit generation which indicates that the 3SCOPE model performs better in this case compared to alternative stochastic benchmark models. Furthermore, the table illustrates a decrease of approximately 14.86\% in profits compared to those generated in the base case. This is because all the costs, as presented Figure \ref{cost_profit_fig}, increased compared to the base case. The 3SCOPE model provides a solution by procuring a higher number of eVTOL parts in $F^P$ compared to the base case and manufacturing the eVTOLs in advance for $S^P$ customers in response to potential disruptions in $S^P$. TICB and TICF experience the most significant increase compared to other cost components, as the EM would need to hold a higher number of eVTOL parts and eVTOLs manufactured in its inventory. In this case, TP decreases by 20\% compared to the TP of the base case, but the 3SCOPE model ensures that the EM can fulfill customer pre-orders even in the midst of disruptions, without delaying deadlines.

\paragraph{Case 3C: Supply Chain Delay}

Considering the same example as in numerical case 3B, another reason for the scarcity of batteries could be the limited sourcing of battery materials, with major player countries serving as the primary suppliers. If geopolitical conflicts arise between these countries and the US, where the EM is located, the EM cannot directly procure batteries from these suppliers. Instead, the EM would need to wait for secondary suppliers from countries without conflict. These secondary suppliers would procure from primary suppliers and sell at a higher price to the EM. Consequently, the lead time for obtaining eVTOL parts increases. If the EM does not have enough batteries in inventory to continue manufacturing eVTOLs and fulfill customer pre-orders for $F^P$ in the third iteration, the EM would need to extend PH for this iteration to allow time for receiving batteries and manufacturing eVTOLs for that period. If the EM could foresee this disruption one period in advance within the planning horizon, the 3SCOPE model would have sufficient time to make proactive decisions for both $F^P$ and $S^P$, as observed in the case of 3B. However, in the same example, if the EM cannot foresee the disruption beforehand, the EM would only become aware of it after the second iteration when the model parameters would be updated for the third iteration of RHA. In this situation, the EM would not have enough time to procure and receive the scarce batteries on time to meet the customer deadline for Period 3 of AH and would need to extend the PH in the third iteration of RHA. We assume that this disruption occurs in Periods 3 and 4 in AH. Consequently, the length of PH in the third iteration of RHA is extended to 380 days from 360 days, with Periods 3 and 4 each having 190 days. As the disruption no longer occurs in Period 5, starting from the fifth iteration, the EM would once again set the length of the period to 360 days. The start and end days for each period in numerical case 3C are listed in Table \ref{period_years_3c}, with AH ranging from day 1 to day 920. In numerical case 3C, we consider that the lead times of suppliers are delayed in the order $s3 < s2 < s1 < s4$. In this case, we assume that the earlier the supplier can provide batteries, the more they charge. The primary supplier $s5$ is excluded from the set of suppliers. Based on their lead times, manufacturing time, transportation time at both ends of the supply chain, the EM would calculate the range of delays and renegotiate with customers and set new deadlines for the pre-orders. Information for all the customers regarding previous deadlines, new deadlines, delivery delays of eVTOL pre-orders, and penalties incurred by the EM due to the disruption in numerical case 3C are listed in Table \ref{delay_penalty_case3C}. The penalties are calculated from the data provided in Table \ref{tab:delay_penalty}.


\begin{table}[t]
\begin{spacing}{1}
\centering
\caption{Period start and end days over AH in numerical case 3C.}
\label{period_years_3c}
\begin{tabular}{lccccc}
\toprule
\textbf{} & \textbf{Period 1} & \textbf{Period 2} & \textbf{Period 3} & \textbf{Period 4} & \textbf{Period 5} \\ \midrule
\textbf{Start Day} & 1 & 181 & 361 & 551 & 741 \\
\textbf{End Day} & 180 & 360 & 550 & 740 & 920 \\

\bottomrule
\end{tabular}
\end{spacing}
\end{table}

\begin{table}[t]
\begin{spacing}{1}
\centering
\caption{Impact of disruption in numerical case 3C on customer deadlines, delivery delays, and penalties.}
\label{delay_penalty_case3C}
\begin{tabular}{ccccc}
\toprule
\textbf{Customers}                   & \textbf{c1} & \textbf{c2} & \textbf{c3} & \textbf{c4} \\
\midrule
$\mathbf{T^{d^1}_{l}}$                             & 490 & 450 & 490 & 510 \\
$\mathbf{T^{d^\Upsilon}_{l}}$                      & 542 & 530 & 540 & 545 \\
$\mathbf{|T^{d^\Upsilon}_{l} - T^{d^1}_{l}|}$    & 52 & 80 & 50 & 35 \\
$\mathbf{D^1_l}$                                  & 8 & 7 & 4 & 4 \\
\textbf{Contract Value (million \$)}                   & 12 & 10.5 & 6 & 6 \\
\textbf{Penalty (million \$)}                    & 1.23 & 3.12 & 0.56 & 0.23 \\
\bottomrule
\end{tabular}
\end{spacing}
\end{table}

\begin{figure}[b!] 
    \centering
    \includegraphics[width=0.9\textwidth]{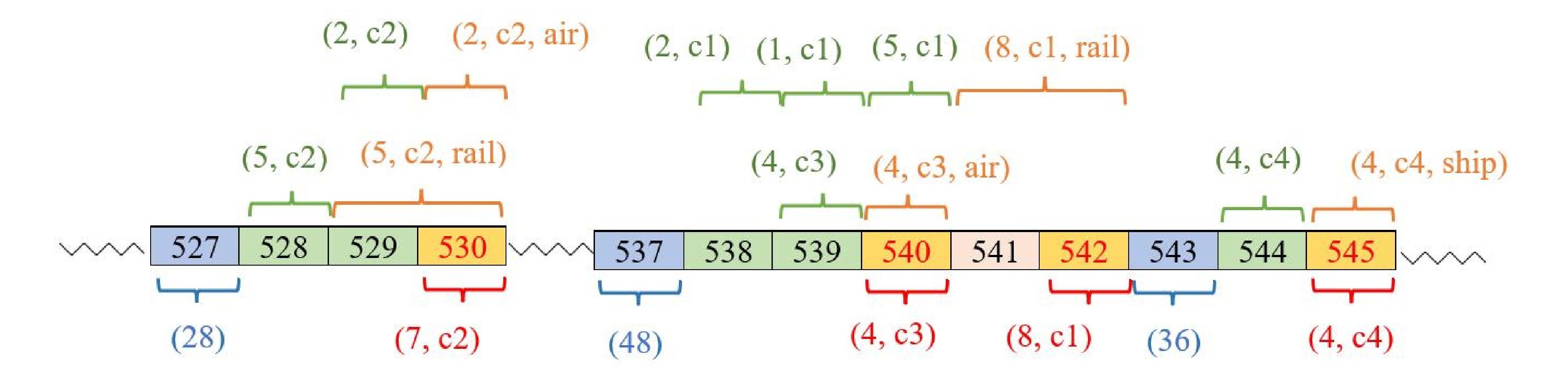}
    \caption{Optimal schedule of the front-end of the supply chain generated in numerical case 3C for Period 3, with blocks ordered by ascending days.}
    \label{schedule_3C}
\end{figure}

The manufacturing schedule for numerical case 3C is displayed in Figure \ref{schedule_3C}. The colors, blocks, and other notations in Figure \ref{schedule_3C} hold the same meaning as described for Figure \ref{schedule}, except that in Figure \ref{schedule_3C}, the numbers in blue represent only the number of batteries. This is because the scarcity of batteries in this numerical case delays the entire supply chain. The EM receives 28 batteries from supplier $s3$ at a higher price to manufacture seven eVTOLs for customer $c2$ and deliver it by the new deadline day 530. As no other supplier can provide these batteries by day 527, the EM should buy from supplier $s3$, even if their price is higher than others. However, by day 537, supplier $s2$ can deliver batteries at a lower price than supplier $s3$, so the EM should receive the batteries from them instead of supplier $s3$. Using these batteries, the EM should manufacture 12 eVTOLs to meet the deadlines of customer $c3$ and customer $c1$. On day 543, the EM receives 36 batteries, although it only needs 16 batteries to manufacture 4 eVTOLs for customer $c4$, who is the only one left for $F^P$. The reason for this surplus is that these batteries are delivered from supplier $s1$, who offers the lowest price among the suppliers available during this disruption. Additionally, there is uncertainty regarding the continuation of this disruption in the next period. Therefore, the EM orders batteries in advance for manufacturing eVTOLs for $S^P$, along with the order placed in $F^P$ to supplier $s1$.

As shown in Table \ref{comparison1}, the total operating profits generated by the 3SCOPE model are approximately 12\% lower than those of the deterministic model in case 3C. Analyzing the deviations in profits across different numerical cases, it becomes evident that the more severe the disruptions, the greater the uncertainty the EM would face in its supply chain planning. When compared with other stochastic models, the 3SCOPE model generates higher profits, indicating better performance in this case than the alternative stochastic benchmark models. Additionally, the table shows a profit decrease of about 28.03\% compared to the base case. This is because all the costs, as presented Figure \ref{cost_profit_fig}, increased significantly compared to the base case, even more than in case 3B. The increase in TPPC arises because the EM is forced to procure batteries from supplier $s3$ at a higher cost compared to the primary supplier $s5$ and other secondary suppliers $s1$ and $s2$ to fulfill eVTOL pre-orders for customer $c2$. This is due to the longer lead times associated with suppliers $s1$ and $s2$ in comparison to $s3$, and the EM cannot directly procure batteries from supplier $s5$ in numerical case 3C. The batteries can be received from supplier $s3$, located outside the US, via two modes: air (transportation time 1 day) and ship (transportation time 22 days). Because TPC for delaying the deadline of the customer is significantly higher in the aviation industry, this higher penalty pushes the EM to obtain the batteries as early as possible in this situation. The EM selects the air transportation to receive the batteries needed for fulfilling the demand of $F^P$, despite its higher cost, because it ensures delivery within just 1 day. Selecting the cheaper ship mode would result in a 22-day delay, causing TPS to increase more than the cost savings achieved in TTCB if the ship mode is chosen instead of the air transportation mode.


\subsubsection{Numerical Case 4: Reshaping eVTOL Manufacturing to Resemble Automotive Manufacturing}

Aircraft manufacturers typically receive pre-orders well in advance, enabling them to plan their manufacturing processes over an extended period. This longer planning horizon is feasible because traditional aircraft production is complex and does not operate at the high manufacturing rates seen in the automotive industry \cite{timmis2020aircraft}. In contrast, the automotive industry relies on high-volume manufacturing processes \cite{wada2020evolution}. Automotive manufacturers produce a significant number of vehicles within a relatively short time frame to meet the high demand in the automotive market. In the base case, we assume that PH for eVTOL manufacturing is similar to that of aircraft manufacturing, where in each PH, we consider $T$ = 360 days. However, the eVTOL manufacturing process would need to adapt to shorter planning horizons, as the demand for manufacturing eVTOLs is expected to align more closely with automotive manufacturing over time \cite{garrow2021urban, cohen2021urban, garrow2019survey}. In this numerical case, we conduct an analysis with a shorter PH to observe a situation where AAM operators place pre-orders for larger quantities of eVTOLs, mirroring the manufacturing rates in the automotive industry. We consider the same initial number of eVTOL pre-orders as shown in Table \ref{demand_deadline}, but we change the length of PH to $T$ = 90 days to perform this analysis, consisting of 45 days in each period. 

\begin{figure}[tb] 
    \centering
    \includegraphics[width=0.5\textwidth]{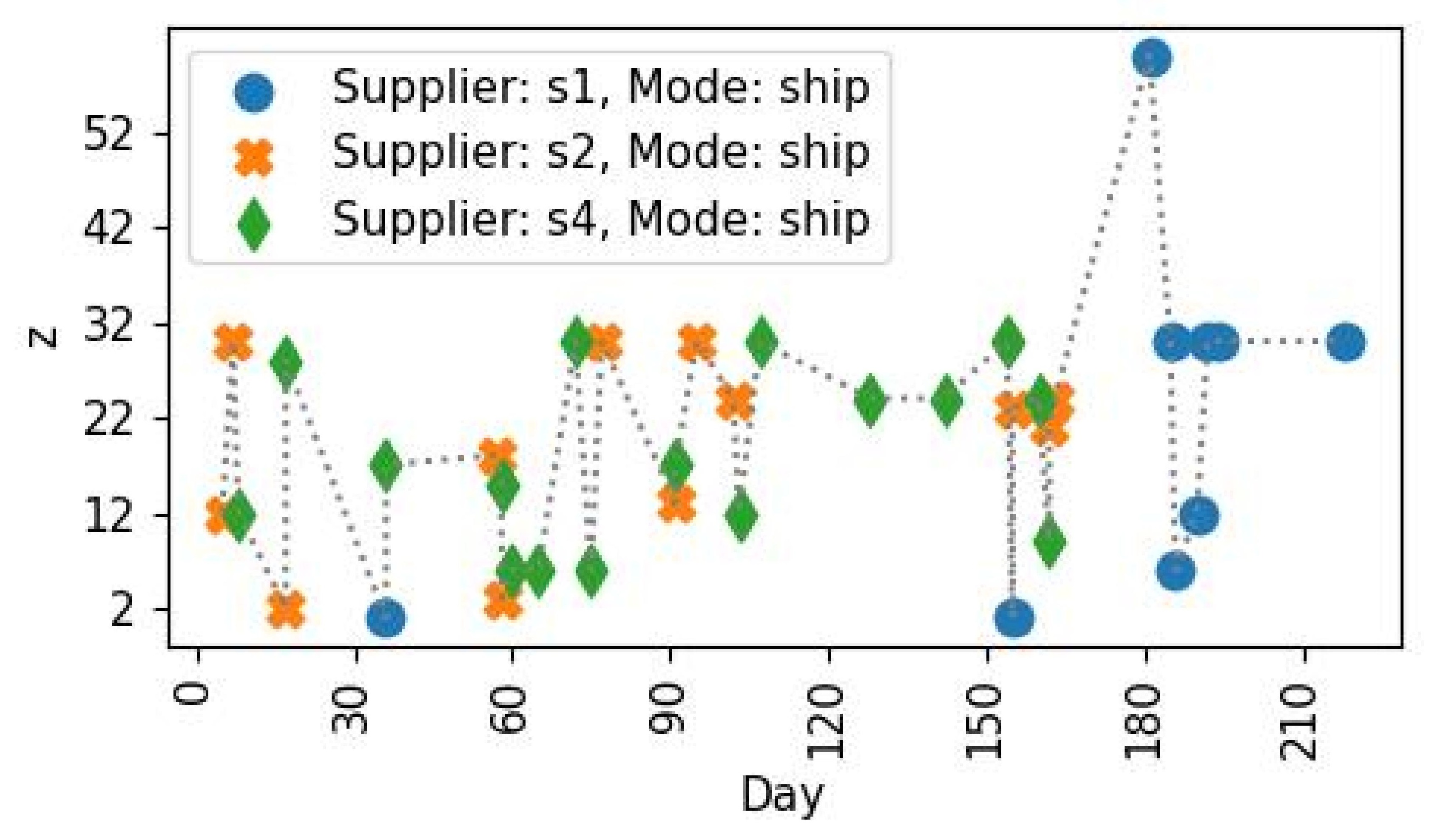}
    \caption{Optimal number of eVTOL parts $\boldsymbol{d}$ received from selected suppliers via the optimal transportation mode, considering that eVTOL manufacturing resembles automotive manufacturing.}
    \label{figSM2d}
\end{figure}

As the final number of eVTOL pre-orders remains the same as in the base case, the number of necessary eVTOL parts also remains the same. Consequently, results such as total revenue, total eVTOL manufacturing cost, total inventory cost for holding eVTOL parts, and total inventory cost for holding manufactured eVTOLs remain consistent with the base case. However, other costs, such as total eVTOL parts procurement cost, total transportation cost in the back-end of the supply chain, total emission cost in the back-end of the supply chain, total transportation cost in the front-end of the supply chain, and total emission cost in the front-end of the supply chain, are higher in comparison to the base case. This is because the chosen suppliers for the eVTOL parts and the transportation mode to deliver the eVTOL parts to the EM are different from those in the base case. The 3SCOPE model selects more than one supplier for a given type of eVTOL part. For example, in Figure \ref{figSM2d}, three suppliers are chosen, whereas in Figure \ref{figSM1d}, only one supplier ($s2$) is selected. In the base case, only supplier $s2$ is selected because it is capable of providing the required number of eVTOL parts throughout the longer PH in that case, considering both its capacity and lead time. However, in this case, only supplier $s2$ is not capable of supplying the required number of eVTOL parts throughout the shorter PH. Therefore, the 3SCOPE model selects supplier $s1$ and supplier $s4$ along with supplier $s2$. In the base case, the EM can afford to wait and procure eVTOL parts at a lower price from suppliers who have comparatively longer lead times. This is because these suppliers are capable of supplying the parts within the required time frame when the EM has a longer PH. However, in a shorter PH with the same number of eVTOLs to be manufactured, the EM cannot afford to procure eVTOL parts solely from these suppliers. In this situation, the EM should also prioritize selecting suppliers $s1$ and $s4$, who can provide the necessary number of eVTOL parts, even if their prices are higher than those offered by supplier $s2$. Since the EM would need to procure from multiple suppliers, both delivery and emission costs also increase compared to the base case.




The total operating profits generated by the EM in numerical case 4 from different models are listed in Tables \ref{comparison1}. An important insight is that, in contrast to the total operating profits generated in the base case, the total operating profits generated in numerical case 4 are lower. This difference in total operating profits arises from the increase in the total cost, as discussed previously. However, the EM can still benefit from this situation if the increased demand for eVTOLs prompts the EM to adapt the manufacturing process, similar to automotive manufacturing. With a 225-day AH, the EM could potentially generate a maximum of \$44.10 million in total operating profit in numerical case 4. In contrast, in the base case where the AH spans 900 days, the maximum total operating profit is \$55.69 million. This implies that in the base case, which has an analysis horizon four times longer than that considered in numerical case 4, the EM can only obtain approximately 13\% higher total operating profit. Therefore, when the demand of eVTOLs would be higher, similar to the automotive industry, then the EM can generate more total operating profit. However, the total operating profit generated from the 3SCOPE model is approximately \$42.90, which is 5.12\% lower than the deterministic model. In numerical case 4, we observe a higher difference in these two solutions than in the base case, where it is 3.35\%. This suggests that while the total operating profit generated would be higher in numerical case 4 compared to the base case, the EM would encounter more uncertainties in its supply chain planning in numerical case 4. As observed in previous cases, the 3SCOPE model also generates higher profits than the stochastic benchmark models in this case. The results indicate that the 3SCOPE model performs better than other existing stochastic models.

\section{Managerial Insights} \label{insights}

This study on EM's supply chain planning offers vital managerial insights to inform practical decision-making. In this section, we provide key insights derived from this study to aid managers in the emerging eVTOL manufacturing industry in optimizing their supply chain planning. 

\begin{itemize}[left=0pt]

\item The 3SCOPE model integrates functions such as supplier selection, procurement planning, manufacturing scheduling, inventory management, transportation mode selection, and quality control into a single framework. Utilizing the 3SCOPE model, managers can streamline decision-making processes, decrease manual intervention, minimize errors, and improve overall efficiency and profitability. 

\item The 3SCOPE model assists managers in identifying optimal suppliers capable of delivering the required quantity of eVTOL parts in alignment with the optimal manufacturing schedule. The managers should consider factors such as price, capacity, quality, and geopolitical risks associated with suppliers of eVTOL parts. Additionally, the results generated from the 3SCOPE model guide managers on when to place orders for how many eVTOL parts from suppliers, as well as when to expect the delivery of these parts. The managers should follow this integrated schedule and receive new pre-orders for eVTOLs accordingly.

\item In the aviation industry, where holding costs are typically higher, the managers should minimize the duration that eVTOL parts and eVTOLs remain in inventory. Our 3SCOPE model suggests implementing a just-in-time policy for the managers, considering various costs in integrated supply chain planning and meeting deadlines on time. However, to mitigate supply chain disruptions, such as potential raw material shortages and supplier unavailability, which could lead to excess penalty costs for delayed delivery to customers, the managers should procure eVTOL parts ahead of time. In this case, the managers should balance the cost of increased inventory against the risk of supply chain disruptions.

\item To continuously update manufacturing schedules based on real-time data, the managers should adopt a rolling horizon approach, as we did to run the 3SCOPE model in this study. This approach enables dynamic adjustments to production plans in response to changes in customer demand, supplier capacities, and other uncertainties. By regularly updating the supply chain plan based on the latest data, the managers can swiftly respond to unforeseen challenges and design a dynamic supply chain plan.

\item  Our study emphasizes the importance of balancing cost, delivery speed, and environmental impact in transportation. The managers should collaborate with logistics partners offering sustainable options and consider emission costs in the overall analysis for improved decision-making. 

\item The findings of this study underscore the positive influence of the quality of eVTOL parts and eVTOLs on customer demand and profitability. Therefore, the managers should prioritize conducting regular audits and performance evaluations of suppliers and supply chain processes. These measures can effectively uphold eVTOL quality and minimize defects, thereby earning customer trust and maximizing profitability.

\item The managers should invest in advanced methods such as a multi-cut Benders decomposition to tackle large-scale supply chain planning problems more effectively. This investment significantly improves the scalability and efficiency of stochastic supply chain models, enabling the managers to manage extensive data and complex scenarios without compromising on the accuracy and timeliness of decisions.

\item This This study demonstrates that our model is applicable to both traditional low-volume aviation manufacturing and high-volume eVTOL production. Consequently, it can aid managers in analyzing their supply chain networks and scaling up eVTOL production to meet the growing demand.

\end{itemize}

\section{Conclusions and Future Works}\label{sec5}

The potential of eVTOLs to revolutionize AAM has generated widespread interest across industry, government, and academia. It is necessary to study the supply chain challenges that the EM is expected to face and develop the 3SCOPE model to enable them to meet customer demands on time by optimizing their supply chain planning. Therefore, this study introduces the 3SCOPE model, a two-stage stochastic optimization model, designed with the primary objective of maximizing the operating profit of the EM. This model functions as a decision-making process for the EM, integrating various challenges in its supply chain planning under various uncertainties including prices offered by suppliers for eVTOL parts, capacities of suppliers, eVTOL manufacturing costs, and customer demand for eVTOLs. To account for the changing market dynamics in the eVTOL industry, the 3SCOPE model is adapted using a rolling horizon approach. Furthermore, a multi-cut Benders decomposition method is utilized to run the 3SCOPE model, addressing the large-scale supply chain planning problem expected to be faced by the EM. Our 3SCOPE model helps EMs integrate and optimize their supply chain planning by addressing several critical questions, such as which suppliers to select for each eVTOL part; how many eVTOL parts to procure from selected suppliers; when to place orders with selected suppliers for these eVTOL parts; when and in what quantities to receive the ordered eVTOL parts; when and how many eVTOLs to manufacture to meet customer pre-orders within given deadlines; which transportation mode to select; when to schedule deliveries, and how many eVTOLs to dispatch in each delivery of eVTOL parts and eVTOLs in the supply chain; how many eVTOL parts and eVTOLs to hold in inventory; and how many eVTOL pre-orders from customers are expected to be affected by the quality of eVTOLs. The results show that the EM would generate approximately \$50 million in total operating profit within a 900-day analysis horizon. The findings indicate that the quality of eVTOL parts and eVTOLs has a proportional impact on both the eVTOL pre-orders and the total operating profit of the EM. Our model also performs well in high-volume eVTOL manufacturing, where the EM needs to match the manufacturing rate of automotive manufacturing. It enables the EM to develop a supply chain plan that can quickly adapt to changes in demand, supplier capacities, and other factors, allowing them to manufacture commercial eVTOLs on a large scale.

Furthermore, we analyze three potential disruptions to understand their impact on the supply chain planning of the EM. In the first disrupted case, where the EM needs to replace an optimal supplier due to geopolitical sanctions, the findings indicate that the 3SCOPE model promptly handles the disruption by selecting another supplier, ensuring smooth supply chain planning throughout the analysis horizon. However, this leads to a 1.76\% decrease in total operating profit, as the replacement supplier provides eVTOL parts at higher prices. In the second disrupted case, the suppliers face scarcity of raw materials, resulting in their inability to produce eVTOL parts for the EM. Consequently, the EM faces a shortage of suppliers in the market. However, our 3SCOPE model effectively manages the disruption without missing customer deadlines. The 3SCOPE model suggests that the EM should procure a sufficient quantity of eVTOL parts and manufacture an adequate number of eVTOLs well in advance to fulfill pre-orders. However, this leads to increased inventory costs and a 14.86\% decrease in total operating profit. Nevertheless, it successfully fulfills customer pre-orders by the deadlines and avoids incurring penalty costs for late deliveries. In the third disrupted case, where the EM's supply chain is delayed due to prolonged lead times of suppliers in the market, the EM faces penalties for missing customer deadlines. This results in a significant 28.03\% decrease in profit. However, the insights also demonstrate that our 3SCOPE model is capable of providing optimal supply chain planning, allowing the EM to make optimal decisions and achieve the highest possible profits even in this worst disrupted case. We consider benchmark models to compare the performance of the 3SCOPE model. Across all numerical cases, the 3SCOPE model consistently yields higher profits, indicating better performance compared to alternative stochastic benchmark models. These insights will help managers in the eVTOL manufacturing industry make informed decisions and optimize their supply chain planning. This study assists the managers in maximizing profitability by offering the capability to handle a range of challenges simultaneously. Additionally, the insights from the study enable managers to respond to market changes and supply chain disruptions in the dynamic eVTOL industry.



This study provides several interesting directions for further research. Firstly, we aim to consider reverse logistics, which involves the process of handling returned goods, recycling, or disposing of products in the supply chain planning of the EM. This would enable the EM to incorporate more sustainable approaches into the 3SCOPE model. Secondly, we intend to strategize the distribution of parts and materials essential for maintenance and repair services to the customers, as well as for general supplies at vertiports. We also aim to determine the locations of the EM for manufacturing eVTOL and the locations of depots for providing after-repair and maintenance services. Thirdly, we aim to explore the EM's adoption of technology for various types of eVTOLs, focusing on those expected to be in demand in the future. Fourthly, we plan to incorporate artificial intelligence to enhance the scalability of our model, outperforming the stochastic heuristic and stochastic sequential models benchmarked in this study. Lastly, our focus is to consider multiple competitors in the eVTOL manufacturing market and develop strategies for the EM to succeed within this competitive landscape.



\noindent \textbf{Declaration of Interests:} None. 
\end{spacing}

\begin{spacing}{0.65}
\bibliography{sample} 

\begin{thebibliography}{138}
\newcommand{\enquote}[1]{``#1''}
\providecommand{\natexlab}[1]{#1}
\providecommand{\url}[1]{\texttt{#1}}
\providecommand{\urlprefix}{URL }
\expandafter\ifx\csname urlstyle\endcsname\relax
  \providecommand{\doi}[1]{\discretionary{}{}{}https://doi.org/#1}\else
  \providecommand{\doi}[1]{\discretionary{}{}{}\urlstyle{rm}\url{https://doi.org/#1}}\fi

\bibitem[{{Federal Aviation Administration}(2022)}]{FAAConOpsv2}
{Federal Aviation Administration}, \enquote{UTM Concept of Operations Version 2.0 (UTM ConOps v2.0),} \url{https://www.faa.gov/researchdevelopment/trafficmanagement/utm-concept-operations-version-20-utm-conops-v20}, 2022.

\bibitem[{Dulia et~al.(2021{\natexlab{a}})Dulia, Sabuj, and Shihab}]{dulia2021benefits}
Dulia, E.~F., Sabuj, M.~S., and Shihab, S.~A., \enquote{Benefits of advanced air mobility for society and environment: A case study of ohio,} \emph{Applied Sciences}, Vol.~12, No.~1, 2021{\natexlab{a}}, p. 207.
\newblock \doi{https://doi.org/10.3390/app12010207}.

\bibitem[{Garrow et~al.(2021)Garrow, German, and Leonard}]{garrow2021urban}
Garrow, L.~A., German, B.~J., and Leonard, C.~E., \enquote{Urban air mobility: A comprehensive review and comparative analysis with autonomous and electric ground transportation for informing future research,} \emph{Transportation Research Part C: Emerging Technologies}, Vol. 132, 2021, p. 103377.
\newblock \doi{https://doi.org/10.1016/j.trc.2021.103377}.

\bibitem[{Vieira et~al.(2019)Vieira, Silva, and Bravo}]{vieira2019electric}
Vieira, D.~R., Silva, D., and Bravo, A., \enquote{Electric VTOL aircraft: the future of urban air mobility (background, advantages and challenges),} \emph{International Journal of Sustainable Aviation}, Vol.~5, No.~2, 2019, pp. 101--118.
\newblock \doi{https://doi.org/10.1504/IJSA.2019.101746}.

\bibitem[{Maes(accessed November 2023)}]{VFS_eVTOL_companies}
Maes, R.~A., \enquote{eVTOLs: Who Will Win the Race to Market?} \url{https://www.avbuyer.com/articles/business-aircraft-development-and-certification/evtols-who-will-win-the-race-to-market-113653#:~:text=According%20to%20the%20Vertical%20Flight%20Society%20%28VFS%29%2C%20the,come%20from%20347%20entities%20spread%20across%2048%20countries.}, accessed November 2023.

\bibitem[{{Joby Aviation}(accessed October 24, 2023)}]{joby_design}
{Joby Aviation}, \enquote{Electric Aerial Ridesharing,} \url{https://www.jobyaviation.com/}, accessed October 24, 2023.

\bibitem[{{Eve Air Mobility}(accessed October 24, 2023)}]{eve}
{Eve Air Mobility}, \enquote{Mobility reimagined,} \url{https://eveairmobility.com/}, accessed October 24, 2023.

\bibitem[{Rothfeld et~al.(2021)Rothfeld, Fu, Bala{\'c}, and Antoniou}]{rothfeld2021potential}
Rothfeld, R., Fu, M., Bala{\'c}, M., and Antoniou, C., \enquote{Potential urban air mobility travel time savings: An exploratory analysis of munich, paris, and san francisco,} \emph{Sustainability}, Vol.~13, No.~4, 2021, p. 2217.
\newblock \doi{https://doi.org/10.3390/su13042217}.

\bibitem[{Yazan and Moniz(2021)}]{yazan2021technology}
Yazan, A., and Moniz, A., \enquote{Technology Assessment of eVTOL Air Transportation System: The Positive Impacts (Potential Benefits),} \emph{International Journal of Advanced Research}, Vol.~9, 2021, pp. 171--180.
\newblock \doi{10.21474/IJAR01/12291}.

\bibitem[{Goyal et~al.(2018{\natexlab{a}})Goyal, Reiche, Fernando, Serrao, Kimmel, Cohen, and Shaheen}]{goyal2018urban}
Goyal, R., Reiche, C., Fernando, C., Serrao, J., Kimmel, S., Cohen, A., and Shaheen, S., \enquote{Urban Air Mobility (UAM) Market Study,} Tech. Rep. HQ-E-DAA-TN65181, Booz Allen Hamilton, USA, and University of California, Berkeley, California, USA, November 2018{\natexlab{a}}.

\bibitem[{Goyal et~al.(2018{\natexlab{b}})Goyal, Reiche, Fernando, Serrao, Kimmel, Cohen, and Shaheen}]{hasan2019urban}
Goyal, R., Reiche, C., Fernando, C., Serrao, J., Kimmel, S., Cohen, A., and Shaheen, S., \enquote{Urban Air Mobility (UAM) Market Study,} Tech. Rep. HQ-E-DAA-TN63717, Booz Allen Hamilton, USA, and University of California, Berkeley, California, USA, October 2018{\natexlab{b}}.

\bibitem[{Doo et~al.(2021)Doo, Pavel, Didey, Hange, Diller, Tsairides, Smith, Bennet, Bromfield, and Mooberry}]{doo2021nasa}
Doo, J.~T., Pavel, M.~D., Didey, A., Hange, C., Diller, N.~P., Tsairides, M.~A., Smith, M., Bennet, E., Bromfield, M., and Mooberry, J., \enquote{NASA Electric Vertical Takeoff and Landing (eVTOL) Aircraft Technology for Public Services – A White Paper,} Tech. Rep. WBS: 330693.01.20.01.02, Ames Research Center, Mountain View, California, United States, August 2021.

\bibitem[{{Deloitte Insights}(accessed November 2023)}]{delloitte0}
{Deloitte Insights}, \enquote{The elevated future of mobility: What’s next on the horizon?} \url{https://www2.deloitte.com/us/en/insights/focus/future-of-mobility/evtol-elevated-future-of-mobility-summary.html}, accessed November 2023.

\bibitem[{Beamon(1998)}]{beamon1998supply}
Beamon, B.~M., \enquote{Supply chain design and analysis:: Models and methods,} \emph{International journal of production economics}, Vol.~55, No.~3, 1998, pp. 281--294.
\newblock \doi{https://doi.org/10.1016/S0925-5273(98)00079-6}.

\bibitem[{GARRETT-GLASER(accessed January 2024)}]{scm_motiv2}
GARRETT-GLASER, B., \enquote{Industry seeks solutions to eVTOL supply chain challenges,} \url{https://verticalmag.com/news/industry-seeks-solutions-supply-chain-challenges/}, accessed January 2024.

\bibitem[{S.~Rangarajan et~al.(2022)S.~Rangarajan, Sunddararaj, Sudhakar, Shiva, Subramaniam, Collins, and Senjyu}]{s2022lithium}
S.~Rangarajan, S., Sunddararaj, S.~P., Sudhakar, A., Shiva, C.~K., Subramaniam, U., Collins, E.~R., and Senjyu, T., \enquote{Lithium-ion batteries—The crux of electric vehicles with opportunities and challenges,} \emph{Clean Technologies}, Vol.~4, No.~4, 2022, pp. 908--930.
\newblock \doi{https://doi.org/10.3390/cleantechnol4040056}.

\bibitem[{{Flight Global}(accessed November 18, 2023)}]{Flight_Global}
{Flight Global}, \enquote{Suppliers seek to renegotiate loss-making contracts with Airbus and Boeing,} \url{https://www.flightglobal.com/airframers/suppliers-seek-to-renegotiate-loss-making-contracts-with-airbus-and-boeing/155881.article}, accessed November 18, 2023.

\bibitem[{Saxena(accessed November 23, 2023)}]{scm_motiv1}
Saxena, S., \enquote{Supply Chain Considerations to Support Urban Air Mobility eVTOLs,} \url{https://www.mrodigestforums.com/supply-chain-considerations-to-support-urban-air-mobility-evtols/}, accessed November 23, 2023.

\bibitem[{Wang et~al.(2024)Wang, Wang, and Chen}]{wang2024impact}
Wang, Z., Wang, H., and Chen, X., \enquote{The impact of delayed fixed-price payment in the decentralised project supply chain,} \emph{International Journal of Systems Science: Operations \& Logistics}, Vol.~11, No.~1, 2024, p. 2308584.
\newblock \doi{https://doi.org/10.1080/23302674.2024.2308584}.

\bibitem[{Jain et~al.(2009)Jain, Benyoucef, and Deshmukh}]{jain2009strategic}
Jain, V., Benyoucef, L., and Deshmukh, S., \enquote{Strategic supplier selection: some emerging issues and challenges,} \emph{International Journal of Logistics Systems and Management}, Vol.~5, No. 1-2, 2009, pp. 61--88.
\newblock \doi{https://doi.org/10.1504/IJLSM.2009.021645}.

\bibitem[{Cheaitou and Khan(2015)}]{cheaitou2015integrated}
Cheaitou, A., and Khan, S.~A., \enquote{An integrated supplier selection and procurement planning model using product predesign and operational criteria,} \emph{International Journal on Interactive Design and Manufacturing (IJIDeM)}, Vol.~9, No.~3, 2015, pp. 213--224.
\newblock \doi{https://doi.org/10.1007/s12008-015-0280-5}.

\bibitem[{Qin et~al.(2020)Qin, Ma, Chan, and Khan}]{qin2020scenario}
Qin, Y., Ma, H.-L., Chan, F.~T., and Khan, W.~A., \enquote{A scenario-based stochastic programming approach for aircraft expendable and rotable spare parts planning in MRO provider,} \emph{Industrial Management \& Data Systems}, Vol. 120, No.~9, 2020, pp. 1635--1657.
\newblock \doi{https://doi.org/10.1108/IMDS-03-2020-0131}.

\bibitem[{Esteso et~al.(2023)Esteso, Peidro, Mula, and D{\'\i}az-Madro{\~n}ero}]{esteso2023reinforcement}
Esteso, A., Peidro, D., Mula, J., and D{\'\i}az-Madro{\~n}ero, M., \enquote{Reinforcement learning applied to production planning and control,} \emph{International Journal of Production Research}, Vol.~61, No.~16, 2023, pp. 5772--5789.
\newblock \doi{https://doi.org/10.1080/00207543.2022.2104180}.

\bibitem[{Avotra and Nawaz(2023)}]{climate_transportation}
Avotra, A. A. R.~N., and Nawaz, A., \enquote{Asymmetric impact of transportation on carbon emissions influencing SDGs of climate change,} \emph{Chemosphere}, Vol. 324, 2023, p. 138301.
\newblock \doi{https://doi.org/10.1016/j.chemosphere.2023.138301}.

\bibitem[{Xu and Xu(2022)}]{xu2022assessing}
Xu, B., and Xu, R., \enquote{Assessing the role of environmental regulations in improving energy efficiency and reducing CO2 emissions: Evidence from the logistics industry,} \emph{Environmental Impact Assessment Review}, Vol.~96, 2022, p. 106831.
\newblock \doi{https://doi.org/10.1016/j.eiar.2022.106831}.

\bibitem[{Liu et~al.(2022)Liu, Qian, Hu, Shang, Li, Zhao, Zhao, and Han}]{liu2022government}
Liu, Z., Qian, Q., Hu, B., Shang, W.-L., Li, L., Zhao, Y., Zhao, Z., and Han, C., \enquote{Government regulation to promote coordinated emission reduction among enterprises in the green supply chain based on evolutionary game analysis,} \emph{Resources, Conservation and Recycling}, Vol. 182, 2022, p. 106290.
\newblock \doi{https://doi.org/10.1016/j.resconrec.2022.106290}.

\bibitem[{Dulia et~al.(2021{\natexlab{b}})Dulia, Ali, Garshasbi, and Kabir}]{green_sc}
Dulia, E.~F., Ali, S.~M., Garshasbi, M., and Kabir, G., \enquote{Admitting risks towards circular economy practices and strategies: An empirical test from supply chain perspective,} \emph{Journal of Cleaner Production}, Vol. 317, 2021{\natexlab{b}}, p. 128420.
\newblock \doi{https://doi.org/10.1016/j.jclepro.2021.128420}.

\bibitem[{{U.S. Department of Transportation}(accessed January 13, 2024)}]{EPA}
{U.S. Department of Transportation}, \enquote{GHG Reduction Strategies,} \url{https://www.transportation.gov/sustainability/climate/ghg-reduction-strategies}, accessed January 13, 2024.

\bibitem[{Baxter(accessed January 10, 2024)}]{regulations}
Baxter, K., \enquote{How New Regulations Impact Freight Transportation in 2023,} \url{https://blog.intekfreight-logistics.com/new-regulations-freight-transportation-2023}, accessed January 10, 2024.

\bibitem[{Klyde et~al.(2023)Klyde, Jones, Kotikalpudi, Lotterio et~al.}]{klyde2023developing}
Klyde, D., Jones, M., Kotikalpudi, A., Lotterio, M., et~al., \enquote{Developing Means of Compliance for eVTOL Vehicles: Phase II Final Report,} Tech. rep., United States. Department of Transportation. Federal Aviation Administration, 2023.
\newblock \doi{https://doi.org/10.21949/1528235}.

\bibitem[{Lina(2022)}]{lina2022improving}
Lina, R., \enquote{Improving Product Quality and Satisfaction as Fundamental Strategies in Strengthening Customer Loyalty,} \emph{AKADEMIK: Jurnal Mahasiswa Ekonomi \& Bisnis}, Vol.~2, No.~1, 2022, pp. 19--26.
\newblock \doi{https://doi.org/10.37481/jmeb.v2i1.245}.

\bibitem[{Silva et~al.(published online 24 Jun. 2018)Silva, Johnson, Solis, Patterson, and Antcliff}]{silva2018vtol}
Silva, C., Johnson, W., Solis, E., Patterson, M., and Antcliff, K., \enquote{VTOL urban air mobility concept vehicles for technology development,} \emph{2018 Aviation Technology, Integration, and Operations Conference}, published online 24 Jun. 2018, p. 3847.
\newblock \doi{https://doi.org/10.2514/6.2018-3847}.

\bibitem[{Thipphavong et~al.(2018)Thipphavong, Apaza, Barmore, Battiste, Burian, Dao, Feary, Go, Goodrich, Homola, Idris, Kopardekar, Lachter, Neogi, Ng, Osequera-Lohr, Patterson, and Verma}]{thipphavong2018urban}
Thipphavong, D., Apaza, R., Barmore, B., Battiste, V., Burian, B., Dao, Q., Feary, M., Go, S., Goodrich, K., Homola, J., Idris, H., Kopardekar, P., Lachter, J., Neogi, N., Ng, H., Osequera-Lohr, R., Patterson, M., and Verma, S., \enquote{Urban air mobility airspace integration concepts and considerations,} \emph{2018 Aviation Technology, Integration, and Operations Conference}, 2018, p. 3676.
\newblock \doi{https://doi.org/10.2514/6.2018-3676}.

\bibitem[{Hasan(Jun. 2018)}]{hasan2018urban}
Hasan, S., \enquote{Urban air mobility market study,} Tech. rep., Crown Consulting, Inc., Washington, DC, Jun. 2018.

\bibitem[{{Federal Aviation Administration (FAA)}(accessed August 16, 2022)}]{FAA2020}
{Federal Aviation Administration (FAA)}, \enquote{Concept of Operations v2.0, Unmanned aircraft System (UAS) Traffic Management (UTM),} \url{https://www.faa.gov/uas/research_development/traffic_management/media/UTM_ConOps_v2.pdf}, accessed August 16, 2022.

\bibitem[{Shihab et~al.(2019)Shihab, Wei, Ramirez, Mesa-Arango, and Bloebaum}]{shihab2019schedule}
Shihab, S. A.~M., Wei, P., Ramirez, D. S.~J., Mesa-Arango, R., and Bloebaum, C., \enquote{By schedule or on demand?-a hybrid operation concept for urban air mobility,} \emph{AIAA Aviation 2019 Forum}, 2019, p. 3522.
\newblock \doi{https://doi.org/10.2514/6.2019-3522}.

\bibitem[{Shihab et~al.(published online 8 Jun. 2020)Shihab, Wei, Shi, and Yu}]{shihab2020optimal}
Shihab, S. A.~M., Wei, P., Shi, J., and Yu, N., \enquote{Optimal evtol fleet dispatch for urban air mobility and power grid services,} \emph{AIAA Aviation 2020 Forum}, published online 8 Jun. 2020, p. 2906.
\newblock \doi{https://doi.org/10.2514/6.2020-2906}.

\bibitem[{Sabziyan~Varnousfaderani et~al.(2023)Sabziyan~Varnousfaderani, Shihab, and Dulia}]{deep_dispatch}
Sabziyan~Varnousfaderani, E., Shihab, S. A.~M., and Dulia, E.~F., \enquote{Deep-Dispatch: A Deep Reinforcement Learning-Based Vehicle Dispatch Algorithm for Advanced Air Mobility,} \url{https://doi.org/10.13140/RG.2.2.14850.86725}, November 2023.

\bibitem[{Dulia and Shihab(2024{\natexlab{a}})}]{dulia2024designing}
Dulia, E.~F., and Shihab, S.~A., \enquote{Designing a Surveillance Sensor Network with Information Clearinghouse for Advanced Air Mobility,} \emph{Sensors}, Vol.~24, No.~3, 2024{\natexlab{a}}, p. 803.
\newblock \doi{https://doi.org/10.3390/s24030803}.

\bibitem[{{Dulia, Esrat Farhana and Shihab, Syed A. M.}(2023)}]{open_framework_standards}
{Dulia, Esrat Farhana and Shihab, Syed A. M.}, \enquote{Open Framework Standards for Combined Aircraft Sensor Network for the State of Ohio to Detect and Track Lower Altitude Aircraft: Cost-Benefit Analysis,} Tech. Rep. Agreement No.: 36496, PID: 114242, SJN: 136337, U.S. Department of Transportation, December 2023.
\newblock \doi{10.13140/RG.2.2.17752.92166}.

\bibitem[{Dulia and Shihab(2024{\natexlab{b}})}]{PPP}
Dulia, E.~F., and Shihab, S. A.~M., \enquote{How to Negotiate with Private Investors for Advanced Air Mobility Infrastructure? An Analysis of Public Private Partnerships using Game Theory,} \url{https://doi.org/10.13140/RG.2.2.25387.86566/1}, March 2024{\natexlab{b}}.

\bibitem[{Mattei et~al.(2024)Mattei, de~Alteriis, Conte, Carotenuto, Rufino, De~Maio, and Accardo}]{mattei2024improving}
Mattei, F., de~Alteriis, G., Conte, C., Carotenuto, V., Rufino, G., De~Maio, A., and Accardo, D., \enquote{Improving Radar Detection of Drones and Air Mobility Systems in Urban Areas,} \emph{AIAA SCITECH 2024 Forum}, 2024, p. 2060.
\newblock \doi{https://doi.org/10.1109/ICCE59016.2024.10444151}.

\bibitem[{Kopardekar(accessed March 2024{\natexlab{a}})}]{nasa_scm}
Kopardekar, P., \enquote{Manufacturing \& Supply Chain for Advanced Air Mobility,} \url{https://ntrs.nasa.gov/api/citations/20220007440/downloads/5.17.22%20aero-supplychain-manufacturing-AvWeek.pdf}, accessed March 2024{\natexlab{a}}.

\bibitem[{{National Aeronautics and Space Administration}(accessed October 24, 2023)}]{AAM_supplychain}
{National Aeronautics and Space Administration}, \enquote{Building Resilient Aerospace Supply Chain,} \url{https://nari.arc.nasa.gov/sites/default/files/attachments/FACTsheetv7.pdf}, accessed October 24, 2023.

\bibitem[{Kopardekar(accessed March 2024{\natexlab{b}})}]{nasa_scm1}
Kopardekar, P., \enquote{Aerospace Supply Chain and Manufacturing,} \url{https://nari.arc.nasa.gov/sites/default/files/attachments/1-pk-feb4-5-SupplyChainManagement.pdf}, accessed March 2024{\natexlab{b}}.

\bibitem[{Mateo-Forn{\'e}s et~al.(2023)Mateo-Forn{\'e}s, Soto-Silva, Gonz{\'a}lez-Araya, Pl{\`a}-Aragon{\`e}s, and Solsona-Tehas}]{mateo2023managing}
Mateo-Forn{\'e}s, J., Soto-Silva, W., Gonz{\'a}lez-Araya, M.~C., Pl{\`a}-Aragon{\`e}s, L.~M., and Solsona-Tehas, F., \enquote{Managing quality, supplier selection, and cold-storage contracts in agrifood supply chain through stochastic optimization,} \emph{International Transactions in Operational Research}, Vol.~30, No.~4, 2023, pp. 1901--1930.
\newblock \doi{https://doi.org/10.1111/itor.13069}.

\bibitem[{Almeida and Concei{\c{c}}{\~a}o(2021)}]{almeida2021decomposition}
Almeida, J. F. d.~F., and Concei{\c{c}}{\~a}o, S.~V., \enquote{A decomposition approach for the two-stage stochastic supply network planning in light of the rolling horizon practice,} \emph{Pesquisa Operacional}, Vol.~41, 2021, p. e234451.
\newblock \doi{https://doi.org/10.1590/0101-7438.2021.041s1.00234451}.

\bibitem[{{Future Flight}(accessed March 14, 2024)}]{delivery_joby_evtol}
{Future Flight}, \enquote{Joby eVTOL - Complete performance Data,} \url{https://www.futureflight.aero/aircraft-program/joby-evtol}, accessed March 14, 2024.

\bibitem[{Hussain and Silver(accessed March 2024)}]{lineberger2021advanced}
Hussain, A., and Silver, D., \enquote{Advanced Air Mobility: Can the United States afford to lose the race?} \url{https://www2.deloitte.com/us/en/insights/industry/aerospace-defense/advanced-air-mobility.html}, accessed March 2024.

\bibitem[{Moses and {\AA}hlstr{\"o}m(2008)}]{moses2008problems}
Moses, A., and {\AA}hlstr{\"o}m, P., \enquote{Problems in cross-functional sourcing decision processes,} \emph{Journal of Purchasing and Supply Management}, Vol.~14, No.~2, 2008, pp. 87--99.
\newblock \doi{https://doi.org/10.1016/j.pursup.2007.11.003}.

\bibitem[{Tao et~al.(2018)Tao, Cheng, Qi, Zhang, Zhang, and Sui}]{tao2018digital}
Tao, F., Cheng, J., Qi, Q., Zhang, M., Zhang, H., and Sui, F., \enquote{Digital twin-driven product design, manufacturing and service with big data,} \emph{The International Journal of Advanced Manufacturing Technology}, Vol.~94, 2018, pp. 3563--3576.
\newblock \doi{https://doi.org/10.1007/s00170-017-0233-1}.

\bibitem[{Chiang et~al.(2012)Chiang, Kocabasoglu-Hillmer, and Suresh}]{chiang2012empirical}
Chiang, C.-Y., Kocabasoglu-Hillmer, C., and Suresh, N., \enquote{An empirical investigation of the impact of strategic sourcing and flexibility on firm's supply chain agility,} \emph{International Journal of Operations \& Production Management}, Vol.~32, No.~1, 2012, pp. 49--78.
\newblock \doi{https://doi.org/10.1108/01443571211195736}.

\bibitem[{Mocenco(2015)}]{mocenco2015supply}
Mocenco, D., \enquote{Supply Chain Features of the Aerospace Industry: Particular Case Airbus and Boeing,} \emph{Scientific Bulletin-Economic Sciences/Buletin Stiintific-Seria Stiinte Economice}, Vol.~14, No.~2, 2015.
\newblock \urlprefix\url{https://www.semanticscholar.org/paper/SUPPLY-CHAIN-FEATURES-OF-THE-AEROSPACE-INDUSTRY-AND-Mocenco/b7925fd5daab2fcb35b17e3db29ba38d58194ff4}.

\bibitem[{Kotabe et~al.(2023)Kotabe, Ayebale, and Murray}]{kotabe2023relationship}
Kotabe, M., Ayebale, D., and Murray, J.~Y., \enquote{Relationship multiplexity, multiple resource acquisition, and export performance of emerging-market firms,} \emph{Journal of International Management}, Vol.~29, No.~3, 2023, p. 101030.
\newblock \doi{https://doi.org/10.1016/j.intman.2023.101030}.

\bibitem[{Dweiri et~al.(2016)Dweiri, Kumar, Khan, and Jain}]{dweiri2016designing}
Dweiri, F., Kumar, S., Khan, S.~A., and Jain, V., \enquote{Designing an integrated AHP based decision support system for supplier selection in automotive industry,} \emph{Expert Systems with Applications}, Vol.~62, 2016, pp. 273--283.
\newblock \doi{https://doi.org/10.1016/j.eswa.2016.06.030}.

\bibitem[{Hashemi et~al.(2015)Hashemi, Karimi, and Tavana}]{hashemi2015integrated}
Hashemi, S.~H., Karimi, A., and Tavana, M., \enquote{An integrated green supplier selection approach with analytic network process and improved Grey relational analysis,} \emph{International Journal of Production Economics}, Vol. 159, 2015, pp. 178--191.
\newblock \doi{https://doi.org/10.1016/j.ijpe.2014.09.027}.

\bibitem[{Bills et~al.(2023)Bills, Fredericks, Sulzer, and Viswanathan}]{bills2023massively}
Bills, A., Fredericks, L., Sulzer, V., and Viswanathan, V., \enquote{Massively Distributed Bayesian Analysis of Electric Aircraft Battery Degradation,} \emph{ACS Energy Letters}, Vol.~8, No.~8, 2023, pp. 3578--3585.
\newblock \doi{https://doi.org/10.1021/acsenergylett.3c01216}.

\bibitem[{Benton~Jr(accessed on March 2024)}]{benton2020purchasing}
Benton~Jr, W., \enquote{Purchasing and supply chain management,} \url{https://books.google.com/books?hl=en&lr=&id=rrD_DwAAQBAJ&oi=fnd&pg=PT18&dq=Purchasing+and+supply+chain+management+benton&ots=QRUtt9bIxS&sig=aVh5YjFe4zk8he9zZBUJffpRFDw#v=onepage&q=Purchasing%20and%20supply%20chain%20management%20benton&f=false}, accessed on March 2024.

\bibitem[{Liang(2023)}]{liang2023production}
Liang, Y., \enquote{Production Scheduling Optimization of An Aviation Bearing Manufacturing Enterprise Based on Teaching-Learning-based Optimization,} \url{https://doi.org/http://dx.doi.org/10.54097/ajst.v6i2.9707}, 2023.

\bibitem[{Ma et~al.(2023)Ma, Huang, Hu, Chen, Qian, Deng, and Hua}]{ma2023multi}
Ma, H., Huang, X., Hu, Z., Chen, Y., Qian, D., Deng, J., and Hua, L., \enquote{Multi-objective production scheduling optimization and management control system of complex aerospace components: a review,} \emph{The International Journal of Advanced Manufacturing Technology}, Vol. 127, No. 11-12, 2023, pp. 4973--4993.
\newblock \doi{https://doi.org/10.1007/s00170-023-11707-4}.

\bibitem[{Bouzembrak et~al.(2011)Bouzembrak, Allaoui, Goncalves, and Bouchriha}]{inproceedings}
Bouzembrak, Y., Allaoui, H., Goncalves, G., and Bouchriha, H., \enquote{Distribution Supply Chain Design under Demand Uncertainty,} \emph{International Conference on Industrial Engineering and Systems Management, France}, 2011.
\newblock \urlprefix\url{https://www.researchgate.net/publication/293815989_Distribution_Supply_Chain_Design_under_Demand_uncertainty}.

\bibitem[{Chigbu and Nekhwevha(2021)}]{chigbu2021future}
Chigbu, B.~I., and Nekhwevha, F.~H., \enquote{The future of work and uncertain labour alternatives as we live through the industrial age of possible singularity: Evidence from South Africa,} \emph{Technology in Society}, Vol.~67, 2021, p. 101715.
\newblock \doi{https://doi.org/10.1016/j.techsoc.2021.101715}.

\bibitem[{Sharma and Sinha(2012)}]{sharma2012production}
Sharma, R., and Sinha, A.~K., \enquote{A production planning model using fuzzy neural network: a case study of an automobile industry,} \emph{International Journal of Computer Applications}, Vol. 975, 2012, p. 8887.
\newblock \doi{10.5120/5032-7183}.

\bibitem[{Khan et~al.(2023)Khan, Zaman, and Khan}]{khan2023relationship}
Khan, M.~I., Zaman, S.~I., and Khan, S.~A., \enquote{Relationship and impact of block chain technology and supply chain management on inventory management,} \emph{Blockchain Driven Supply Chain Management: A Multi-dimensional Perspective}, Springer, 2023, pp. 53--74.
\newblock \doi{https://doi.org/10.1007/978-981-99-0699-4_4}.

\bibitem[{Bhattacharya(accessed March 2024)}]{bhattacharya2021working}
Bhattacharya, H., \emph{Working capital management: Strategies and techniques (4th edition)}, PHI Learning Pvt. Ltd., accessed March 2024.
\newblock \urlprefix\url{https://books.google.com/books?id=jBw7EAAAQBAJ}.

\bibitem[{Derhami et~al.(2021)Derhami, Montreuil, and Bau}]{derhami2021assessing}
Derhami, S., Montreuil, B., and Bau, G., \enquote{Assessing product availability in omnichannel retail networks in the presence of on-demand inventory transshipment and product substitution,} \emph{Omega}, Vol. 102, 2021, p. 102315.
\newblock \doi{https://doi.org/10.1016/j.omega.2020.102315}.

\bibitem[{Singh and Verma(2018)}]{singh2018inventory}
Singh, D., and Verma, A., \enquote{Inventory management in supply chain,} \emph{Materials Today: Proceedings}, Vol.~5, No.~2, 2018, pp. 3867--3872.
\newblock \doi{https://doi.org/10.1016/j.matpr.2017.11.641}.

\bibitem[{Dillon(2019)}]{dillon2019study}
Dillon, A.~P., \emph{A study of the Toyota production system: From an Industrial Engineering Viewpoint}, Routledge, New York, 2019.
\newblock \doi{https://doi.org/10.4324/9781315136509}.

\bibitem[{Masoud and Mason(2016)}]{masoud2016integrated}
Masoud, S.~A., and Mason, S.~J., \enquote{Integrated cost optimization in a two-stage, automotive supply chain,} \emph{Computers \& Operations Research}, Vol.~67, 2016, pp. 1--11.
\newblock \doi{https://doi.org/10.1016/j.cor.2015.08.012}.

\bibitem[{Tomic et~al.(2012)Tomic, Spasojević-Brkić, and Klarin}]{article}
Tomic, B., Spasojević-Brkić, V., and Klarin, M., \enquote{Quality management system for the aerospace industry,} \emph{Journal of Engineering Management and Competitiveness}, Vol.~2, 2012.
\newblock \doi{10.5937/jemc1201011T}.

\bibitem[{Findlay and Harrison(2002)}]{FINDLAY200218}
Findlay, S., and Harrison, N., \enquote{Why aircraft fail,} \emph{Materials Today}, Vol.~5, No.~11, 2002, pp. 18--25.
\newblock \doi{https://doi.org/10.1016/S1369-7021(02)01138-0}, \urlprefix\url{https://www.sciencedirect.com/science/article/pii/S1369702102011380}.

\bibitem[{Roca et~al.(2017)Roca, Vaishnav, Morgan, Mendon{\c{c}}a, and Fuchs}]{roca2017risks}
Roca, J.~B., Vaishnav, P., Morgan, M.~G., Mendon{\c{c}}a, J., and Fuchs, E., \enquote{When risks cannot be seen: Regulating uncertainty in emerging technologies,} \emph{Research Policy}, Vol.~46, No.~7, 2017, pp. 1215--1233.
\newblock \doi{https://doi.org/10.1016/j.respol.2017.05.010}.

\bibitem[{{National Aeronautics and Space Administration}(accessed March 2024{\natexlab{a}})}]{q1}
{National Aeronautics and Space Administration}, \enquote{NASA Aviation Safety: Procurement Quality Assurance,} \url{https://sma.nasa.gov/vids/video-item/nasa-aviation-safety-procurement-quality-assurance#:~:text=Each%20center%20that%20operates%20aircraft%20currently%20handles%20parts,the%20acquisition%20of%20quality%20aircraft%20parts%20and%20services}, accessed March 2024{\natexlab{a}}.

\bibitem[{{National Aeronautics and Space Administration}(accessed March 2024{\natexlab{b}})}]{q2}
{National Aeronautics and Space Administration}, \enquote{Quality,} \url{https://sma.nasa.gov/sma-disciplines/quality}, accessed March 2024{\natexlab{b}}.

\bibitem[{{National Aeronautics and Space Administration}(accessed February 2024)}]{q3}
{National Aeronautics and Space Administration}, \enquote{NASA Procedural Requirements,} \url{https://nodis3.gsfc.nasa.gov/displayDir.cfm?t=NPR&c=8000&s=4B}, accessed February 2024.

\bibitem[{Alkahtani et~al.(2023)Alkahtani, Abidi, Obaid, and Alotaik}]{alkahtani2023modified}
Alkahtani, M., Abidi, M.~H., Obaid, H. S.~B., and Alotaik, O., \enquote{Modified Gannet Optimization Algorithm for Reducing System Operation Cost in Engine Parts Industry with Pooling Management and Transport Optimization,} \emph{Sustainability}, Vol.~15, No.~18, 2023, p. 13815.
\newblock \doi{https://doi.org/10.3390/su151813815}.

\bibitem[{Barzinpour and Taki(2018)}]{barzinpour2018dual}
Barzinpour, F., and Taki, P., \enquote{A dual-channel network design model in a green supply chain considering pricing and transportation mode choice,} \emph{Journal of Intelligent Manufacturing}, Vol.~29, 2018, pp. 1465--1483.
\newblock \doi{https://doi.org/10.1007/s10845-015-1190-x}.

\bibitem[{Rad and Nahavandi(2018)}]{rad2018novel}
Rad, R.~S., and Nahavandi, N., \enquote{A novel multi-objective optimization model for integrated problem of green closed loop supply chain network design and quantity discount,} \emph{Journal of cleaner production}, Vol. 196, 2018, pp. 1549--1565.
\newblock \doi{https://doi.org/10.1016/j.jclepro.2018.06.034}.

\bibitem[{Liberti and Kucherenko(2005)}]{liberti2005comparison}
Liberti, L., and Kucherenko, S., \enquote{Comparison of deterministic and stochastic approaches to global optimization,} \emph{International Transactions in Operational Research}, Vol.~12, No.~3, 2005, pp. 263--285.

\bibitem[{Pinsky and Karlin(2010)}]{pinsky2010introduction}
Pinsky, M., and Karlin, S., \emph{An introduction to stochastic modeling}, Academic press, 2010.

\bibitem[{Koirala et~al.(2022)Koirala, Van~Acker, D’hulst, and Van~Hertem}]{koirala2022hosting}
Koirala, A., Van~Acker, T., D’hulst, R., and Van~Hertem, D., \enquote{Hosting capacity of photovoltaic systems in low voltage distribution systems: A benchmark of deterministic and stochastic approaches,} \emph{Renewable and Sustainable Energy Reviews}, Vol. 155, 2022, p. 111899.

\bibitem[{Foroozesh et~al.(2021)Foroozesh, Jolai, Mousavi, and Karimi}]{foroozesh2021new}
Foroozesh, N., Jolai, F., Mousavi, S.~M., and Karimi, B., \enquote{A new fuzzy-stochastic compromise ratio approach for green supplier selection problem with interval-valued possibilistic statistical information,} \emph{Neural Computing and Applications}, Vol.~33, 2021, pp. 7893--7911.

\bibitem[{Kannegiesser and G{\"u}nther(2014)}]{kannegiesser2014sustainable}
Kannegiesser, M., and G{\"u}nther, H.-O., \enquote{Sustainable development of global supply chains—part 1: sustainability optimization framework,} \emph{Flexible Services and Manufacturing Journal}, Vol.~26, 2014, pp. 24--47.
\newblock \doi{https://doi.org/10.1007/s10696-013-9176-5}.

\bibitem[{Ghadimi et~al.(2023)Ghadimi, Aouam, and Uzsoy}]{ghadimi2023safety}
Ghadimi, F., Aouam, T., and Uzsoy, R., \enquote{Safety stock placement with market selection under load-dependent lead times,} \emph{IISE Transactions}, Vol.~55, No.~3, 2023, pp. 314--328.
\newblock \doi{https://doi.org/10.1080/24725854.2022.2074578}.

\bibitem[{Gruler et~al.(2018)Gruler, Panadero, de~Armas, P{\'e}rez, and Juan}]{gruler2018combining}
Gruler, A., Panadero, J., de~Armas, J., P{\'e}rez, J. A.~M., and Juan, A.~A., \enquote{Combining variable neighborhood search with simulation for the inventory routing problem with stochastic demands and stock-outs,} \emph{Computers \& Industrial Engineering}, Vol. 123, 2018, pp. 278--288.

\bibitem[{Fattahi and Govindan(2022)}]{fattahi2022data}
Fattahi, M., and Govindan, K., \enquote{Data-driven rolling horizon approach for dynamic design of supply chain distribution networks under disruption and demand uncertainty,} \emph{Decision Sciences}, Vol.~53, No.~1, 2022, pp. 150--180.
\newblock \doi{https://doi.org/10.1111/deci.12481}.

\bibitem[{Kaur et~al.(2020)Kaur, Singh, Garza-Reyes, and Mishra}]{kaur2020sustainable}
Kaur, H., Singh, S.~P., Garza-Reyes, J.~A., and Mishra, N., \enquote{Sustainable stochastic production and procurement problem for resilient supply chain,} \emph{Computers \& Industrial Engineering}, Vol. 139, 2020, p. 105560.

\bibitem[{Dnistran(accessed February 2024)}]{pre_order1}
Dnistran, I., \enquote{Alef Flying Car Racks Up 2,500 Pre-Orders Worth \$750M,} \url{https://insideevs.com/news/678657/alef-evtol-2500-pre-orders/}, accessed February 2024.

\bibitem[{Alock(accessed February 2024)}]{pre_order2}
Alock, C., \enquote{Airlines Place Biggest eVTOL Orders to Date as Vertical Goes Public,} \url{https://www.futureflight.aero/news-article/2021-06-11/airlines-place-biggest-evtol-orders-date-vertical-goes-public}, accessed February 2024.

\bibitem[{Motamed~Nasab and Li(2021)}]{motamed2021multistage}
Motamed~Nasab, F., and Li, Z., \enquote{Multistage adaptive stochastic mixed integer optimization under endogenous and exogenous uncertainty,} \emph{AIChE Journal}, Vol.~67, No.~10, 2021, p. e17333.
\newblock \doi{https://doi.org/10.1002/aic.17333}.

\bibitem[{Cummins and Danzon(1997)}]{cummins1997price}
Cummins, J.~D., and Danzon, P.~M., \enquote{Price, financial quality, and capital flows in insurance markets,} \emph{Journal of financial intermediation}, Vol.~6, No.~1, 1997, pp. 3--38.
\newblock \doi{https://doi.org/10.1006/jfin.1996.0205}.

\bibitem[{Paulley et~al.(2006)Paulley, Balcombe, Mackett, Titheridge, Preston, Wardman, Shires, and White}]{paulley2006demand}
Paulley, N., Balcombe, R., Mackett, R., Titheridge, H., Preston, J., Wardman, M., Shires, J., and White, P., \enquote{The demand for public transport: The effects of fares, quality of service, income and car ownership,} \emph{Transport policy}, Vol.~13, No.~4, 2006, pp. 295--306.
\newblock \doi{https://doi.org/10.1016/j.tranpol.2005.12.004}.

\bibitem[{S{\"o}derlund(1998)}]{soderlund1998customer}
S{\"o}derlund, M., \enquote{Customer satisfaction and its consequences on customer behaviour revisited: The impact of different levels of satisfaction on word-of-mouth, feedback to the supplier and loyalty,} \emph{International journal of service industry management}, Vol.~9, No.~2, 1998, pp. 169--188.
\newblock \doi{https://doi.org/10.1108/09564239810210532}.

\bibitem[{Masten(1984)}]{masten1984organization}
Masten, S.~E., \enquote{The organization of production: Evidence from the aerospace industry,} \emph{The Journal of Law and Economics}, Vol.~27, No.~2, 1984, pp. 403--417.
\newblock \doi{https://doi.org/10.1086/467071}.

\bibitem[{Dimitri et~al.(2006)Dimitri, Piga, and Spagnolo}]{dimitri2006handbook}
Dimitri, N., Piga, G., and Spagnolo, G., \emph{Handbook of procurement}, Cambridge University Press, Cambridge, UK, 2006.
\newblock \urlprefix\url{https://books.google.com/books?hl=en&lr=&id=K3ChhBdIMxwC&oi=fnd&pg=PR5&dq=Handbook+of+procurement&ots=ah3zeKUjnH&sig=R5bnC3lQDaa_ea-HrDkkK8Ee-Qw#v=onepage&q=Handbook%20of%20procurement&f=false}.

\bibitem[{B{\"u}y{\"u}kda{\u{g}} et~al.(2020)B{\"u}y{\"u}kda{\u{g}}, Soysal, and Kitapci}]{buyukdaug2020effect}
B{\"u}y{\"u}kda{\u{g}}, N., Soysal, A.~N., and Kitapci, O., \enquote{The effect of specific discount pattern in terms of price promotions on perceived price attractiveness and purchase intention: An experimental research,} \emph{Journal of Retailing and Consumer Services}, Vol.~55, 2020, p. 102112.
\newblock \doi{https://doi.org/10.1016/j.jretconser.2020.102112}.

\bibitem[{Alford and Biswas(2002)}]{alford2002effects}
Alford, B.~L., and Biswas, A., \enquote{The effects of discount level, price consciousness and sale proneness on consumers' price perception and behavioral intention,} \emph{Journal of Business research}, Vol.~55, No.~9, 2002, pp. 775--783.
\newblock \doi{https://doi.org/10.1016/S0148-2963(00)00214-9}.

\bibitem[{Adenso-D{\'\i}az and Lozano(2023)}]{adenso2023metafrontier}
Adenso-D{\'\i}az, B., and Lozano, S., \enquote{A metafrontier analysis approach for assessing the efficiency of freight service providers,} \emph{International Journal of Systems Science: Operations \& Logistics}, Vol.~10, No.~1, 2023, p. 2177896.
\newblock \doi{https://doi.org/10.1080/23302674.2023.2177896}.

\bibitem[{{Federal Motor Carrier Safety Administration}(accessed January 21, 2024)}]{hours1}
{Federal Motor Carrier Safety Administration}, \enquote{Hours of Service (HOS),} \url{https://www.fmcsa.dot.gov/regulations/hours-of-service}, accessed January 21, 2024.

\bibitem[{Mohammadi et~al.(2014)Mohammadi, Soleymani, and Mozafari}]{mohammadi2014scenario}
Mohammadi, S., Soleymani, S., and Mozafari, B., \enquote{Scenario-based stochastic operation management of microgrid including wind, photovoltaic, micro-turbine, fuel cell and energy storage devices,} \emph{International Journal of Electrical Power \& Energy Systems}, Vol.~54, 2014, pp. 525--535.
\newblock \doi{https://doi.org/10.1016/j.ijepes.2013.08.004}.

\bibitem[{Bornapour et~al.(2019)Bornapour, Hooshmand, and Parastegari}]{bornapour2019efficient}
Bornapour, M., Hooshmand, R.-A., and Parastegari, M., \enquote{An efficient scenario-based stochastic programming method for optimal scheduling of CHP-PEMFC, WT, PV and hydrogen storage units in micro grids,} \emph{Renewable energy}, Vol. 130, 2019, pp. 1049--1066.
\newblock \doi{https://doi.org/10.1016/j.renene.2018.06.113}.

\bibitem[{Ahmadi et~al.(2016)Ahmadi, Charwand, Siano, Nezhad, Sarno, Gitizadeh, and Raeisi}]{ahmadi2016novel}
Ahmadi, A., Charwand, M., Siano, P., Nezhad, A.~E., Sarno, D., Gitizadeh, M., and Raeisi, F., \enquote{A novel two-stage stochastic programming model for uncertainty characterization in short-term optimal strategy for a distribution company,} \emph{Energy}, Vol. 117, 2016, pp. 1--9.
\newblock \doi{https://doi.org/10.1016/j.energy.2016.10.067}.

\bibitem[{Ahmadi et~al.(2013)Ahmadi, Charwand, and Aghaei}]{ahmadi2013risk}
Ahmadi, A., Charwand, M., and Aghaei, J., \enquote{Risk-constrained optimal strategy for retailer forward contract portfolio,} \emph{International Journal of Electrical Power \& Energy Systems}, Vol.~53, 2013, pp. 704--713.
\newblock \doi{https://doi.org/10.1016/j.ijepes.2013.05.051}.

\bibitem[{Kumar et~al.(2021{\natexlab{a}})Kumar, Singh, Bilga, Singh, Singh, Scutaru, Pruncu et~al.}]{EWM}
Kumar, R., Singh, S., Bilga, P.~S., Singh, J., Singh, S., Scutaru, M.-L., Pruncu, C.~I., et~al., \enquote{Revealing the benefits of entropy weights method for multi-objective optimization in machining operations: A critical review,} \emph{Journal of materials research and technology}, Vol.~10, 2021{\natexlab{a}}, pp. 1471--1492.
\newblock \doi{https://doi.org/10.1016/j.jmrt.2020.12.114}.

\bibitem[{Kumar et~al.(2021{\natexlab{b}})Kumar, Dubey, Singh, Singh, Prakash, Nirsanametla, Kr{\'o}lczyk, and Chudy}]{EWM1}
Kumar, R., Dubey, R., Singh, S., Singh, S., Prakash, C., Nirsanametla, Y., Kr{\'o}lczyk, G., and Chudy, R., \enquote{Multiple-criteria decision-making and sensitivity analysis for selection of materials for knee implant femoral component,} \emph{Materials}, Vol.~14, No.~8, 2021{\natexlab{b}}, p. 2084.
\newblock \doi{https://doi.org/10.3390/ma14082084}.

\bibitem[{Chodha et~al.(2022)Chodha, Dubey, Kumar, Singh, and Kaur}]{EWM2}
Chodha, V., Dubey, R., Kumar, R., Singh, S., and Kaur, S., \enquote{Selection of industrial arc welding robot with TOPSIS and Entropy MCDM techniques,} \emph{Materials Today: Proceedings}, Vol.~50, 2022, pp. 709--715.
\newblock \doi{https://doi.org/10.1016/j.matpr.2021.04.487}.

\bibitem[{Borjalilu et~al.(2021)Borjalilu, Sazvar, and Nayeri}]{EWM3}
Borjalilu, N., Sazvar, Z., and Nayeri, S., \enquote{An integrated method for airline company supplier selection based on the entropy and vikor methods: a real case study,} \emph{International Journal of Aviation, Aeronautics, and Aerospace}, Vol.~8, No.~4, 2021, p.~1.
\newblock \doi{https://doi.org/10.15394/ijaaa.2021.1626}.

\bibitem[{Boran et~al.(2009)Boran, Gen{\c{c}}, Kurt, and Akay}]{boran2009multi}
Boran, F.~E., Gen{\c{c}}, S., Kurt, M., and Akay, D., \enquote{A multi-criteria intuitionistic fuzzy group decision making for supplier selection with TOPSIS method,} \emph{Expert systems with applications}, Vol.~36, No.~8, 2009, pp. 11363--11368.
\newblock \doi{https://doi.org/10.1016/j.eswa.2009.03.039}.

\bibitem[{Nasir et~al.(2022)Nasir, Kansal, Alshaltone, Barneih, Sameer, Shanableh, and Al-Shamma'a}]{w1}
Nasir, N., Kansal, A., Alshaltone, O., Barneih, F., Sameer, M., Shanableh, A., and Al-Shamma'a, A., \enquote{Water quality classification using machine learning algorithms,} \emph{Journal of Water Process Engineering}, Vol.~48, 2022, p. 102920.
\newblock \doi{https://doi.org/10.1016/j.jwpe.2022.102920}.

\bibitem[{Sun et~al.(2022)Sun, Min, Lu, and Zhai}]{w}
Sun, W., Min, X., Lu, W., and Zhai, G., \enquote{A deep learning based no-reference quality assessment model for ugc videos,} \emph{Proceedings of the 30th ACM International Conference on Multimedia}, 2022, pp. 856--865.
\newblock \doi{https://doi.org/10.1145/3503161.3548329}.

\bibitem[{You and Grossmann(2013)}]{you2013multicut}
You, F., and Grossmann, I.~E., \enquote{Multicut Benders decomposition algorithm for process supply chain planning under uncertainty,} \emph{Annals of Operations Research}, Vol. 210, 2013, pp. 191--211.
\newblock \doi{https://doi.org/10.1007/s10479-011-0974-4}.

\bibitem[{Laporte and Louveaux(1993)}]{laporte1993integer}
Laporte, G., and Louveaux, F.~V., \enquote{The integer L-shaped method for stochastic integer programs with complete recourse,} \emph{Operations research letters}, Vol.~13, No.~3, 1993, pp. 133--142.
\newblock \doi{https://doi.org/10.1016/0167-6377(93)90002-X}.

\bibitem[{Benders(2005)}]{benders2005partitioning}
Benders, J.~F., \enquote{Partitioning procedures for solving mixed-variables programming problems,} \emph{Computational Management Science}, Vol.~2, No.~1, 2005, pp. 3--19.
\newblock \doi{https://doi.org/10.1007/s10287-004-0020-y}.

\bibitem[{Akash et~al.(2021)Akash, Raj, Sushmitha, Prateek, Aditya, and Sreehari}]{akash2021design}
Akash, A., Raj, V. S.~J., Sushmitha, R., Prateek, B., Aditya, S., and Sreehari, V.~M., \enquote{Design and analysis of VTOL operated intercity electrical vehicle for urban air mobility,} \emph{Electronics}, Vol.~11, No.~1, 2021, p.~20.
\newblock \doi{https://doi.org/10.3390/electronics11010020}.

\bibitem[{Erden and Yayla(2021)}]{seat}
Erden, S., and Yayla, P., \enquote{Finite Element Stress Analysis of Airplane Seat,} \emph{European Mechanical Science}, Vol.~5, No.~1, 2021, pp. 6--13.
\newblock \doi{https://doi.org/10.26701/ems.799180}.

\bibitem[{{Spin Life}(accessed 15 August, 2023)}]{seat1}
{Spin Life}, \enquote{Spirit Car Seat,} \url{https://www.spinlife.com/Inspired-by-Drive-Spirit-Car-Seat-Car-Seats-and-Boosters/spec.cfm?productID=176234&adv=bingshopping&utm_medium=CSE&utm_source=bingShopping}, accessed 15 August, 2023.

\bibitem[{Cachon and Fisher(2000)}]{cachon2000supply}
Cachon, G.~P., and Fisher, M., \enquote{Supply chain inventory management and the value of shared information,} \emph{Management science}, Vol.~46, No.~8, 2000, pp. 1032--1048.
\newblock \doi{https://doi.org/10.1287/mnsc.46.8.1032.12029}.

\bibitem[{Hess and Romanoff(1987)}]{airframe}
Hess, R.~W., and Romanoff, H., \enquote{Aircraft Airframe Cost Estimating Relationships: Study Approach and Conclusions,} Tech. Rep. R-3255-AF, Rand Corporation, 1987.
\newblock \urlprefix\url{https://www.rand.org/pubs/reports/R3255.html}.

\bibitem[{Alba-Maestre et~al.(2021)Alba-Maestre, Beyne, Buszek, Cuadrat-G, López, Poliakov, Reine, Santamaría, Schoser, Wadia, and Castro}]{weigon}
Alba-Maestre, J., Beyne, E., Buszek, M., Cuadrat-G, M., López, N., Poliakov, N., Reine, K., Santamaría, A., Schoser, J., Wadia, K., and Castro, S., \enquote{Final Report - Multi-Disciplinary Design and Optimisation of a Long-Range eVTOL Aircraft,} Tech. rep., 11 2021.
\newblock \doi{10.5281/zenodo.5576103}.

\bibitem[{{Marketplace - Aviation Week}(accessed July 7, 2023)}]{p_e_1}
{Marketplace - Aviation Week}, \enquote{Starr Aircraft Products, Inc.} \url{https://marketplace.aviationweek.com/company/starr-aircraft-products-inc}, accessed July 7, 2023.

\bibitem[{{Statista}(accessed August 18, 2023)}]{market1}
{Statista}, \enquote{Projected global urban air mobility market in 2023 and 2035,} \url{https://www.statista.com/statistics/1046436/global-urban-air-mobility-market-size/}, accessed August 18, 2023.

\bibitem[{{Research and Markets}(accessed February 2024)}]{market3}
{Research and Markets}, \enquote{United States Urban Air Mobility (UAM) Market Report 2021,} \url{https://www.globenewswire.com/news-release/2021/10/26/2320467/28124/en/United-States-Urban-Air-Mobility-UAM-Market-Report-2021-Market-is-Estimated-to-Reach-18-81-Billion-in-2035-at-a-CAGR-of-23-12-during-2023-2035.html}, accessed February 2024.

\bibitem[{{Lilium}(accessed August 17, 2023)}]{lilium}
{Lilium}, \enquote{SEC Filings,} \url{https://investors.lilium.com/financial-information/sec-filings}, accessed August 17, 2023.

\bibitem[{{Investor Relations}(accessed August 17, 2023)}]{joby}
{Investor Relations}, \enquote{SEC Filings,} \url{https://ir.jobyaviation.com/sec-filings}, accessed August 17, 2023.

\bibitem[{{CNBC}(accessed August 17, 2023)}]{sell1}
{CNBC}, \enquote{Vertical Aerospace to go public in \$2.2 billion SPAC, gets American Airlines investment,} \url{https://www.cnbc.com/2021/06/10/vertical-aerospace-to-go-public-in-2point2-billion-spac-gets-american-airlines-investment.html}, accessed August 17, 2023.

\bibitem[{{Reuters}(accessed August 17, 2023)}]{sell2}
{Reuters}, \enquote{Vertical Aerospace to go public in \$2.2 billion SPAC deal,} \url{https://www.reuters.com/business/aerospace-defense/american-airlines-invest-electric-aircraft-maker-vertical-aerospace-2021-06-10/}, accessed August 17, 2023.

\bibitem[{{Kirsten Korosec}(accessed February 10, 2023)}]{sell3}
{Kirsten Korosec}, \enquote{Archer lands \$1B order from United Airlines and a SPAC deal,} \url{https://techcrunch.com/2021/02/10/archer-lands-1-1b-order-from-united-airlines-and-a-spac-deal/}, accessed February 10, 2023.

\bibitem[{Odedairo et~al.(2020)Odedairo, Alaba, and Edem}]{odedairo2020system}
Odedairo, B.~O., Alaba, E.~H., and Edem, I., \enquote{A System Dynamics Model To determine the value of inventory holding cost,} \emph{Journal of Engineering Studies and Research}, Vol.~26, No.~3, 2020, pp. 112--123.
\newblock \doi{https://doi.org/10.29081/jesr.v26i3.213}.

\bibitem[{Azzi et~al.(2014)Azzi, Battini, Faccio, Persona, and Sgarbossa}]{azzi2014inventory}
Azzi, A., Battini, D., Faccio, M., Persona, A., and Sgarbossa, F., \enquote{Inventory holding costs measurement: a multi-case study,} \emph{The International Journal of Logistics Management}, Vol.~25, No.~1, 2014, pp. 109--132.
\newblock \doi{https://doi.org/10.1108/IJLM-01-2012-0004}.

\bibitem[{{Eco TransIT World}(accessed August 15, 2023)}]{eco-calculator}
{Eco TransIT World}, \enquote{Emission calculator for greenhouse gases and exhaust emissions,} \url{https://www.ecotransit.org/en/emissioncalculator/}, accessed August 15, 2023.

\bibitem[{{Bureau of Transportation Statistics}(accessed August 17, 2023)}]{allowed_speed}
{Bureau of Transportation Statistics}, \enquote{Average Freight Revenue per Ton-Mile,} \url{https://www.bts.gov/content/average-freight-revenue-ton-mile}, accessed August 17, 2023.

\bibitem[{{Bureau of Transportation Statistics}(accessed December 2023)}]{unit_cost_mode}
{Bureau of Transportation Statistics}, \enquote{Freight Transportation System Extent \& Use,} \url{https://data.bts.gov/stories/s/Freight-Transportation-System-Extent-Use/r3vy-npqd/}, accessed December 2023.

\bibitem[{{Carbon Offset Guide}(accessed August 17, 2023)}]{emission_cost}
{Carbon Offset Guide}, \enquote{Understanding carbon offsets,} \url{https://www.offsetguide.org/understanding-carbon-offsets/other-instruments-for-claiming-emission-reductions/allowances/#:~:text=California:%20US%2414-18%20EU,ETS:%20around%20%E2%82%AC20-30}, accessed August 17, 2023.

\bibitem[{{California Air Resources Board}(accessed August 17, 2023)}]{emission_cost1}
{California Air Resources Board}, \enquote{Cap-and-Trade Program - Cost Containment Informations,} \url{https://ww2.arb.ca.gov/our-work/programs/cap-and-trade-program/cost-containment-information#:~:text=Pursuant%20to%20section%2095913%20%28h%29%20%285%29%20of%20the,since%202021%20are%20provided%20in%20the%20table%20below.}, accessed August 17, 2023.

\bibitem[{Timmis(2020)}]{timmis2020aircraft}
Timmis, A., \emph{Aircraft Manufacturing and Technology}, 2\textsuperscript{nd} ed., Routledge, 2020.
\newblock \urlprefix\url{https://www.taylorfrancis.com/chapters/edit/10.4324/9780429299445-17/aircraft-manufacturing-technology-andrew-timmis}.

\bibitem[{Wada(2020)}]{wada2020evolution}
Wada, K., \emph{The Evolution of the Toyota Production System}, Springer, Gateway East, Singapore, 2020.
\newblock \urlprefix\url{https://link.springer.com/content/pdf/10.1007/978-981-15-4928-1.pdf}.

\bibitem[{Cohen et~al.(2021)Cohen, Shaheen, and Farrar}]{cohen2021urban}
Cohen, A.~P., Shaheen, S.~A., and Farrar, E.~M., \enquote{Urban air mobility: History, ecosystem, market potential, and challenges,} \emph{IEEE Transactions on Intelligent Transportation Systems}, Vol.~22, No.~9, 2021, pp. 6074--6087.
\newblock \doi{https://doi.org/10.1109/TITS.2021.3082767}.

\bibitem[{Garrow et~al.(2019)Garrow, German, Mokhtarian, and Glodek}]{garrow2019survey}
Garrow, L.~A., German, B., Mokhtarian, P., and Glodek, J., \enquote{A survey to model demand for eVTOL urban air trips and competition with autonomous ground vehicles,} \emph{AIAA Aviation 2019 Forum}, 2019, p. 2871.
\newblock \doi{https://doi.org/10.2514/6.2019-2871}.

\end{thebibliography}
\end{spacing}

\end{document}